\documentclass{article}
\usepackage{latexsym,amsmath,amssymb,amscd}
\usepackage[mathscr]{eucal}
\usepackage[dvips]{color}
\begin{document}

\title{On Linear Adiabatic Perturbations of Spherically Symmetric Gaseous Stars Governed by the Euler-Poisson Equations  }
\author{Tetu Makino \footnote{Professor Emeritus at Yamaguchi University, Japan  E-mail: makino@yamaguchi-u.ac.jp}}
\date{June 21, 2021}
\maketitle

\newtheorem{Lemma}{Lemma}
\newtheorem{Proposition}{Proposition}
\newtheorem{Theorem}{Theorem}
\newtheorem{Definition}{Definition}
\newtheorem{Remark}{Remark}
\newtheorem{Corollary}{Corollary}
\newtheorem{Assumption}{Assumption}
\newtheorem{Notation}{Notation}

\numberwithin{equation}{section}

\begin{abstract}
The linearized operator for non-radial oscillations of spherically symmetric self-gravitating gaseous stars is analyzed in view of the functional analysis. The evolution of the star is supposed to be governed by the Euler-Poisson equations under the equation of state of the ideal gas, and the motion is supposed to be adiabatic.  We consider the case of not necessarily isentropic, that is, not barotropic motions. Basic theory of self-adjoint realization of the linearized operator is established. Some problems in the investigation of the concrete properties of the spectrum of the linearized operator are proposed. The existence of eigenvalues which accumulate to 0 is proved in a mathematically rigorous fashion.
The absence of continuous spectra and the completeness of eigenfunctions for the operators reduced by spherical harmonics is discussed. \\

Key Words and Phrases. Gaseous star. Adiabatic Oscillation. Self-adjoint operator. Friedrichs extension. Spectrum of Sturm-Liouville type. Brunt-V\"{a}is\"{a}l\"{a} frequency. Gravity mode.

2010 Mathematical Subject Classification Numbers. 35P05, 35L51, 35Q31, 35Q85, 46N20, 76N15.
\end{abstract}

\section{Introduction}

We consider the adiabatic hydrodynamic evolution of a self-gravitating gaseous star governed by the Euler-Poisson equations 
\begin{subequations}
\begin{align}
&\frac{\partial\rho}{\partial t}+\sum_{k=1}^3\frac{\partial}{\partial x^k}(\rho v^k)=0, \label{EPa} \\
&\rho\Big(\frac{\partial v^j}{\partial t}+\sum_{k=1}^3 v^k\frac{\partial v^j}{\partial x^k}\Big)
+\frac{\partial P}{\partial x^j}+\rho\frac{\partial\Phi}{\partial x^j}=0, \quad j=1,2,3,
 \label{EPb} \\
&\rho\Big(\frac{\partial S}{\partial t}+\sum_{k=1}^3 v^k\frac{\partial S}{\partial x^k}\Big)=0, \label{EPc} \\
&\triangle \Phi =4\pi\mathsf{G}\rho. \label{EPd}
\end{align}
\end{subequations} 
Here $t\geq 0, \mbox{\boldmath$x$}=(x^1,x^2,x^3) \in \mathbb{R}^3$. The unknowns $\rho \geq 0, P, S, \Phi \in \mathbb{R}$ are the density, the pressure, the specific entropy, the gravitational potential, and $\mbox{\boldmath$v$}=(v^1,v^2,v^3) \in \mathbb{R}^3$ is the velocity fields. $\mathsf{G}$ is a positive constant, the gravitation constant. 

In this article the pressure $P$ is supposed to be a prescribed function of $\rho, S$. But for the sake of simplicity, we assume the equation of state of the ideal fluid, that is, we assume

\begin{Assumption} $P$ is the function of $(\rho, S) \in [0,+\infty[ \times \mathbb{R}$ given by
\begin{equation}
P=\rho^{\gamma}\exp\Big(\frac{S}{\mathsf{C}_V}\Big), \label{DefP}
\end{equation}
where $\gamma$ and $\mathsf{C}_V$ are positive constants such that
\begin{equation}
1<\gamma <2.
\end{equation}
\end{Assumption}

The constant $\gamma$ is the adiabatic exponent and $\mathsf{C}_V$ is the specific heat per unit mass at constant volume. 

Since we are concerned with compactly supported density distribution $\rho$ in this article, the Poisson equation
\eqref{EPd} will be replaced by the Newtonian potential
\begin{equation}
\Phi(t,\mbox{\boldmath$x$})=-4\pi\mathsf{G}\mathcal{K}\rho(t,\cdot)(\mbox{\boldmath$x$}), \label{NPot}
\end{equation}
where
\begin{equation}
\mathcal{K}f(\mbox{\boldmath$x$}):=\frac{1}{4\pi}\int\frac{f(\mbox{\boldmath$x$}')}{|\mbox{\boldmath$x$}-\mbox{\boldmath$x$}'|}d\mbox{\boldmath$x$}'. \label{defK}
\end{equation}

We suppose that there is fixed a spherically symmetric equilibrium $\bar{\rho}, \bar{P},
\bar{S}, \bar{\Phi}$, which satisfy \eqref{EPa}, \eqref{EPb}, \eqref{EPc}, \eqref{NPot},
such that $\bar{\rho}(\mbox{\boldmath$x$})>0 \Leftrightarrow r=|\mbox{\boldmath$x$}| <R$ with a finite positive number $R$, the radius of the equilibrium.

We consider the perturbation $\mbox{\boldmath$\xi$}=\delta\mbox{\boldmath$x$}, \delta\rho, \delta P, \delta S, \delta\Phi$ at this fixed equilibrium. We use the Lagrangian co-ordinate which will be dented by the diversion of the letter $\mbox{\boldmath$x$}$ of the Eulerian co-ordinate. So, $\mbox{\boldmath$x$}$ runs on the fixed domain
$B_R:=\{ \mbox{\boldmath$x$} \in \mathbb{R}^3 | r=|\mbox{\boldmath$x$}| < R\}$, while $\{\rho >0\} $ described by the Eulerian co-ordinate may move along $t$. 

Then the linearized equation which governs the perturbations turns out to be
\begin{equation}
\frac{\partial^2\mbox{\boldmath$\xi$}}{\partial t^2}+\mbox{\boldmath$L$}\mbox{\boldmath$\xi$}=0,
\label{Eqxi}
\end{equation}
where
\begin{align}
\mbox{\boldmath$L$}\mbox{\boldmath$\xi$}&=\frac{1}{\bar{\rho}}\mathrm{grad}\delta P+\frac{\delta\rho}{\bar{\rho}}\mathrm{grad}\bar{\Phi}+
\mathrm{grad}\delta\Phi \nonumber \\
&=\frac{1}{\bar{\rho}}\mathrm{grad}\delta P
-\frac{\delta\rho}{\bar{\rho}^2}\mathrm{grad}\bar{P}+
\mathrm{grad}\delta\Phi. \label{L}
\end{align}

We have
\begin{align}
&\delta\rho=-\mathrm{div}(\bar{\rho}\mbox{\boldmath$\xi$}), \label{drho} \\
&\delta\Phi=-4\pi\mathsf{G}\mathcal{K}(\delta\rho). \label{dPhi}
\end{align}

Here we use the following notation:

\begin{Notation}

The symbol  $\delta$ denotes the Eulerian perturbation, while $\Delta$ will denote the Lagrangian perturbation, which are defined by
\begin{align*}
&\Delta Q(t,\mbox{\boldmath$x$})=Q(t, \mbox{\boldmath$\varphi$}(t,\mbox{\boldmath$x$}))-\bar{Q}(\mbox{\boldmath$x$}), \\
&\delta Q(t,\mbox{\boldmath$x$})=Q(t,\mbox{\boldmath$\varphi$}(t,\mbox{\boldmath$x$}))-\bar{Q}
(\mbox{\boldmath$\varphi$}(t,\mbox{\boldmath$x$})),
\end{align*}
where $\mbox{\boldmath$\varphi$}(t,\mbox{\boldmath$x$})=\mbox{\boldmath$x$}+\mbox{\boldmath$\xi$}(t,\mbox{\boldmath$x$})$ is the steam line given by 
$$\frac{\partial}{\partial t}\mbox{\boldmath$\varphi$}(t,\mbox{\boldmath$x$})=\mbox{\boldmath$v$}(t,\mbox{\boldmath$\varphi$}(t,\mbox{\boldmath$x$})),
\quad \mbox{\boldmath$\varphi$}(0,\mbox{\boldmath$x$})=\mbox{\boldmath$x$}. $$
\end{Notation}

Note that
$$\Delta Q=\delta Q+(\mbox{\boldmath$\xi$}|\mathrm{grad}\bar{Q}) $$
in the linearized approximation, for any quantity $Q$.\\

Supposing that the initial perturbation of the density vanishes, that is, $\Delta\rho|_{t=0}=0$, the equation \eqref{EPa} implies $\Delta\rho+\bar{\rho}\mathrm{div}\mbox{\boldmath$\xi$}=0$ always, which is \eqref{drho}.
Supposing $\Delta S|_{t=0} =0$, the equation \eqref{EPc} implies $\Delta S=0$ always, therefore
$$\Delta P=\overline{\frac{\gamma P}{\rho}}\Delta\rho. $$
This implies
\begin{equation}
\delta P=\overline{\frac{\gamma P}{\rho}}\delta\rho+
\gamma\mathscr{A}\bar{P}(\mbox{\boldmath$\xi$}|\mbox{\boldmath$e$}_r) \label{dP}.
\end{equation}\\

Here we denote by  $\mbox{\boldmath$e$}_r$ the unit vector $\partial/\partial r$ and we  define the `Schwarzschild's discriminant of convective stability' $\mathscr{A}$ with the `Brunt-V\"{a}is\"{a}l\"{a} frequency' $\mathscr{N}$ by

\begin{Definition}
We put
\begin{equation}
\mathscr{A}:=\overline{\frac{1}{\rho}\frac{d\rho}{dr}}-\overline{\frac{1}{\gamma P}\frac{dP}{dr}}\quad \Big(=-\frac{1}{\gamma \mathsf{C}_V}\frac{d\bar{S}}{dr}\Big) 
\end{equation}
and
\begin{equation}
\mathscr{N}^2:=\mathscr{A}\overline{\frac{1}{\rho}\frac{dP}{dr}}
=-\mathscr{A}\overline{\frac{d\Phi}{dr}}.
\end{equation}
\end{Definition}

For the physical meaning of these quantities, see \cite{LedouxW} or \cite[Chapter III, Section 17]{Cox}. When $\overline{dP/dr}<0$, $\mathscr{N}$ is real if and only if
$\mathscr{A} \leq 0$. The condition $\mathscr{A}<0$ is that of the convective stability.\\

Note that $\Delta S=0$ for $\forall t \geq 0$ means $S(t, \mbox{\boldmath$\varphi$}(t,\mbox{\boldmath$x$}))=\bar{S}(\mbox{\boldmath$x$})$ for $\forall t, \forall \mbox{\boldmath$x$}$. Then, under the linearized approximation, we have
\begin{align*}
\delta S&=\Delta S-(\mbox{\boldmath$\xi$}|\mathrm{grad}\bar{S}) = -(\mbox{\boldmath$\xi$}|\mathrm{grad}\bar{S}) = \\
&=-(\mbox{\boldmath$\xi$}|\mbox{\boldmath$e$}_r)\frac{d\bar{S}}{dr}=\gamma\mathsf{C}_V\mathscr{A}(\mbox{\boldmath$\xi$}|\mbox{\boldmath$e$}_r).
\end{align*}\\

Thus by \eqref{drho}, \eqref{dPhi}, \eqref{dP} we can see the right-hand side of \eqref{L} 
is an integro-differential operator acting on the unknown $\mbox{\boldmath$\xi$}$, provided that the spherically symmetric equilibrium
$(\bar{\rho}, \bar{S}, \bar{P})$ is fixed.

For the derivation of $\mbox{\boldmath$L$}$, see e.g., \cite{LedouxW}, \cite{Cox} or
\cite{Unno}.\\

The purpose of this article is to clarify the functional analysis properties of this
integro-differential operator $\mbox{\boldmath$L$}$.\\

Nonlinear evolution of spherically symmetric perturbations has been investigated sufficiently well in \cite{TM.OJM} and \cite{JJ.APDE}. In these studies spectral properties of the linearized operator for spherically symmetric perturbations, which was established by \cite{Beyer1995} and independently by \cite{Lin}, are fully presupposed.
Its spectrum was proved to be actually of the Sturm-Liouville type, and it was not obvious because of the singularity of the coefficients, caused by the physical vacuum boundary of the equilibrium.  Therefore if we want to study nonlinear evolution of not necessarily spherically symmetric perturbations around a spherically symmetric equilibrium, we should prepare a sufficiently strong functional analysis study of spectral properties of the linearized operator for general, not necessarily spherically  symmetric, perturbations. As for barotropic case, we have attacked this task, and have gotten sufficiently strong results in \cite{JJTM}. Thus here we consider the case of
not necessarily barotropic motions. Unfortunately the results which we have established is little bit weaker than the barotropic case. There remains some open problems. But mathematically rigorous treatment of the problem is quiet new.\\

This article is organized as follows.

In Section 2, we discuss on the existence of spherically symmetric equilibrium for prescribed entropy distribution. The concept of the `admissible' equilibrium will play a crucial r\^{o}le throughout the mathematically rigorous investigations of this article. In Section 3, we prove the self-adjoint realization of the operator $\mbox{\boldmath$L$}$ as the Friedrichs extension in the Hilbert space $\mathfrak{H}=L^2(\bar{\rho}d\mbox{\boldmath$x$})$. Astrophysical texts lacked mathematically rigorous proof. But such a  strong assertion on the concrete  form of the spectrum as that of the barotropic case given in \cite{JJTM} is not yet obtained. In order to investigate the specified concrete form of the spectrum we investigate eigenfunctions represented by spherical harmonics in Section 4. The situation is clarified to be quite different from the barotropic case. That is, it may be impossible to reduce the problem to that of Sturm-Liouville type.  But the justification of the self-adjoint
realization of the associated operator $\vec{L}_l$ for each degree $l$ of the harmonics $Y_{lm}$ 
in the Hilbert space $\mathfrak{X}_l=L^2(\bar{\rho}r^2dr)$ can be done with success. A strong guess that the form of the spectrum of
$\vec{L}_l$ is quite different from that of the barotropic case is suggested by the so-called `g-modes', say, a
sequence of eigenvalues accumulating to $0$. Section 5 is devoted to a mathematically rigorous proof of the existence of
 the g-modes  proposed by astrophysicists for the so called Cowling approximation given by neglecting the perturbation of the gravitational potential.
 We shall give a proof of the existence of g-modes under the assumption that 
$\inf \mathscr{N}^2 /r^2>0$
and a set of restrictions on the configuration of the back ground equilibrium which guarantees the smallness of the effect of the self-gravitational perturbation. Also, the existence of `p-modes', say, the existence of eigenvalues which accumulate to $+\infty$ will be proved supposing neither $\inf\mathscr{N}^2 /r^2>0$
nor other restrictions. 
 The last Section 6 is devoted to examination of the arguments in the work \cite{Eisenfeld} by J. Eisenfeld.
Besides the assumption for 
 $1/(\gamma-1) $ to be an integer done without reasoning,  it seems that the proof of the existence and completeness of eigenvalues found in \cite{Eisenfeld} is not so complete. 
 Therefore 
we try to give a rigorous  proof of the absence of continuous spectrum
of the self-adjoint operator $\vec{L}_l$. This is done by considering the operator not in $\mathfrak{X}_l$ but in the subspace $\mathfrak{W}_l$, which is dense in $\mathfrak{X}_l$. \\

We shall use the following notations:

\begin{Notation}
We denote
\begin{subequations}
\begin{align}
&B_R:=\{ \mbox{\boldmath$x$}\in\mathbb{R}^3 \  |\  r=|\mbox{\boldmath$x$}| <R \}, \\
& \overline{B_R}:=\{ \mbox{\boldmath$x$}\in\mathbb{R}^3 \  |\  r=|\mbox{\boldmath$x$}| \leq R \}.
\end{align}
\end{subequations}
\end{Notation}

\begin{Notation}\label{Not.1}
1) A function $F$ on a subset of $\mathbb{R}^3$ is said to be spherically symmetric if there exists a function $f$ on a subset of $[0,+\infty[$ such that
$F(\mbox{\boldmath$x$})=f(|\mbox{\boldmath$x$}|)$ for $\forall \mbox{\boldmath$x$}$ in the domain of $F$. Then we shall denote $f=F^{\sharp}$.

2) For a function $f$ on a subset of $[0,+\infty[$, we shall denote by $f^{\flat}$ the function on a subset of $\mathbb{R}^3$ such that $f^{\flat}(\mbox{\boldmath$x$})=
f(|\mbox{\boldmath$x$}|)$ for $\forall \mbox{\boldmath$x$}$ such that $r=|\mbox{\boldmath$x$}|$ is in the domain of $f$.

3) When it is expected that no confusion may occur, we shall divert the symbols $f$ or $F$
instead of $f^{\flat}$ or $F^{\sharp}$.

\end{Notation}

Here let us note the following lemma, which can be verified easily:

\begin{Lemma}\label{Lem.0}
If a function $f$ defined on $[0,R[$ satisfies $f^{\flat}\in C^{k+2}(B_R), k=0,1,2,\cdots$, then 
$f \in C^{k+2}([0,R[), \quad \displaystyle \frac{df}{dr}\Big|_{r=+0}=0$, \  
and $\displaystyle \frac{1}{r}\frac{df}{dr} \in C^k([0,R[)$.
\end{Lemma}

Proof. We can show inductively that
$$D^kh(r)=\frac{1}{r^{k+1}}\int_0^rD^{k+2}f(s)s^kds $$
for $\displaystyle h=\frac{1}{r}\frac{df}{dr}$. Therefore
$$D^kh(r) \rightarrow \frac{D^{k+2}f(0)}{k+1}\quad\mbox{as}\quad r \rightarrow +0.$$
$\square$\\

\begin{Notation}
We denote the unit vectors
\begin{equation}
\mbox{\boldmath$e$}_r=\frac{\partial}{\partial r},\quad \mbox{\boldmath$e$}_{\vartheta}=\frac{1}{r}\frac{\partial}{\partial\vartheta},
\quad \mbox{\boldmath$e$}_{\phi}=\frac{1}{r\sin\phi}\frac{\partial}{\partial\phi}
\end{equation}
for the spherical polar co-ordinates
\begin{align}
&x^1=r\sin\vartheta\cos\phi, \nonumber \\
&x^2=r\sin\vartheta\sin\phi, \nonumber \\
&x^3=r\cos\vartheta.
\end{align}
\end{Notation}

\section{Existence of spherically symmetric equilibrium for prescribed entropy distribution}
 
In this section we establish the existence of spherically symmetric equilibria which enjoy good properties used in the following consideration on $\mbox{\boldmath$L$}$.\\

Let us put the following 

\begin{Definition}\label{Def.2}
A pair of $t$-independent spherically symmetric functions $(\bar{\rho}, \bar{S})
\in C_0^1(\mathbb{R}^3; [0,+\infty[)\times C^1(\mathbb{R}^3; \mathbb{R})$
which satisfies \eqref{EPa}\eqref{EPb}\eqref{EPc} with $\mbox{\boldmath$v$}=0$ and $\Phi$, $P$ determined by \eqref{NPot},
 \eqref{DefP}  is called an admissible spherically symmetric equilibrium, if there is a finite positive number $R$ such that

1) $\{\bar{\rho} >0\}=B_R $;

2) $\bar{\rho}^{\gamma-1}, \bar{S}, \bar{\Phi}=-4\pi\mathsf{G}\mathcal{K}\bar{\rho} \in C^{\infty}(B_R) \cap
C^{3,\alpha}(\overline{B_R})$, $\alpha$ being a positive number such that
$0<\alpha <\min(\frac{1}{\gamma-1}-1, 1)$ ;

3) $\displaystyle \frac{d\bar{\rho}^{\sharp}}{dr}, \frac{dP^{\sharp}}{dr} <0$ for $0<r<R$ and
\begin{equation}
\frac{1}{r}\frac{d\bar{\rho}^{\sharp}}{dr}\Big|_{r=+0}<0, 
\quad \frac{1}{r}\frac{dP^{\sharp}}{dr}\Big|_{r=+0} <0; \label{PON}
\end{equation}

4) The boundary $\partial B_R$, on which $\bar{\rho}=0$, is a physical vacuum boundary, that is,
\begin{equation}
-\infty <\frac{d}{dr}(\bar{\rho}^{\gamma-1})^{\sharp}\Big|_{r=R-0}<0, \label{PhysVacBd}
\end{equation}
which means
$$ -\infty < \frac{d}{dr}\overline{\frac{\gamma P}{\rho}}^{\sharp}\Big|_{r=R-0} <0, $$
where $\overline{\gamma P/\rho}=\overline{(\partial P/\partial\rho)_{S=\mbox{Const}}}$ is the square of the sound speed. 
\end{Definition}

Note that $\gamma <2$ implies $1<\frac{1}{\gamma-1}$, therefore such an $\alpha$ exists. Moreover we see that $\bar{\rho}^{\gamma-1}\in C^{3,\alpha}(\overline{B_R})$
implies $\bar{\rho}^{\gamma}, \bar{P}, \bar{S} \in 
C^{2,\alpha}(\overline{B_R})$,
since $\frac{\gamma}{\gamma-1}>2+\alpha$.. However $\bar{\rho}\in C^{1,\alpha}(\overline{B_R})$ and 
$\bar{\rho} \not\in C^2(\overline{B_R})$ unless $\frac{1}{\gamma-1}\geq 2 \Leftrightarrow \gamma \leq \frac{3}{2}$. \\

Note that 4) of Definition \ref{Def.2} is equivalent to 
\begin{equation}
-\infty <\frac{1}{\bar{\rho}^{\sharp}}\frac{d\bar{P}^{\sharp}}{dr}\Big|_{r=R-0} <0.
\end{equation}
In fact, we have
$$\frac{1}{\rho}\frac{dP}{dr}=\Big[\frac{\gamma}{\gamma-1}\frac{d}{dr}\rho^{\gamma-1}+
\frac{\rho^{\gamma-1}}{\mathsf{C}_V}\frac{dS}{dr}\Big]\exp\Big(\frac{S}{\mathsf{C}_V}\Big),$$
and $\rho=\bar{\rho}^{\sharp}\rightarrow 0$ as $r\rightarrow R-0$ and
$S=\bar{S}\in C^1(\bar{B}_R)$.\\

We claim

\begin{Theorem}\label{Th.1}
 Let a smooth function $\Sigma$ on $\mathbb{R}$ and a positive number $\rho_{\mathsf{O}}$ be given. Assume that it holds, for $\omega>0$, that
\begin{equation}
\gamma +\frac{\gamma-1}{\mathsf{C}_V}\omega\frac{d}{d\omega}\Sigma(\omega) >0. \label{pDP}
\end{equation}
Either if $\displaystyle\frac{4}{3}<\gamma <2$ or if $\displaystyle \frac{6}{5}<\gamma \leq \frac{4}{3}$ and $\rho_{\mathsf{O}}$ is sufficiently small, then there exists an
admissible spherically symmetric equilibrium
$(\bar{\rho}, \bar{S})$ such that $\bar{S}=\Sigma(\bar{\rho}^{\gamma-1})$ and 
$\bar{\rho}(O)=\rho_{\mathsf{O}}$.
\end{Theorem}

Proof . Consider the functions $f^P, f^u$ defined by 
\begin{align}
&f^P(\rho):=\rho^{\gamma}\exp\Big[\frac{\Sigma(\rho^{\gamma-1})}{\mathsf{C}_V}\Big], \\
&f^u(\rho):=\int_0^{\rho}
\frac{Df^P(\rho')}{\rho'}d\rho'
\end{align}
for $\rho >0$. Thanks to the assumption \eqref{pDP} we have
$$Df^P(\rho)=
\Big[\gamma+\frac{\gamma-1}{\mathsf{C}_V}\omega\frac{d\Sigma}{d\omega}\Big]_{\omega=\rho^{\gamma-1}}\frac{f^{P}(\rho)}{\rho}
 >0$$
for $\rho >0$, and there exists a smooth function $\Lambda$ on $\mathbb{R}$ such that
$\Lambda(0)=0$ and
\begin{equation}
f^P(\rho)=\mathsf{A}\rho^{\gamma}(1+\Lambda(\rho^{\gamma-1}))
\end{equation}
for $\rho >0$. Here $\mathsf{A}:=\exp(\Sigma(0)/\mathsf{C}_V)$ is a positive constant. 
Then we have
\begin{equation}
u=f^u(\rho)=\frac{\gamma\mathsf{A}}{\gamma-1}\rho^{\gamma-1}(1+\Lambda_u(\rho^{\gamma-1}))
\end{equation}
for $\rho >0$, where $\Lambda_u$ is a smooth function on $\mathbb{R}$ such that 
$\Lambda_u(0)=0$, and  the inverse function $f^{\rho}$ of $f^u$ 
\begin{equation}
f^{\rho}(u)=\Big(\frac{\gamma-1}{\gamma\mathsf{A}}\Big)^{\frac{1}{\gamma-1}}
(u \vee 0)^{\frac{1}{\gamma-1}}(1+\Lambda_{\rho}(u))
\end{equation}
is given so that
$ \rho=f^{\rho}(u) \Leftrightarrow u=f^u(\rho)$ for $u>0 (\rho >0)$. Here $u\vee 0$ stands for $\max( u, 0)$ and  $ \Lambda_{\rho}$ are smooth functions on
$\mathbb{R}$ such that $\Lambda_{\rho}(0)=0$.

Therefore the problem is reduced to that for barotropic case to solve
$$-\frac{1}{r^2}\frac{d}{dr}r^2\frac{du}{dr}=4\pi\mathsf{G}f^{\rho}(u), \quad
u=u_{\mathsf{O}}+O(r^2) \qquad ( r \rightarrow +0 ) \label{Shoot}
$$
by the shooting method.
Here $u_{\mathsf{O}}=f^u(\rho_{\mathsf{O}})$ is given. Then $u$ is a monotone decreasing function of  $r \in ]0,+\infty[$ and the proof of the existence of the finite zero $r=R$  can be found in  Appendix. 
Put $\bar{u}=u^{\flat}$, $u$ being the solution of \eqref{Shoot}, and put $\bar{\rho}=f^{\rho}(\bar{u}),
\bar{S}=\Sigma(\bar{\rho}^{\gamma-1}) $. 
Let us verify the regularity of $\bar{\rho}^{\gamma-1}, \bar{S}, \bar{\Phi}$ including the vacuum boundary
$\{r=  R\}$. Let us start from $\bar{u}\in C^1(\mathbb{R}^3) \subset
C^{\alpha}(\overline{B_R})$. Then we have $f^{\rho}(\bar{u}) \in C^{\alpha}(\overline{B_R})$, since $\frac{1}{\gamma-1}>1$. The regularity theorem ( e.g., \cite[Theorem 4.5]{GilbergT} )
guarantees that the solution $\bar{u}$ of
$-\triangle \bar{u}=4\pi\mathsf{G}f^{\rho}(\bar{u})$ belongs to $C^{2,\alpha}(\overline{B_R})$. Then we have $f^{\rho}(\bar{u}) \in C^{1,\alpha}(\overline{B_R})$, since $\frac{1}{\gamma-1}>1+\alpha$. Again the regularity theorem guarantees that $\bar{u} \in C^{3,\alpha}(\overline{B_R})$.
Since $\bar{\rho}^{\gamma-1}$ is a smooth function of $\bar{u}$, we have
$\bar{\rho}^{\gamma-1} \in C^{3,\alpha}(\overline{B_R})$. Of course $\bar{\Phi}=\bar{u}+\mbox{Const.}  \in C^{3,\alpha}(\overline{B_R})$. Since $\bar{S}$ is a smooth function of
$\bar{\rho}^{\gamma-1}$ , we have $\bar{S}\in C^{3,\alpha}(\overline{B_R})$. Summing up, we have verified the condition 2) of Definition \ref{Def.2}.
As for the condition 2), we note that
$$\frac{1}{r}\frac{du}{dr}=\frac{1}{r}\frac{1}{\rho}\frac{dP}{dr}=-4\pi\mathsf{G}
\frac{1}{r^3}\int_0^r
\rho(\acute{r})\acute{r}^2d\acute{r} \rightarrow -\frac{4\pi}{3}\mathsf{G}\rho_{\mathsf{O}}$$
as $r \rightarrow +0$ so that
$$-\lim_{r\rightarrow +0}\frac{1}{r}\frac{dP}{dr}=\frac{4\pi}{3}\mathsf{G}\rho_{\mathsf{O}}^2$$
and
$$
-\lim_{r\rightarrow+0}\frac{1}{r}\frac{d\rho}{dr}
=\Big[\gamma+\frac{\gamma-1}{\mathsf{C}_V}\omega\frac{d\Sigma}{d\omega}\Big]_{\omega=\rho_{\mathsf{O}}^{\gamma-1}}^{-1}\frac{4\pi}{3}\mathsf{G}\rho_{\mathsf{O}}
\bar{P}(O).$$
$\square$

\begin{Remark}
In the barotropic case, the quantity $u$ means the specific enthalpy. But in the general baroclinic case, $u$ is not the specific enthalpy $\chi$ which should be defined as
$$\chi:=\mathsf{C}_VT+\frac{P}{\rho}=\frac{\gamma}{\gamma-1}\frac{P}{\rho},$$
$T=P/((\gamma-1)\mathsf{C}_V\rho)$ being the absolute temperature. In fact we have
$$\frac{d}{dr}(u-\chi)=-\frac{1}{\mathsf{C}_V}\omega\frac{d}{d\omega}\Sigma(\omega)
\Big|_{\omega=\bar{\rho}^{\gamma-1}}\overline{\frac{P}{\rho^2}\frac{d\rho}{dr} }$$
does not vanish if $\Sigma$ is not constant, that is, $\bar{S}$ is not constant, for
$d\bar{\rho}/dr <0$.
\end{Remark}

Hereafter in this article we fix an admissible spherically symmetric equilibrium $(\bar{\rho}, \bar{S})$, and denote
\begin{equation}
\rho_{\mathsf{O}}=\bar{\rho}(O)(>0),
\quad
\rho_{\mathsf{O}1}=-\lim_{r\rightarrow +0}\frac{1}{r}\frac{d\bar{\rho}}{dr}(>0), 
\quad
P_{\mathsf{O}1}=-\lim_{r\rightarrow +0}\frac{1}{r}\frac{d\bar{P}}{dr}(>0).
\end{equation}\\

As for the Schwarzschild's discriminant, since
$\displaystyle \mathscr{A}=-\frac{1}{\gamma \mathsf{C}_V}\frac{d\bar{S}}{dr} $
and we are assuming 
$ \bar{S} \in 
C^{2,\alpha}(\overline{B_R})$,
 we have
$$
\mathscr{A}=O(r) 
$$
as $r\rightarrow +0$. (Recall Lemma \ref{Lem.0}.)

If $\bar{S}=\Sigma(\bar{\rho}^{\gamma-1})$, then we have
$$
\mathscr{A}=-\frac{1}{\gamma \mathsf{C}_V}\frac{d\bar{S}}{dr}
=-\frac{1}{\gamma \mathsf{C}_V}\frac{d}{d\eta}\Sigma(\eta)\Big|_{\eta=\bar{\rho}^{\gamma-1}}
\frac{d}{dr}\bar{\rho}^{\gamma-1}. 
$$\\

\begin{Remark}

In Theorem \ref{Th.1} we used the auxiliary relation $\bar{S}=\Sigma(\bar{\rho}^{\gamma-1})$ in order to construct an admissible spherically symmetric equilibrium $(\bar{\rho}, \bar{S})$. In fact, if $\bar{\rho}(\mbox{\boldmath$x$})=\bar{\rho}^{\sharp}(r), \bar{S}(\mbox{\boldmath$x$})=\bar{S}^{\sharp}(r)$  is an admissible equilibrium, there should exist a function $\Sigma$ such that
$\bar{S}(\mbox{\boldmath$x$})=\Sigma((\bar{\rho}(\mbox{\boldmath$x$}))^{\gamma-1})$ for $\forall \mbox{\boldmath$x$} \in B_R$, since $\bar{\rho}^{\sharp}$ is supposed to be monotone.
 Actually it is sufficient to put $\Sigma=\bar{S}^{\sharp}\circ \eta^{-1}$, where
$\eta^{-1}$ is the inverse function of the monotone function $\eta:
r \mapsto (\bar{\rho}^{\sharp}(r))^{\gamma-1}$.
But, once we have constructed and fixed the equilibrium $(\bar{\rho}, \bar{S})$, we should forget this relation $\bar{S}=\Sigma(\bar{\rho}^{\gamma-1})$, namely, the constraint supposed for the perturbed state variables $\rho=\bar{\rho}+\Delta\rho, P=\bar{P}+\Delta P, S=\bar{S}+\Delta S$ is nothing other than the relation
$$P=\rho^{\gamma}\exp\Big(\frac{S}{\mathsf{C}_{V}}),$$
and the relation $S=\Sigma(\rho^{\gamma-1})$ is never supposed for the perturbed motion. In fact, if the relation $S=\Sigma(\rho^{\gamma-1})$ was still supposed for the perturbed motion, then 
the relation 
$$P=\rho^{\gamma}\exp\Big(\frac{\Sigma(\rho^{\gamma-1})}{\mathsf{C}_V}\Big)$$
would give a relation of 
the form $P=F(\rho)=\mathsf{A}\rho^{\gamma}(1+\Lambda(\rho^{\gamma-1}))$ of a barotropic motion, which is not the subject of the present discussion. 
In other words, if the relation
$S = \Sigma(\rho^{\gamma-1})$ holds throughout the motion and 
if $d\Sigma(\eta)/d\eta \not=0$ a.e.-$\eta$, then \eqref{EPc} implies 
$$\frac{D\rho}{Dt}=\frac{\partial\rho}{\partial t}+\sum_{k=1}^3v^k\frac{\partial \rho}{\partial x^k}=0, $$
that is, the flow is incompressible, 
or $\rho(t,\mbox{\boldmath$\varphi$}(t,\mbox{\boldmath$x$}))=
\rho(0,\mbox{\boldmath$x$})=
\bar{\rho}(\mbox{\boldmath$x$})$,
 and
$$\mathrm{div}\mbox{\boldmath$v$}=\sum_{k=1}^3\frac{\partial v^k}{\partial x^k}=0.$$
This is not the situation we are concerned with.
\end{Remark}

\section{Self-adjoint realization of $\mbox{\boldmath$L$}$}

We are considering the integro-differential operator
\begin{equation}
\mbox{\boldmath$L$}\mbox{\boldmath$\xi$}=\frac{1}{\rho}\mathrm{grad}\delta P-
\frac{\delta\rho}{\rho^2}\mathrm{grad} P +\mathrm{grad}\delta\Phi,
\end{equation}
where
\begin{subequations}
\begin{align}
\delta\rho&=-\mathrm{div}(\rho\mbox{\boldmath$\xi$}), \\ 
\delta P&=\frac{\gamma P}{\rho}\delta \rho+\gamma \mathscr{A}P(\mbox{\boldmath$\xi$}|\mbox{\boldmath$e$}_r) ,\\ 
\delta \Phi&=-4\pi\mathsf{G}\mathcal{K}(\delta\rho). 
\end{align}
\end{subequations}
Here and hereafter we use the following 

\begin{Notation}
The bars to denote the quantities evaluated at the fixed equilibrium are omitted, that is, $\rho, P$ etc stand for $\bar{\rho}, \bar{P}$ etc.
\end{Notation}

Let us consider the operator $\mbox{\boldmath$L$}$ in the Hilbert space $\mathfrak{H}=L^2((B_R, \rho d\mbox{\boldmath$x$}), \mathbb{C}^3)$ endowed with the norm $\|\mbox{\boldmath$\xi$}\|_{\mathfrak{H}}$ defined by
\begin{equation}
\|\mbox{\boldmath$\xi$}\|_{\mathfrak{H}}^2=\int_{B_R}|\mbox{\boldmath$\xi$}(\mbox{\boldmath$x$})|^2\rho(\mbox{\boldmath$x$})d\mbox{\boldmath$x$}
\end{equation}

We shall use

\begin{Notation}
For complex number $z=x+\sqrt{-1}y, x, y \in \mathbb{R}$, the complex conjugate is denoted by
$z^*=x-\sqrt{-1}y$. Thus, for $\displaystyle \mbox{\boldmath$\xi$}=\sum \xi^k\frac{\partial}{\partial x^k}, \xi^k \in
\mathbb{C}$, we denote $\displaystyle \mbox{\boldmath$\xi$}^*=\sum (\xi^k)^*\frac{\partial}{\partial x^k}$.
\end{Notation}

First we observe $\mbox{\boldmath$L$}$ restricted on $C_0^{\infty}(B_R, \mathbb{C}^3)$. Let us write
\begin{align}
\mbox{\boldmath$L$}\mbox{\boldmath$\xi$}&=\mbox{\boldmath$L$}_0\mbox{\boldmath$\xi$} +4\pi\mathsf{G}\mbox{\boldmath$L$}_1\mbox{\boldmath$\xi$}, \\
\mbox{\boldmath$L$}_0\mbox{\boldmath$\xi$}&=
\mathrm{grad}\Big(-\frac{\gamma P}{\rho^2}\mathrm{div}(\rho\mbox{\boldmath$\xi$})+
\frac{\gamma\mathscr{A}P}{\rho}(\mbox{\boldmath$\xi$}|\mbox{\boldmath$e$}_r)\Big) + \nonumber \\
&+\frac{\gamma\mathscr{A}P}{\rho^2}\Big(-\mathrm{div}(\rho\mbox{\boldmath$\xi$})+
\frac{d\rho}{dr}(\mbox{\boldmath$\xi$}|\mbox{\boldmath$e$}_r)\Big)\mbox{\boldmath$e$}_r, \\
\mbox{\boldmath$L$}_1\mbox{\boldmath$\xi$}&=\mathrm{grad}\mathcal{K}(\mathrm{div}(\rho\mbox{\boldmath$\xi$})).
\end{align}

Using this expression for $\mbox{\boldmath$\xi$}_{(\mu)} \in C_0^{\infty}(B_R), \mu=1,2$, we have the following formula by integration by parts:
\begin{align*}
(\mbox{\boldmath$L$}_0\mbox{\boldmath$\xi$}_{(1)}|\mbox{\boldmath$\xi$}_{(2)})_{\mathfrak{H}}&=
\int\frac{\gamma P}{\rho^2}\mathrm{div}(\rho\mbox{\boldmath$\xi$}_{(1)})\mathrm{div}(\rho\mbox{\boldmath$\xi$}_{(2)}^*) + \\
&+\int\frac{\gamma\mathscr{A}P}{\rho}
\Big[(\mbox{\boldmath$\xi$}_{(1)}|\mbox{\boldmath$e$}_r)\cdot\mathrm{div}(\rho\mbox{\boldmath$\xi$}_{(2)}^*)
-\mathrm{div}(\rho\mbox{\boldmath$\xi$}_{(1)})\cdot (\mbox{\boldmath$\xi$}_{(2)}|\mbox{\boldmath$e$}_r)^*\Big] + \\
&+\int\frac{\gamma\mathscr{A}P}{\rho}\frac{d\rho}{dr}(\mbox{\boldmath$\xi$}_{(1)}|\mbox{\boldmath$e$}_r)
(\mbox{\boldmath$\xi$}_{(2)}|\mbox{\boldmath$e$}_r)^*, \\
(\mbox{\boldmath$L$}_1\mbox{\boldmath$\xi$}_{(1)}|\mbox{\boldmath$\xi$}_{(2)})_{\mathfrak{H}}&=
-\int\mathcal{K}(\mathrm{div}(\rho\mbox{\boldmath$\xi$}_{(1)})\cdot \mathrm{div}(\rho\mbox{\boldmath$\xi$}_{(2)}^*).
\end{align*}

Thanks to the symmetry of $\mathcal{K}$, we have
$$ (\mbox{\boldmath$L$}\mbox{\boldmath$\xi$}_{(1)}|\mbox{\boldmath$\xi$}_{(2)})_{\mathfrak{H}}=
(\mbox{\boldmath$\xi$}_{(1)}|\mbox{\boldmath$L$}\mbox{\boldmath$\xi$}_{(2)})_{\mathfrak{H}},
$$
that is, $\mbox{\boldmath$L$}$ restricted on $C_0^{\infty}(B_R)$ is a symmetric operator.
Of course $C_0^{\infty}(B_R)$ is dense in $\mathfrak{H}$.

Moreover we have
\begin{align*}
(\mbox{\boldmath$L$}_0\mbox{\boldmath$\xi$}|\mbox{\boldmath$\xi$})&=
\int\frac{\gamma P}{\rho^2}|\mathrm{div}(\rho\mbox{\boldmath$\xi$})|^2 + \\
&+2\sqrt{-1}\mathfrak{Im}\Big[
\int\frac{\gamma\mathscr{A}P}{\rho}(\mbox{\boldmath$\xi$}|\mbox{\boldmath$e$}_r)\cdot\mathrm{div}(\rho\mbox{\boldmath$\xi$}^*)\Big] +
\int\frac{\gamma\mathscr{A}P}{\rho}\frac{d\rho}{dr}|(\mbox{\boldmath$\xi$}|\mbox{\boldmath$e$}_r)|^2.
\end{align*}

Since $\mathscr{A} \in C^{1,\alpha}(\overline{B_R})$, we have
$$|\mathscr{A}|\sqrt{\frac{\gamma P}{\rho}}\leq C_1$$
on $0<r\leq R$, for $P/\rho =O(R-r)$. Therefore
\begin{align*}
\Big|\int
\frac{\gamma\mathscr{A}P}{\rho}(\mbox{\boldmath$\xi$}|\mbox{\boldmath$e$}_r)
\mathrm{div}(\rho\mbox{\boldmath$\xi$})^*\Big|
&\leq C_1\int
\sqrt{ \frac{\gamma P}{\rho} }|(\mbox{\boldmath$\xi$}|\mbox{\boldmath$e$}_r)|
|\mathrm{div}(\rho\mbox{\boldmath$\xi$})| \\
&\leq \frac{C_1}{2}\Big[
\frac{1}{\epsilon}\int \rho |(\mbox{\boldmath$\xi$}|\mbox{\boldmath$e$}_r)|^2+
\epsilon\int
\frac{\gamma P}{\rho^2}|\mathrm{div}(\rho\mbox{\boldmath$\xi$})|^2\Big] \\
&\leq \frac{C_1}{2}\Big[
\frac{1}{\epsilon}\|\mbox{\boldmath$\xi$}\|_{\mathfrak{H}}^2+
\epsilon\int\frac{\gamma P}{\rho^2}|\mathrm{div}\rho\mbox{\boldmath$\xi$})|^2\Big].
\end{align*}

Since 
$\displaystyle \frac{P}{\rho}\frac{d\rho }{dr}=O(\rho)$, we have 
$$\Big|\frac{\gamma\mathscr{A}P}{\rho}\frac{d\rho}{dr}\Big|\leq C \rho $$
Therefore we have
$$
\Big|\int\frac{\gamma\mathscr{A}P}{\rho}\frac{d\rho}{dr}|(\mbox{\boldmath$\xi$}|\mbox{\boldmath$e$}_r)|^2\Big|\leq C_2\|\mbox{\boldmath$\xi$}\|_{\mathfrak{H}}^2.
$$

Thus
$$
(\mbox{\boldmath$L$}_0\mbox{\boldmath$\xi$}|\mbox{\boldmath$\xi$})_{\mathfrak{H}}\geq
\Big(1-\frac{\epsilon C_1}{2}\Big)\int
\frac{\gamma P}{\rho^2}|\mathrm{div}(\rho\mbox{\boldmath$\xi$})|^2 
-\Big(\frac{C_1}{2\epsilon}+C_2\Big)\|\mbox{\boldmath$\xi$}\|_{\mathfrak{H}}^2. $$
Taking $\epsilon$ so small that $1-\frac{\epsilon C_1}{2} \geq 0$, we get
$$
(\mbox{\boldmath$L$}_0\mbox{\boldmath$\xi$}|\mbox{\boldmath$\xi$})_{\mathfrak{H}}\geq
-\Big(\frac{C_1}{2\epsilon}+C_2\Big)\|\mbox{\boldmath$\xi$}\|_{\mathfrak{H}}^2. $$

On the other hand, it is known that
$$(\mbox{\boldmath$L$}_1\mbox{\boldmath$\xi$}|\mbox{\boldmath$\xi$})_{\mathfrak{H}} \geq -\rho_{\mathsf{O}}\|\mbox{\boldmath$\xi$}\|_{\mathfrak{H}}
^2.$$
For a proof , see \cite[Proof of Proposition 2]{JJTM}. 

Summing up, $\mbox{\boldmath$L$}$ is bounded from below in $\mathfrak{H}$. Therefore, thanks to \cite[Chapter VI, Section 2.3]{Kato}, we have

\begin{Theorem}
The integro-differential operator $\mbox{\boldmath$L$}$ on $C_0^{\infty}(B_R, \mathbb{C}^3)$ admits the Friedrichs extension, which is a self-adjoint operator, in $\mathfrak{H}$.
\end{Theorem}

We want to clarify the spectral property of the self-adjoint operator $\mbox{\boldmath$L$}$. But this task has not yet been completely done. 

At least we can claim that the spectrum of $\mbox{\boldmath$L$}$ cannot be of the Sturm-Liouville type in the sense defined in \cite{JJTM},
(that is, the spectrum $\sigma(T)$ of a self-adjoint operator $T$ in a Hilbert space $\mathsf{X}$ is said to be of the Sturm-Liouville type if $\sigma(T)$ consists of isolated eigenvalues with finite multiplicities,)
since $\mathrm{dim}\mathsf{N}(\mbox{\boldmath$L$})=\infty$, where $\mathsf{N}(T)$ denotes the kernel of the operator $T$ that is, $\mathsf{N}(T)=\{ x \in \mathsf{D}(T) | Tx =0\}$. In fact, if we consider a scalar field $a$ on $B_R$ given by a function $a^{\sharp}: (r,\vartheta,\phi) \mapsto a(\mbox{\boldmath$x$})$ which belongs to $C_0^{\infty}(]0,R[\times ]0, \pi[\times \mathbb{S}^1)$, then the field
$$\mbox{\boldmath$\xi$}=\frac{1}{\sin\vartheta}\frac{\partial a}{\partial \phi}\mbox{\boldmath$e$}_r
-\frac{\partial a}{\partial\vartheta}\mbox{\boldmath$e$}_{\phi} $$
belongs to $C_0^{\infty}(B_R)$ and satisfies $\mathrm{div}(\rho\mbox{\boldmath$\xi$})=0$ and
$(\mbox{\boldmath$\xi$}|\mbox{\boldmath$e$}_r)=0$, therefore it belongs to $\mathsf{N}(\mbox{\boldmath$L$})$. Since the dimension of spaces of such $a^{\sharp}$ is infinite, we see
$\mathrm{dim}\mathsf{N}(\mbox{\boldmath$L$})=\infty$. 

In the work \cite{JJTM}, we proved that, when $S$ is constant so that $\mathscr{A}=0$, then the spectrum of the operator $\mbox{\boldmath$L$}$ is $\{0\}\cup \{\lambda_n , n=1,2,\cdots \}$, where 0 is an essential spectrum and $\lambda_n$ are eigenvalues of finite multiplicities,
$\lambda_n \rightarrow +\infty$ as $n \rightarrow \infty$, provided that $\mbox{\boldmath$L$}$ is considered in the Hilbert space 
$$\mathfrak{F}=\mathfrak{H}\cap\{ \mbox{\boldmath$\xi$} | 
\mathrm{div}(\rho\mbox{\boldmath$\xi$})\in \mathfrak{G}\}, $$
while 
$$\mathfrak{G}=
L^2(B_R, \frac{1}{\rho}\frac{dP}{d\rho}d\mbox{\boldmath$x$})\cap\{g | \int gd\mbox{\boldmath$x$}=0\}.$$

However, when $S$ is not constant and $\mathscr{A}$ does not identically vanish, such a situation
cannot be expected, but the so called g-modes can appear, that is, there  can exist a sequence of eigenvalues which accumulates to 0. In Section 5 we  shall prove in a mathematically rigorous way that actually it is the case  under the 
assumption that $\inf\mathscr{N}^2 /r^2>0$ and a set of 
 conditions which guarantees the smallness of the effect of the perturbation of the self-gravitation.

\section{Solutions represented by spherical harmonics}

In this section we consider the perturbation $\mbox{\boldmath$\xi$}$ of the particular form
\begin{equation}
\mbox{\boldmath$\xi$}=V^r(t,r)Y_{lm}(\vartheta,\phi)\mbox{\boldmath$e$}_r+
V^h(t,r)\nabla_sY_{lm}+V^v(t,r)\nabla_s^{\perp}Y_{lm}.
\end{equation}
Here $l, m \in \mathbb{Z}, 0 \leq l, |m|\leq l$, and $Y_{lm}$ is the spherical harmonics:
\begin{align*}
&Y_{lm}(\vartheta,\phi)=\sqrt{\frac{2l+1}{4\pi}\frac{(l-m)!}{(l+m)!}}
P_l^m(\cos\vartheta)e^{\sqrt{-1}m\phi}, \\
&Y_{l,-m}=(-1)^mY_{lm}^*
\end{align*}
for $m \geq 0$, while $P_l^m$ is the associated Legendre function given by
$$P_l^m(\zeta)=\frac{(-1)^m}{2^ll!}(1-\zeta^2)^{m/2}\Big(\frac{d}{d\zeta}\Big)^{m+l}(\zeta^2-1)^l.$$
$l$ is called the 'degree' of the mode and $m$ is called the 'azimuthal order' of the mode. See \cite{Jackson}. We use

\begin{Notation}
We denote
\begin{subequations}
\begin{align}
&\nabla_sf:=\frac{\partial f}{\partial\vartheta}\mbox{\boldmath$e$}_{\vartheta}+
\frac{1}{\sin\vartheta}
\frac{\partial f}{\partial\phi}\mbox{\boldmath$e$}_{\phi}, \\
&\nabla_s^{\perp}f:=\frac{1}{\sin\vartheta}
\frac{\partial f}{\partial\phi}\mbox{\boldmath$e$}_{\vartheta}-\frac{\partial f}{\partial\vartheta}\mbox{\boldmath$e$}_{\phi}.
\end{align}
\end{subequations}
\end{Notation}

Note that
\begin{subequations}
\begin{align}
&\mathrm{div}(\psi(r)Y_{lm}\mbox{\boldmath$e$}_r)=
\frac{1}{r^2}\frac{d}{dr}(r^2\psi)Y_{lm}, \\
&\mathrm{div}(\psi(r)\nabla_sY_{lm})=-\frac{l(l+1)}{r}\psi(r)Y_{lm}, \\
&\mathrm{div}(\psi(r)\nabla_s^{\perp}Y_{lm})=0.
\end{align}
\end{subequations}\\

Then \eqref{drho}, \eqref{dP}, \eqref{dPhi} read
\begin{subequations}
\begin{align}
\delta\rho&=\delta\check{\rho}(t,r)Y_{lm}(\vartheta,\phi), \\
\delta P&=\delta\check{P}(t,r)Y_{lm}(\vartheta,\phi), \\
\delta\Phi&=\delta\check{\Phi}(t,r)Y_{lm}(\vartheta, \phi),
\end{align}
\end{subequations}
where
\begin{subequations}
\begin{align}
\delta\check{\rho}&=-\frac{1}{r^2}\frac{d}{dr}(r^2\rho V^r)+\frac{l(l+1)}{r}\rho V^h, \label{4.6a}\\
\delta\check{P}&=\frac{\gamma P}{\rho}\delta\check{\rho}+\gamma
\mathscr{A}P V^r \nonumber \\
&=-\frac{\gamma P}{r^2\rho}\frac{d}{dr}(r^2\rho V^r)+\gamma \mathscr{A}PV^r+
l(l+1)\frac{\gamma P}{r}V^h, \label{4.6b}\\
\delta\check{\Phi}&=-4\pi\mathsf{G}\mathcal{H}_l(\delta\check{\rho}) \nonumber \\
&=4\pi\mathsf{G}\mathcal{H}_l\Big(
\frac{1}{r^2}\frac{d}{dr}(r^2\rho V^r)-\frac{l(l+1)}{r}\rho V^h\Big).\label{4.6c}
\end{align}
\end{subequations}
Here the integral operator $\mathcal{H}_l$ is defined by

\begin{Definition}
We put
\begin{equation}
\mathcal{H}_lf(r)=\frac{1}{2l+1}\Big[
\int_r^{\infty}f(\acute{r})\Big(\frac{r}{\acute{r}}\Big)^l\acute{r}d\acute{r}+
\int_0^rf(\acute{r})\Big(\frac{r}{\acute{r}}\Big)^{-(l+1)}\acute{r}d\acute{r}\Big],
\end{equation}
provided that $f \in L^2([0,+\infty[, r^2dr)$ and
$f(r)=0$ for $r\geq R$. 
\end{Definition}

The equation \eqref{Eqxi} reads
\begin{equation}
\frac{\partial^2V^r}{\partial t^2}+L_l^r=0, \quad \frac{\partial^2V^h}{\partial t^2}+L_l^h=0, \quad \frac{\partial^2V^v}{\partial t^2} =0, \label{EqVl}
\end{equation}
where
\begin{subequations}
\begin{align}
&L_l^r=\frac{1}{\rho}\frac{d}{dr}\delta\check{P}-\frac{1}{\rho^2}\frac{dP}{dr}\delta\check{\rho}+\frac{d}{dr}\delta\check{\Phi}, \\
&L_l^h=\frac{1}{r}\Big(\frac{\delta\check{P}}{\rho}+\delta\check{\Phi}\Big).
\end{align}
\end{subequations}
We mean
\begin{equation}
\mbox{\boldmath$L$}(V^rY_{lm}\mbox{\boldmath$e$}_r+V^h\nabla_sY_{lm}
+V^v\nabla_s^{\perp}Y_{lm})=L_l^rY_{lm}\mbox{\boldmath$e$}_r+
L_l^h\nabla_sY_{lm}.
\end{equation}

We are going to analyze the operator 
\begin{equation}
\vec{L}_l=\begin{bmatrix}
L_l^r \\
\\
L_l^h
\end{bmatrix}
\end{equation}
which acts on
\begin{equation}
\vec{V}=\begin{bmatrix}
V^r \\
V^h
\end{bmatrix}.
\end{equation}

Actually we can neglect the component $V^v(t,r)$, which should be an arbitrary affine function of $t$ in order to satisfy  \eqref{EqVl}, and we are going to consider the eigenvalue problem
\begin{equation}
\vec{L}_l\vec{V}=\lambda\vec{V}. 
\end{equation}
\\

As for the integral operator $\mathcal{H}_l$, we shall keep in mind the following lemma, which is easy to prove:
\begin{Lemma}\label{Lem.2}
1) Let $l =0$ and $f \in L^2([0,+\infty[, r^2dr), f(r)=0$ for $r \geq R$. Then $H=\mathcal{H}_0f \in C^1(]0,+\infty[)$ and satisfies
\begin{align*}
&H(r)=O(1), \quad \frac{d}{dr}H(r) =O(r^{-\frac{1}{2}})\quad \mbox{as}\quad r \rightarrow +0, \\
&H(r)=\frac{C}{r}\quad\mbox{as}\quad r \geq R,
\end{align*}
\begin{equation}
-\frac{1}{r^2}\frac{d}{dr}\Big(r^2\frac{dH}{dr}\Big)=f. \label{EqH0}
\end{equation}
Here $C$ is the constant given by
$$C=\int_0^Rf(r)r^{2}dr.
$$

Conversely, if $H$ is absolutely continuous and satisfies \eqref{EqH0}
on $]0,+\infty[$, there exist  constants 
$C_1,C_2$ such that
$$H=\mathcal{H}_0f +C_1  +\frac{C_2}{r},
$$
therefore,  $H\in L^{\infty}([0,+\infty[)$and $ H \rightarrow  0 $
as $ r \rightarrow +\infty$ if and only if $H=\mathcal{H}_0f$.

2) Let $l \geq 1$ and $f \in L^2([0,+\infty[, r^2dr), f(r)=0$ for $r \geq R$. Then $H=\mathcal{H}_lf \in C^1(]0,+\infty[)$ and satisfies
\begin{align*}
&H(r)=O(r^{\frac{1}{2}}), \quad \frac{d}{dr}H(r)=O(r^{-\frac{1}{2}})\quad \mbox{as}\quad r \rightarrow +0, \\
&H(r)=Cr^{-(l+1)}\quad\mbox{as}\quad r \geq R,
\end{align*}
\begin{equation}
-\frac{1}{r^2}\frac{d}{dr}\Big(r^2\frac{dH}{dr}\Big)+
\frac{l(l+1)}{r^2}H=f. \label{EqHl}
\end{equation}
Here $C$ is the constant given by
$$C=\frac{1}{2l+1}\int_0^Rf(r)r^{l+2}dr.
$$

Conversely, if $H$ is absolutely continuous and satisfies \eqref{EqHl}
on $]0,+\infty[$, there exists  constants 
$C_1,C_2$ such that
$$H=\mathcal{H}_lf +C_1 r^l +C_2 r^{-(l+1)},
$$ and $H\in L^2([0,+\infty[)$ if and only if $H=\mathcal{H}_lf$.
\end{Lemma}

\subsection{Case $l=0$}

First let us consider the case $l=0$, when only $m=0$ is possible, and 
$\displaystyle Y_{00}=\frac{1}{\sqrt{4\pi}}$. We are considering
\begin{equation}
\mbox{\boldmath$\xi$}=\frac{V(r)}{\sqrt{4\pi}}\mbox{\boldmath$e$}_r,
\end{equation}
where we write $V$ instead of $V^r$, while we need not consider $V^h$. 

We are concerned with the operator
\begin{equation}
L_0=\frac{1}{\rho}\frac{d}{dr}\delta\check{P}-\frac{1}{\rho^2}\frac{dP}{dr}\delta\check{\rho}+\frac{d}{dr}\delta\check{\Phi}
\end{equation}
with
\begin{subequations}
\begin{align}
\delta\check{\rho}&=-\frac{1}{r^2}\frac{d}{dr}(r^2\rho V), \\
\delta\check{P}&=\frac{\gamma P}{\rho}\delta\check{\rho}+\gamma
\mathscr{A}P V \nonumber \\
&=-\frac{\gamma P}{r^2\rho}\frac{d}{dr}(r^2\rho V)+\gamma \mathscr{A}PV, \\
\delta\check{\Phi}&=-4\pi\mathsf{G}\mathcal{H}_0(\delta\check{\rho}) \nonumber \\
&=4\pi\mathsf{G}\mathcal{H}_0\Big(
\frac{1}{r^2}\frac{d}{dr}(r^2\rho V)\Big)
\end{align}
\end{subequations}
so that
\begin{equation}
\frac{d}{dr}\delta\check{\Phi}=-4\pi\mathsf{G}\rho V.
\end{equation}
We mean
$$\mbox{\boldmath$L$}\Big(\frac{V}{\sqrt{4\pi}}\mbox{\boldmath$e$}_r\Big)=\frac{L_0V}{\sqrt{4\pi}}\mbox{\boldmath$e$}_r.
$$\\

Introducing the variable $\psi$ by
\begin{equation}
V=r\psi,
\end{equation}
and putting
\begin{equation}
L^{\mathsf{ss}}\psi=\frac{1}{r}L_0(r\psi),
\end{equation}
we analyze the differential operator operator $L^{\mathsf{ss}}$ in the Hilbert space
$L^2([0,R], \rho r^4dr), \mathbb{C})$, since
$$\|\mbox{\boldmath$\xi$}\|_{\mathfrak{H}}^2=\int_0^R|\psi(r)|^2\rho(r)r^4dr $$
for $\displaystyle \mbox{\boldmath$\xi$}=\frac{r\psi}{\sqrt{4\pi}}\mbox{\boldmath$e$}_r$.\\

We claim

\begin{Theorem}
The operator $L^{\mathsf{ss}}$ on $C_0^{\infty}(]0,R[)$ admits the Friedrichs extension, a self-adjoint operator bounded from below in $L^2([0,R], \rho r^4dr)$, and its spectrum consists of simple eigenvalues $\lambda_1^{\mathsf{ss}}<\lambda_2^{\mathsf{ss}}<\cdots <
\lambda_n^{\mathsf{ss}}<\cdots \rightarrow +\infty$.
\end{Theorem}

Proof. First we write $L^{\mathsf{ss}}$ as
\begin{equation}
L^{\mathsf{ss}}\psi=-\frac{1}{r^4\rho}
\frac{d}{dr}\Big(\gamma r^4P\frac{d\psi}{dr}\Big)+q_{00}(r)\psi,
\end{equation}
where
\begin{align}
q_{00}(r)&=-\frac{\gamma P}{\rho^2}\frac{d^2\rho}{dr^2}
-\frac{\gamma-1}{\rho^2}\frac{dP}{dr}\frac{d\rho}{dr}
+\frac{\gamma P}{\rho^3}\Big(\frac{d\rho}{dr}\Big)^2 \nonumber \\
&-\frac{\gamma P}{r\rho^2}\frac{d\rho}{dr}
-\frac{3(\gamma-1)}{r\rho}\frac{dP}{dr}
+\frac{1}{r\rho}\frac{d}{dr}(\gamma r\mathscr{A}P)-4\pi\mathsf{G}\rho.
\label{q00}
\end{align}

We see that 
$$ |q_{00}(r)|\leq C\quad\mbox{for}\quad 0<r<R.$$
In fact, although each term in the first line of the right-hand side of
\eqref{q00} is of order $(R-r)^{-1}$, these singularities are canceled
after the summation, which turns out to be
$$-\frac{\gamma}{\gamma-1}\Big[
\frac{d^2\eta}{dr^2}+\frac{\gamma^1}{\gamma \mathsf{C}_V}\eta^{\frac{2-\gamma}{\gamma-1}}
\frac{dS}{dr}\frac{d\eta}{dr}\Big]\exp\Big(\frac{S}{\mathsf{C}_V}\Big)=O(1),
$$
where $\eta:=\rho^{\gamma-1} \in C^{2,\alpha}(\overline{B_R})$.

So, we perform the Liouville transformation of $L^{\mathsf{ss}}\psi=\lambda\psi$ to
$$-\frac{d^2y}{dx^2}+q_0(x)y=\lambda y,
$$
where
\begin{align*}
x&=\int\sqrt{\frac{c}{a}}dr,\qquad y=(ac)^{\frac{1}{4}}\psi, \\
q_0&=q_{00}+\frac{1}{4}\frac{a}{c}Q, \\
Q&=\frac{d^2}{dr^2}\log(ac)-\frac{1}{4}\Big(\frac{d}{dr}
\log(ac)\Big)^2+\Big(\frac{d}{dr}\log a\Big)\Big(
\frac{d}{dr}\log(ac)\Big)
\end{align*}
with
$$a=\gamma r^4P,\qquad c=r^4\rho.$$
(See \cite[p. 275, Theorem 6]{BirkhoffR} or \cite[p.110]{Yosida}.)

Since 
$$\sqrt{\frac{c}{a}}=\sqrt{\frac{\rho}{\gamma P}},$$
we can put
$$x=\int_0^r\sqrt{\frac{\rho}{\gamma P}}(\acute{r})d\acute{r} $$
so that $x$ runs on the interval $[0,x_+]$, where
$x_+=\int_0^R\sqrt{\rho/\gamma P}dr$. We have
\begin{align*}
&x \sim \sqrt{\frac{\rho}{\gamma P}}\Big|_Or\quad\mbox{as}\quad r \rightarrow +0 \\
&x_+-x \sim C_R(R-r)^{\frac{1}{2}}\quad\mbox{as}\quad r \rightarrow R-0
\end{align*}
with a positive constant $C_R$. 

It can be verified that
\begin{align*}
&q_0 \sim \frac{2}{x^2} \quad\mbox{as}\quad x\rightarrow +0 \\
&q_0 \sim
\frac{(\gamma+1)(3-\gamma)}{4(\gamma-1)^2}\frac{1}{(x_+-x)^2}
\quad\mbox{as}\quad x \rightarrow x_+-0.
\end{align*}
Hence the assertion follows from \cite[p.159, Theorem X.10]{ReedS}. $\square$

\subsection{Case $l \geq 1$}

Suppose $\l \geq1$. \\

Let us consider the Hilbert space $\mathfrak{X}_l$ of functions $\vec{f}=(f^r, f^h)^{\top}$ defined on $[0,R[$ endowed with the norm $\|\vec{f}\|_{\mathfrak{X}_l}$ given by

\begin{equation}
\|\vec{f}\|_{\mathfrak{X}_l}^2=\int_0^R
(|f^r(r)|^2+l(l+1)|f^h(r)|^2)\rho r^2dr.
\end{equation}
Of course $\vec{f} =(f^r, f^h)^{\top} \in \mathfrak{X}_l$ if and only if
$$\mbox{\boldmath$\xi$}=\vec{f}^{\flat/lm}:=
f^r(r)Y_{lm}\mbox{\boldmath$e$}_r+f^h(r)\nabla_sY_{lm} \in \mathfrak{H}$$ for $|m|\leq l$, and
$$\|\vec{f}^{\flat/lm}\|_{\mathfrak{H}}=\sqrt{4\pi}\|\vec{f}\|_{\mathfrak{X}_l}. $$
we consider the operator $\vec{L}_l$ in $\mathfrak{X}_l$. We claim

\begin{Theorem}\label{Th.3}
The integro-differential operator $\vec{L}_l$ on $C_0^{\infty}([0,R[, \mathbb{C}^2)$ admits the Friedrichs extension, which is a self-adjoint operator bounded from below, in $\mathfrak{X}_l$.
\end{Theorem}

Proof.  First we look at $L_l^r, L_l^h$ by writing them as

\begin{subequations}
\begin{align}
L_l^r&=-\frac{\gamma P^{\frac{1}{\gamma}}}{\rho}\frac{d}{dr}\Big(\frac{P^{1-\frac{2}{\gamma}}}{r^2}
\frac{d}{dr}(r^2P^{\frac{1}{\gamma}}V^r)\Big)
+\frac{\mathscr{A}}{\rho}\frac{dP}{dr}V^r+ \nonumber \\
&+l(l+1)\frac{\gamma P^{\frac{1}{\gamma}}}{\rho}\frac{d}{dr}
\Big(\frac{P^{1-\frac{1}{\gamma}}}{r}V^h\Big)+\frac{d}{dr}\delta\check{\Phi}, \\
L_l^h&=-\frac{\gamma P^{1-\frac{1}{\gamma}}}{r^3\rho}
\frac{d}{dr}(r^2P^{\frac{1}{\gamma}}V^r)+
\frac{l(l+1)}{r^2}\frac{\gamma P}{\rho}V^h+\frac{1}{r}\delta\check{\Phi},
\end{align}
\end{subequations}
where
\begin{align}
\delta\check{\Phi}&=-4\pi\mathsf{G}\mathcal{H}_l(\delta\check{\rho}) \nonumber \\
&=4\pi\mathsf{G}\mathcal{H}_l
\Big(\frac{1}{r^2}\frac{d}{dr}(r^2\rho V^r)-\frac{l(l+1)}{r}\rho V^h\Big).
\end{align}

Using this expression, we see that the operator $\vec{L}_l$ restricted on
$C_0^{\infty}([0,R[, \mathbb{C}^2)$ is symmetric and bounded from below in $\mathfrak{X}_l$. 

In fact, if
$\vec{V}_{(\mu)} \in C_0^{\infty}, \mu=1,2$, then the integration by parts leads us to
\begin{align}
(\vec{L}_l\vec{V}_{(1)}|\vec{V}_{(2)})_{\mathfrak{X}_l} &=\gamma \int U_{(1)}U_{(2)}^*dr +\int\mathscr{A}\frac{dP}{dr}V_{(1)}^rV_{(2)}^{r*}r^2dr + \nonumber \\
&-4\pi\mathsf{G}\int\mathcal{H}_l(\delta\check{\rho}_{(1)})(\delta\check{\rho}_{(2)})^*r^2dr, \label{Ad.4.26}
\end{align}
where
\begin{equation}
U_{(\mu)}:=\frac{P^{\frac{1}{2}-\frac{1}{\gamma}}}{r}
\frac{d}{dr}(r^2P^{\frac{1}{\gamma}}V_{(\mu)}^r)
-l(l+1)P^{\frac{1}{2}}V_{(\mu)}^h.
\end{equation}
Since the integral operator $\mathcal{H}_l$ is symmetric, we see that the restriction of $\vec{L}_l$ onto
$C_0^{\infty}$ is symmetric. 

Let us estimate $$(\vec{L}_l\vec{V}|\vec{V})_{\mathfrak{X}_l}
=\gamma\int |U|^2dr +\int\mathscr{A}\frac{dP}{dr}|V^r|^2r^2dr
-4\pi\mathsf{G}\int\mathcal{H}_l(\delta\check{\rho})(\delta\check{\rho})^*r^2dr$$ from below
for $\vec{V} \in C_0^{\infty}$. 

Since 
$$\Big|\mathscr{A}\frac{dP}{dr}\Big|\leq C_1 \rho ,$$
we have
$$\Big|\int
\mathscr{A}\frac{dP}{dr}|V^r|^2r^2dr\Big|
\leq C_1\|\vec{V}\|_{\mathfrak{X}_l}^2.
$$
On the other hand, we know
$$\Big|\int\mathcal{H}_l(\delta\check{\rho})(\delta\check{\rho})^*r^2dr\Big|\leq
\rho_{\mathsf{O}}\|\vec{V}\|_{\mathfrak{X}_l}^2.$$
For a proof, see \cite[Section 5.2]{JJTM}.
Therefore we have
\begin{equation}
(\vec{L}_l\vec{V}|\vec{V})_{\mathfrak{X}_l}\geq
\gamma \int |U|^2dr -(C_1+4\pi\mathsf{G}\rho_{\mathsf{O}})\|\vec{V}\|_{\mathfrak{X}_l}^2, \label{Ad.4.28}
\end{equation}
that is, $\vec{L}_l$ is bounded from below. 

Therefore, thanks to \cite[Chapter VI, Section 2.3]{Kato}, the restriction of $\vec{L}_l$ onto
$C_0^{\infty}$ admits the Friedrichs extension. This completes the proof of Theorem \ref{Th.3}.
$\square$\\

We see
\begin{equation}
U=-\frac{1}{\gamma}rP^{-\frac{1}{2}}\delta\check{P}
\end{equation}
by a direct calculation.\\

Hereafter we denote by $\vec{L}_l$ the self-adjoint operator in $\mathfrak{X}_l$. 

Note that the domain $\mathsf{D}(\vec{L}_l)$ of the Friedrichs extension is given by
\begin{equation}
\mathsf{D}(\vec{L}_l)=\overset{\circ}{\mathfrak{W}}_l\cap\{\vec{V} \  | \  \vec{L}_l\vec{V} \in \mathfrak{X}_l\quad\mbox{in distribution sense}\quad \}.
\end{equation}
Here $\overset{\circ}{\mathfrak{W}}_l$ is the closure of
$C_0^{\infty}([0,R[,\mathbb{C}^2)$ in the Hilbert space $\mathfrak{W}_l$  endowed with the norm
$\|\cdot\|_{\mathfrak{W}_l}$ given by
\begin{equation}
\|\vec{V}\|_{\mathfrak{W}_l}^2=
\|\vec{V}\|_{\mathfrak{X}_l}^2+\|
\delta\check{\rho}\|_{L^2(\frac{\gamma P}{\rho^2}r^2dr)}^2, \label{NormWl}
\end{equation}
where 
\begin{equation}
\|\delta\check{\rho}\|_{L^2(\frac{\gamma P}{\rho^2}r^2dr)}^2=
\int_0^R
\Big|-\frac{1}{r^2}\frac{d}{dr}(r^2\rho V^r)+
\frac{l(l+1)}{r}\rho V^h\Big|^2\frac{\gamma P}{\rho^2}r^2dr.
\end{equation}\\

Now, actually \eqref{NormWl} is equivalent to
$$\|\vec{V}\|_{\mathfrak{X}_l}^2+\gamma \int_0^R|U|^2dr $$
with
$$
U=\frac{P^{\frac{1}{2}-\frac{1}{\gamma}}}{r}
\frac{d}{dr}(r^2P^{\frac{1}{\gamma}}V^r)
-l(l+1)P^{\frac{1}{2}}V^h 
=-\frac{1}{\gamma}rP^{-\frac{1}{2}}\delta\check{P},
$$
that is, \eqref{NormWl} is equivalent to $((\vec{L}_l+\kappa)\vec{V}|\vec{V})_{\mathfrak{X}_l}$ with a sufficiently large constant $\kappa$.
In fact, 
$$U=\frac{1}{\gamma}rP^{-\frac{1}{2}}\delta\check{P}=
r\rho^{-1}P^{\frac{1}{2}}\delta\check{\rho} +r \mathscr{A}P^{\frac{1}{2}}V^r $$
implies
$$\frac{1}{2}
\int_0^R|\delta\check{\rho}|^2\frac{P}{\rho^2}r^2dr
-B\leq
\int_0^R|U|^2dr \leq 2
\int_0^R|\delta\check{\rho}|^2\frac{P}{\rho^2}r^2dr
+2B $$
with
$\displaystyle B=\int_0^R|\mathscr{A}V^r|^2Pr^2dr$, while
$B \leq C\|\vec{V}\|_{\mathfrak{X}_l}^2$ since
$|\mathscr{A}|^2P \leq C\rho$. \\

We can claim the following 

\begin{Lemma}
The operator $\vec{L}_l$ can be considered as a self-adjoint operator in the Hilbert space $\mathfrak{W}_l$, which is a proper dense subspace of $\mathfrak{X}_l$.
\end{Lemma}

Proof by a direct calculation using the inner  product $(\cdot | \cdot)_{\mathfrak{W}_l}$ is too complicated. So we introduce another inner product $(\cdot|\cdot)_{\clubsuit}$ defined by
\begin{equation}
(\vec{V}_{(1)}|\vec{V}_{(2)})_{\clubsuit}=
(\vec{L}_l\vec{V}_{(1)}|\vec{V}_{(2)})_{\mathfrak{X}_l}+
\kappa (\vec{V}_{(1)}|\vec{V}_{(2)})_{\mathfrak{X}_l},
\end{equation}
where $\kappa$ is a sufficiently large positive number. Thanks to
\eqref{Ad.4.26},\eqref{Ad.4.28}, we can assume that the corresponding norm
$\|\cdot\|_{\clubsuit}$, which can be defined by
$$\|\vec{V}\|_{\clubsuit}^2=(\vec{V}|\vec{V})_{\clubsuit}, $$
is equivalent to the norm $\|\cdot\|_{\mathfrak{W}_l} $, that is,
$$\frac{1}{C}\|\vec{V}\|_{\mathfrak{W}_l} \leq
\|\vec{V}\|_{\clubsuit} \leq C\|\vec{V}\|_{\mathfrak{W}_l}. $$
As shown in \eqref{Ad.4.26}, $\vec{L}_l$ resricted to $C_0^{\infty}$ is symmetric with respect to this inner product $(\cdot|\cdot)_{\clubsuit}$, and is bounded from below as shown in \eqref{Ad.4.28}. Thus the Friedrichs extension is a self-adjoint operator in
$\mathfrak{W}_l$. $\square$.\\

Here let us note that we can claim

\begin{Proposition}\label{Prop.2}
Let $\vec{V}\in \mathfrak{W}_l$. Then $\vec{V}$
belongs to
$\overset{\circ}{\mathfrak{W}}_l$ if $V^r =O((R-r)^{\kappa})$ with
$\displaystyle \kappa > -\frac{\nu}{2} $ for $0 < R-r \ll 1$.
\end{Proposition}

Proof can be done by taking $\varphi_n\in C^{\infty}([0,R])$
such that $\varphi_n(r)=1$ for $0 \leq r \leq R-\frac{1}{n},
\varphi_n(r)=0$ for  $R-\frac{1}{2n} \leq r \leq R$, and
$0 \leq \varphi_n\leq 1, |d\varphi_n/dr|\leq Cn$, and considering 
$\varphi_n\cdot \vec{V}$ for $\vec{V} \in \mathfrak{W}_l$. Let us omit the details.\\

As for the dimension of the kernel of $\vec{L}_l$, we have the following

\begin{Theorem}
1) Let $l \geq1$. Suppose that $\mathscr{A}=0$ identically on $]0,R[$. Then
$\mathrm{dim}\mathsf{N}(\vec{L}_l)=\infty$.

2) Suppose that $\mathscr{A}\not=0$ almost  everywhere on $]0,R[$. Then $\mathrm{dim}\mathsf{N}(\vec{L}_l)=0$ when $l \geq 2$ and $\mathrm{dim}\mathsf{N}(\vec{L}_1)=1$ when $l=1$.
\end{Theorem}

Proof. 1) Let $l \geq 1$ and suppose that $\mathscr{A}=0$ identically. 
Then, if 
\begin{equation}
-\frac{1}{r^2}\frac{d}{dr}(r^2\rho V^r)+\frac{l(l+1)}{r}\rho V^h=0, \label{4.29}
\end{equation}
then $\delta\check{\rho}=\delta\check{P}=\delta\check{\Phi}=0$ by
\eqref{4.6a}, \eqref{4.6b}, \eqref{4.6c}, therefore $L_l^r=L_l^h=0$, that
is, $\vec{V}\in \mathsf{N}(\vec{L}_l)$. But, 
if $a \in C_0^{\infty}(]0,R[)$, then $\vec{V}=(V^r,V^h)^T$ given by
$$V^r(r)=\frac{1}{\rho}\cdot\frac{l(l+1)}{r^2}\int_0^ra(\acute{r})\acute{r}d\acute{r},\quad
V^h(r)=\frac{1}{\rho}a(r)
$$
belongs to $C_0^{\infty}(]0,R[)$ and satisfies \eqref{4.29} so that $\vec{L}_l\vec{V}=0$. Since $a \in C_0^{\infty}$ is arbitrary, we see $\mathrm{dim}\mathsf{N}(\vec{L}_l)=\infty$.

2) Suppose that $\mathscr{A} \not=0$ almost everywhere. Let us consider $\vec{V}\in \mathsf{N}(\vec{L}_l)$, $l\geq 1$. Of course $\vec{V} \in \mathfrak{X}_l$, and moreover, since $\vec{V} \in
\mathsf{D}(\vec{L}_l)$, we have
$$\delta\check{\rho} \in L^2([0,R], P\rho^{-2}r^2dr) \subset L^2([0,R], r^2dr), $$
since $P\rho^{-2}\geq 1/C$. Thus $\delta\check{\Phi}=-4\pi\mathsf{G}\mathcal{H}_l(\delta\check{\rho})$ enjoys the properties listed in Lemma \ref{Lem.2}.
Now $\vec{V}\in \mathsf{N}(\vec{L}_l)$ means
\begin{subequations}
\begin{align}
&\frac{1}{\rho}\frac{d}{dr}\delta\check{P}-\frac{1}{\rho^2}\frac{dP}{dr}\check{\rho}+
\frac{d}{dr}\delta\check{\Phi}=0, \label{N4.33a}\\
&\frac{\delta\check{P}}{\rho}+\delta\check{\Phi}=0. \label{N4.33b}
\end{align}
\end{subequations}
It follows from \eqref{N4.33a}\eqref{N4.33b} that
\begin{equation}
\Phi'\delta\check{\rho}=\frac{d\rho}{dr}\delta\check{\Phi}, \label{4.31}
\end{equation}
where
 we put
\begin{equation}
\Phi':=\frac{d\Phi}{dr}=-\frac{1}{\rho}\frac{dP}{dr}. 
\end{equation}
On the other hand, $\delta\check{\Phi}=-4\pi\mathsf{G}
\mathcal{H}_l(\delta\check{\rho})$ implies
\begin{equation}
\Big[\frac{1}{r^2}\frac{d}{dr}r^2\frac{d}{dr}-\frac{l(l+1)}{r^2}\Big]\delta\check{\Phi}=
4\pi\mathsf{G}\delta\check{\rho}. \label{4.33}
\end{equation}
Let us consider
\begin{equation}
Z:=\frac{\delta\check{\Phi}}{\Phi'},
\end{equation}
keeping in mind that $$\Phi'
=\frac{d\Phi}{dr}=\frac{4\pi\mathsf{G}}{r^2}\int_0^r\rho(\acute{r})\acute{r}^2d\acute{r}>0
\quad\mbox{for}\quad r>0.$$ 
Then $\delta\check{\Phi}=\Phi'Z$ and \eqref{4.31} reads
 $\displaystyle \delta\check{\rho}=\frac{d\rho}{dr}Z$. Thus, eliminating $\delta\check{\rho}$ from
\eqref{4.31} and \eqref{4.33}, we can derive the equation
\begin{equation}
-\frac{d}{dr}\Big(r^2(\Phi')^2\frac{dZ}{dr}\Big)+(l+2)(l-1)(\Phi')^2Z=0, \label{4.35}
\end{equation}
using $\displaystyle \frac{1}{r^2}\frac{d}{dr}\Big(r^2\Phi'\Big)=4\pi\mathsf{G}\rho$.
We owe this trick to N. R. Lebovitz, \cite{Lebovitz}, but we are considering that the equation 
\eqref{4.35} holds for $0<r<+\infty$ in view of $\displaystyle \delta\check{\rho}=\frac{d\rho}{dr}=0$ for $r \geq R$ so that we do not care the `boundary condition' of $\delta\check{\Phi}$ at $r=R$.

 Let us multiply \eqref{4.35} by $Z$ and integrate it on $[0,+\infty[$. We are going to perform the integration by parts using the following observations:
For $r \geq R$ we have
\begin{align*}
&\Phi'=\frac{C_0}{r^2},\quad C_0=\mathrm{Const.}>0, \\
&\delta\check{\Phi}=\frac{C_1}{r^{l+1}},\quad C_1=\mathrm{Const.}, \\
&Z=\frac{C_1}{C_0}\frac{1}{r^{l-1}},
\end{align*}
and, on the other hand, as $r\rightarrow +0$, we have
\begin{align*}
&\Phi'=O(r), \quad\frac{d\Phi'}{dr}=\frac{d^2\Phi}{dr^2}=O(1), \\
&\delta\check{\Phi}=O(r^{\frac{1}{2}}), \quad \frac{d}{dr}\delta\check{\Phi}=O(r^{-\frac{1}{2}}), \\
\mbox{therefore}& \\
&\frac{dZ}{dr}=O(r^{-\frac{3}{2}}),
\quad
r^2(\Phi')^2\frac{dZ}{dr}Z=O(r^2).
\end{align*}
Thus the contributions from the boundaries vanish, that is,
$$ r^2(\Phi')^2\frac{dZ}{dr}Z \rightarrow 0 $$
both as $r \rightarrow +0$ and as $r \rightarrow +\infty$. So the integration by parts gives
\begin{equation}
\int_0^{+\infty}r^2(\Phi')^2\Big(\frac{dZ}{dr}\Big)^2dr+
(l+2)(l-1)\int_0^{+\infty}(\Phi')^2Z^2dr=0. \label{4.36}
\end{equation}

Let $l \geq 2$, that is, $(l+2)(l-1)>0$. Then 
\eqref{4.36} implies $Z=0$, therefore, 
$\delta\check{\Phi}=\delta\check{\rho}=\delta\check{P}=0$. Then we have
$$ \gamma\mathscr{A}PV^r=\delta\check{P}-\frac{\gamma P}{\rho}\delta\check{\rho}=0$$
by \eqref{4.6c}. Since $\mathscr{A}\not=0$ almost everywhere, we have $V^r=0$. Since
$\delta\check{\rho}=0$, this implies $V^h=0$ in view of \eqref{4.6a}. Thus $\vec{V}=\vec{0}$, and $\mathrm{dim}\mathsf{N}(\vec{L}_l)=0$.

Let $l=1$. Then $\displaystyle \frac{dZ}{dr}=0$, that is, $Z$ is a constant $\kappa$. Then we have
$$\delta\check{\Phi}=-\kappa\frac{1}{\rho}\frac{dP}{dr},\quad
\delta\check{P}=\kappa\frac{dP}{dr}, \quad
\delta\check{\rho}=\kappa\frac{d\rho}{dr},
$$
which imply $V^r=V^h=-\kappa$. Conversely, if $V^r=V^h=1$, then
we have
$$\delta\check{\rho}=-\frac{d\rho}{dr}, \quad
\delta\check{P}=-\frac{dP}{dr}, \quad
\delta\check{\Phi}=4\pi\mathsf{G}\mathcal{H}_1\Big(-\frac{d\rho}{dr}\Big)=-\frac{d\Phi}{dr},
$$
and $L_1^r=L_1^h=0$. But $(1,1)^{\top} \in \mathsf{D}(\vec{L}_1)$, 
since $(1,1)^{\top} \in \mathfrak{W}_l$.
Recall Proposition \ref{Prop.2}.

Summing up, we can claim
$\mathrm{dim}\mathsf{N}(\vec{L}_1)=1$. This completes the proof.
$\square$\\

As noted in \cite{Lebovitz}, the identity
$$Y_{10}\mbox{\boldmath$e$}_r+\nabla_sY_{10}=\sqrt{\frac{3}{4\pi}}\frac{\partial}{\partial x^3} $$
leads us to the interpretation that the eigenfunction $\vec{V}=(1,1)^{\top}$ for $l=1$ means a uniform translation, and it can be eliminated by requiring that the center of mass remains fixed in space, or, by requiring $\displaystyle \int_{B_R}(\delta\rho)x^3d\mbox{\boldmath$x$}=0$.

\section{Existence of g-modes an p-modes}

We are considering the eigenvalue problem
\begin{equation}
\vec{L}_l\vec{V}=\lambda \vec{V}. \label{G01}
\end{equation}

In this Section we suppose $l \geq 1$ and consider $\lambda >0$. We are going to prove the existence of a sequence $(\lambda_{-n})_{n\in\mathbb{N}}$ of positive eigenvalues such that $\lambda_{-n} \rightarrow 0$ as $ n \rightarrow \infty$, so called `g-modes',
and a sequence $(\lambda_n)_{n\in\mathbb{N}}$ of positive eigenvalues such that $\lambda_n \rightarrow +\infty$ as $ n\rightarrow \infty$, so called `p-modes'. 
In order to do it, the formulation of the equations given by D. O. Gough, \cite{Gough}, is useful. Let us adopt it as follows.\\

First of all let us recall the eigenvalue problem \eqref{G01} is
\begin{subequations}
\begin{align}
&\frac{1}{\rho}\frac{d}{dr}\delta \check{P}-\frac{1}{\rho^2}\frac{dP}{dr}\delta\check{\rho}+\frac{d}{dr}\delta\check{\Phi}=\lambda V^r, \\
&\frac{1}{r}\Big(\frac{\delta\check{P}}{\rho}+\delta\check{\Phi}\Big)=\lambda V^h.
\end{align}
\end{subequations}
We are taking the abbreviation  $\rho, P, S, \Phi$ for
$\bar{\rho}, \bar{P}, \bar{S}, \bar{\Phi}$. We have
\begin{subequations}
\begin{align}
&\delta\check{\rho}=-\frac{1}{r^2}\frac{d}{dr}(r^2\rho V^r)+\frac{l(l+1)}{r}\rho V^h, \\
&\delta\check{P}=\mathsf{c}^2\delta\check{\rho}+\gamma\mathscr{A}PV^r
=\mathsf{c}^2\delta\check{\rho}+\Big(\mathsf{c}^2\frac{d{\rho}}{dr}+g{\rho}\Big)V^r, \\
&\delta\check{\Phi}=-4\pi\mathsf{G}\mathcal{H}_l(\delta\check{\rho}).
\end{align}
\end{subequations}
Here we put
\begin{equation}
\mathsf{c}^2=\frac{\gamma P}{\rho},\qquad g=-\frac{1}{\rho}\frac{dP}{dr}=\frac{d\Phi}{dr}.
\end{equation}
The integral operator $\mathcal{H}_l$ is defined by
\begin{equation}
\mathcal{H}_lf(r)=\frac{1}{2l+1}\Big[
\int_r^{R}f(\acute{r})\Big(\frac{r}{\acute{r}}\Big)^l\acute{r}d\acute{r}+
\int_0^rf(\acute{r})\Big(\frac{r}{\acute{r}}\Big)^{-(l+1)}\acute{r}d\acute{r}\Big],
\end{equation}

Following \cite{Gough}, we introduce the variables
\begin{equation}
\xi=V^r,\qquad \eta=\delta\check{P}+\frac{dP}{dr}V^r
=\delta\check{P}-\rho g V^r. \label{Defxieta}
\end{equation}
We note that, in the context of the  linearized approximation, $\eta$ is nothing but $\Delta\check{P}$, say, the Lagrange perturbation. Then, supposing that $\lambda >0$, we have the system of equations, which is equivalent to \eqref{G01}, 
\begin{subequations}
\begin{align}
&\frac{d\xi}{dr}+A_{11}\xi+A_{12}\eta+A_{10}=0, \label{G11a}\\
&\frac{d\eta}{dr}+A_{21}\xi+A_{22}\eta+A_{20}=0, \label{G11b}
\end{align}
\end{subequations}
with
\begin{align}
&A_{11}=\frac{2}{r}-\frac{l(l+1)g}{\lambda r^2},\quad A_{12}=\Big(1-\frac{l(l+1)\mathsf{c}^2}{\lambda r^2}\Big)\frac{1}{\mathsf{c}^2\rho}, \nonumber \\
&A_{21}=\rho\Big[
\frac{l(l+1)}{\lambda}\Big(\frac{g}{r}\Big)^2
-g\Big(\frac{1}{\mathsf{H}[g]}+\frac{2}{r}\Big)-\lambda\Big],
\quad
A_{22}=\frac{l(l+1)g}{\lambda r^2}, \nonumber \\
&A_{10}=-\frac{l(l+1)}{\lambda r^2}\delta\check{\Phi}, \nonumber \\
&A_{20}=\rho\Big[\frac{d}{dr}\delta\check{\Phi}+
\frac{l(l+1)g}{\lambda r^2}\delta\check{\Phi}\Big]. \label{Aij}
\end{align}

Here and hereafter we use

\begin{Definition}
For any quantity $Q$ which is a function of $r$, we put
\begin{equation}
\frac{1}{\mathsf{H}[Q]}=-\frac{d}{dr}\log Q=-\frac{1}{Q}\frac{dQ}{dr}.
\end{equation}
\end{Definition}
$\mathsf{H}[Q]$ is the so called `scale height' of $Q$. \\

Since we have 
\begin{equation}
\delta\check{\rho}=\frac{\rho}{\mathsf{H}[\rho]}\xi+\frac{1}{\mathsf{c}^2}\eta, \label{G15}
\end{equation}
we can write
\begin{subequations}
\begin{align}
&A_{10}=4\pi\mathsf{G}\frac{l(l+1)}{\lambda r^2}\mathcal{H}_l
\Big(\frac{\rho}{\mathsf{H}[\rho]}\xi+\frac{1}{\mathsf{c}^2}\eta\Big), \\
&A_{20}=-4\pi\mathsf{G}\rho\Big[
\frac{1}{r}\dot{\mathcal{H}}_l
\Big(\frac{\rho}{\mathsf{H}[\rho]}\xi+\frac{1}{\mathsf{c}^2}\eta\Big)
+\frac{l(l+1)g}{\lambda r^2}\mathcal{H}_l
\Big(\frac{\rho}{\mathsf{H}[\rho]}\xi+\frac{1}{\mathsf{c}^2}\eta\Big)
\Big],
\end{align}
\end{subequations}
where
\begin{align}
\dot{\mathcal{H}}_lf&:=r\frac{d}{dr}\mathcal{H}_lf=  \nonumber \\
&=\frac{1}{2l+1}
\Big[
l\int_r^Rf(\acute{r})\Big(\frac{r}{\acute{r}}\Big)^l\acute{r}d\acute{r}-
(l+1)\int_0^rf(\acute{r})\Big(\frac{r}{\acute{r}}\Big)^{-(l+1)}\acute{r}d\acute{r}
\Big].
\end{align}\\

We are going to analyze the system of equations \eqref{G11a}, \eqref{G11b}. \\

Keeping in mind that it should hold
\begin{align*}
& \mathsf{c}^2\delta\check{\rho}+\mathsf{c}^2\frac{d\rho}{dr}V^r=\eta, \\
&l(l+1)V^h=r\frac{1}{\rho\mathsf{c}^2}\eta+\Big(r\frac{d}{dr}+2\Big)\xi,
\end{align*}
we can claim the following

\begin{Proposition} \label{EqVxieta}
We consider the correspondence between the variable $\vec{V}=(V^r,V^h)$ and $(\xi, \eta)$ defined by \eqref{Defxieta}. Then 1) $\vec{V} \in \mathfrak{W}_l$ if and only if
$\displaystyle \xi \in L^2(\rho r^2dr), \eta \in L^2\Big(\frac{1}{\mathsf{c}^2\rho}r^2dr\Big)$
and
\begin{equation}
r\frac{1}{\mathsf{c}^2\rho}\eta+\Big(r\frac{d}{dr}+2\Big)\xi \quad\in\quad 
L^2(\rho r^2dr). \label{AX}
\end{equation}
2) When $(\xi, \eta )$ satisfies the equation \eqref{G11a}, then \eqref{AX} reduces
\begin{equation}
\frac{g}{r}\xi + \frac{1}{r\rho}\eta+\frac{1}{r}\delta\check{\Phi} \quad\in\quad L^2(\rho r^2dr). \label{AX2}
\end{equation}
\end{Proposition}

We suppose the following assumptions:

\begin{Assumption}\label{Ass.G2}
$\bar{\rho}$ and $\bar{S}$ are analytic functions of $r^2$ near $r=0$.
\end{Assumption}

\begin{Assumption}\label{Ass.G3}
As $r \rightarrow R-0$, it holds that
\begin{align}
&\bar{\rho}=C_{\rho}(R-r)^{\nu}(1+[R-r, (R-r)^{\nu+1}]_1), \\
&\bar{S}=[R-r, (R-r)^{\nu+1}]_0,
\end{align}
where $C_{\rho}$ is a positive constant and
\begin{equation}
\nu:=\frac{1}{\gamma-1}.
\end{equation}
Here  and hereafter we use the following notation.
\end{Assumption}

\begin{Notation}
Fir a non-negative integer $K$, 
the symbol $[X]_K$ stands for various convergent power series of the form
$\displaystyle \sum_{k\geq K}a_kX^k$, and
the symbol $[X_1, X_2]_N$ stands for various convergent double power series of the form 
$ \displaystyle \sum_{k_1+k_2 \geq N}a_{k_1k_2}X_1^{k_1}X_2^{k_2} $.
\end{Notation}

Assumption \ref{Ass.G3} is satisfied if $\bar{S}=\Sigma(\bar{\rho}^{\gamma-1})$ with a function $\Sigma$ which is analytic at $\rho_{\mathsf{O}}^{\gamma-1}$, and 
Assumption \ref{Ass.G3} is satisfied if the function $\Sigma$ is analytic at $0$. For a
proof see Appendix. \\

\subsection{ g-modes}

We suppose 
\begin{Assumption}\label{Ass.G1}
There exists a positive number $\delta$ such that 
\begin{equation}
\frac{1}{r}\frac{d\bar{S}}{dr}\geq \delta \label{G02}
\end{equation}
for $0 < r< R$. Or, equivalently, there exists a positive number $\delta_{\mathbf{A}}$ such that 
$$\mathscr{N}^2\Big(=-g\mathscr{A}=\frac{g}{\gamma \mathsf{C}_V}\frac{d\bar{S}}{dr} \Big)
\geq \delta_{\mathbf{A}} r^2$$
uniformly on $0<r<R$. 
\end{Assumption}

If the equilibrium is given through $S=\Sigma(\rho^{\gamma-1})$ by Theorem \ref{Th.1}, 
this assumption means
that there is a positive number $\delta_{\Sigma}$ such that
\begin{equation}
\frac{d}{d\omega}\Sigma(\omega) \leq -\delta_{\Sigma}
\end{equation}
for $0< \omega \leq \rho_{\mathsf{O}}^{\gamma-1}$.\\

Keeping in mind small $\lambda$, we  put
\begin{subequations}
\begin{align}
&E=1+(E_1+E_2\lambda)\frac{\lambda}{l(l+1)}, \label{DefEa}\\
&E_1=-\Big(\frac{r}{\mathsf{H}[g]}+2\Big)\frac{r}{g}=
-4\frac{r}{g}+
4\pi\mathsf{G}\rho \Big(\frac{r}{g}\Big)^2, \label{DefEb}\\
&E_2=-\Big(\frac{r}{g}\Big)^2, \label{DefEc}
\end{align}
\end{subequations}
to write
\begin{equation}
A_{21}=\frac{l(l+1)\rho g^2}{\lambda r^2}E
\end{equation}
in \eqref{Aij}.

Eliminating $\xi$ from the system \eqref{G11a}\eqref{G11b}, we get the single second order equation for $\eta$:
\begin{align}
&\frac{d^2\eta}{dr^2}+\frac{1}{\mathsf{H}[\mathfrak{Q}]}\frac{d\eta}{dr} + \nonumber \\
&+\Big[
\frac{1}{\mathsf{c}^2}\Big(\lambda+g\Big(\frac{1}{\mathsf{H}[g]}+\frac{2}{r}\Big)\Big)
-\frac{l(l+1)}{r^2}\Big(1-\frac{\mathfrak{N}^2}{\lambda}\Big)
\Big]\eta
+\rho F=0, \label{G18}
\end{align}
where
\begin{align}
&F=\frac{1}{\rho}\Big(\frac{d}{dr}+\frac{1}{\mathsf{H}[\mathfrak{Q}]}-\frac{\mathscr{S}_l^2g}{\mathsf{c}^2\lambda}\Big)
\rho\Big(\frac{d}{dr}+\frac{\mathscr{S}_l^2g}{\mathsf{c}^2\lambda}\Big)\delta\check{\Phi}+
\Big(\frac{\mathscr{S}_l^2g}{\mathsf{c}^2\lambda}\Big)^2 E\delta\check{\Phi}, \\
& \mathscr{S}_l^2=\frac{l(l+1)\mathsf{c}^2}{r^2}, \\
&\mathfrak{Q}=\frac{\rho g^2}{r^4}E, \\
&\mathfrak{N}^2=g\Big(\frac{1}{\mathsf{H}[\mathfrak{Q}]}-
\frac{g}{\mathsf{c}^2}-\frac{2}{\mathsf{H}[g]}-\frac{4}{r}\Big)=
\mathscr{N}^2+\frac{g}{\mathsf{H}[E]}. \label{G22}
\end{align}\\

{\it Here and hereafter we consider $\lambda \in [0,\lambda_0], \mathsf{G}\in [0,\mathsf{G}_0]$, where $\lambda_0$ is a fixed small positive number and $\mathsf{G}_0$ is a fixed  large number such that $\mathsf{G}_0 >1$. }\\

Taking $\lambda_0$ sufficiently small, we can suppose that $E \geq 1-\varepsilon >0$ for $0<r<R, 0\leq \lambda \leq \lambda_0$, with $0< \varepsilon \ll 1$, since $E_1, E_2$ are bounded on $0<r<R$. \\

Calculating $F$, we reduce \eqref{G18} to
\begin{equation}
\frac{d^2\eta}{dr^2}+A\frac{d\eta}{dr}+B\eta+\rho K=0, \label{G24}
\end{equation}
where
\begin{subequations}
\begin{align}
A&=\frac{1}{\mathsf{H}[\mathfrak{Q}]}
-4\pi\mathsf{G}\frac{\rho }{\mathsf{H}[\rho]E}\frac{\mathsf{c}^2\lambda}{\mathscr{S}_l^2g^2} \nonumber \\
&=\frac{1}{\mathsf{H}[\rho g^2/r^4]}+\frac{1}{\mathsf{H}[E]}
-4\pi\mathsf{G}\frac{\rho }{\mathsf{H}[\rho]E}\frac{\mathsf{c}^2\lambda}{\mathscr{S}_l^2g^2}  \label{G24a}\\
B&=\frac{1}{\mathsf{c}^2}\Big(\lambda +g
\Big(\frac{1}{\mathsf{H}[g]}+\frac{2}{r}\Big)\Big)-
\frac{l(l+1)}{r^2}
\Big(1-\frac{\mathfrak{N}^2}{\lambda}\Big)+ \nonumber \\
&+4\pi\mathsf{G}\rho\Big(-\frac{1}{\mathsf{H}[\rho]gE}+\frac{1}{\mathsf{c}^2}\Big),  \label{G24b}\\
K&=\Big(\frac{1}{\mathsf{H}[g^2E/r^2]}
-4\pi\mathsf{G}\frac{\rho }{\mathsf{H}[\rho]E}\frac{\mathsf{c}^2\lambda}{\mathscr{S}_l^2g^2}\Big)\frac{d}{dr}\delta\check{\Phi} + \nonumber \\
&+\Big(
\frac{1}{\mathsf{H}[E]}\frac{\mathscr{S}_l^2g}{\mathsf{c}^2\lambda}
-4\pi\mathsf{G}
\frac{\rho}{\mathsf{H}[\rho]gE}\Big)\delta\check{\Phi}. \label{G24c}
\end{align}
\end{subequations}

Here we have used the equation
\begin{equation}
\Big[\frac{d^2}{dr^2}+\frac{2}{r}\frac{d}{dr}-\frac{l(l+1)}{r^2}\Big]\delta\check{\Phi}=4\pi\mathsf{G}\delta\check{\rho}
\end{equation}
and the identity
\begin{equation}
\frac{\lambda}{g}+
\frac{1}{\mathsf{H}[gE/r^2]}+
(E-1)\frac{\mathscr{S}_l^2g}{\mathsf{c}^2\lambda}=\frac{1}{\mathsf{H}[E]}.
\end{equation}
Note that
$$
\frac{1}{\mathsf{H}[E]}\frac{\mathscr{S}_l^2g}{\mathsf{c}^2\lambda}=-\frac{1}{E}\Big(\frac{dE_1}{dr}+\lambda\frac{dE_2}{dr}\Big)\frac{g}{r^2}.
$$

Moreover \eqref{G15} reads
\begin{align}
\delta\check{\rho}&=-\frac{1}{\mathsf{H}[\rho]E}\frac{\mathsf{c}^2\lambda}{\mathscr{S}_l^2g^2}\frac{d\eta}{dr}+
\Big(-\frac{1}{\mathsf{H}[\rho]gE}+\frac{1}{\mathsf{c}^2}\Big)\eta + \nonumber \\
&-\frac{\rho }{\mathsf{H}[\rho]E}\frac{\mathsf{c}^2\lambda}{\mathscr{S}_l^2g^2}
\frac{d}{dr}\delta\check{\Phi}-\frac{\rho}{\mathsf{H}[\rho]gE}
\delta\check{\Phi},
\end{align}
which should be inserted to
\begin{subequations}
\begin{align}
\delta\check{\Phi}&=-4\pi\mathsf{G}\mathcal{H}_l\delta\check{\rho}, \\
r\frac{d}{dr}\delta\check{\Phi}&=-4\pi\mathsf{G}
\dot{\mathcal{H}}_l\delta\check{\rho}.
\end{align}
\end{subequations}

This means that the variables
\begin{equation}
X:=\delta\check{\Phi}, \qquad \dot{X}:=r\frac{d}{dr}\delta\check{\Phi} \label{G29}
\end{equation}
should satisfy the system of equations
\begin{subequations}
\begin{align}
X&=4\pi\mathsf{G}\mathcal{H}_l
\frac{\rho}{\mathsf{H}[\rho]gE}X
+4\pi\mathsf{G}\mathcal{H}_l
\frac{\rho }{\mathsf{H}[\rho]Er}\frac{\mathsf{c}^2\lambda}{\mathscr{S}_l^2g^2}\dot{X} 
+ \nonumber \\
&+4\pi\mathsf{G}\mathcal{H}_lY, \label{G30a}\\
\dot{X}&=4\pi\mathsf{G}\dot{\mathcal{H}}_l
\frac{\rho}{\mathsf{H}[\rho]gE}X
+4\pi\mathsf{G}\dot{\mathcal{H}}_l
\frac{\rho }{\mathsf{H}[\rho]Er}\frac{\mathsf{c}^2\lambda}{\mathscr{S}_l^2g^2}\dot{X} 
+ \nonumber \\
&+4\pi\mathsf{G}\dot{\mathcal{H}}_lY, \label{G30b}
\end{align}
\end{subequations}
where
\begin{equation}
Y:=\frac{1}{\mathsf{H}[\rho]Er}\frac{\mathsf{c}^2\lambda}{\mathscr{S}_l^2g^2}\cdot r\frac{d\eta}{dr}+
\Big(\frac{1}{\mathsf{H}[\rho]gE}-\frac{1}{\mathsf{c}^2}\Big)\eta. \label{G31}
\end{equation}

Here let us recall that under Assumptions \ref{Ass.G1}, \ref{Ass.G2}, \ref{Ass.G3}, we have the following asymptotic behaviors:

$$
\mathsf{c}^2=\frac{\gamma P}{\rho}=\begin{cases}
\displaystyle \frac{\gamma P_{\mathsf{O}}}{\rho_{\mathsf{O}}}(1+[r^2]_1)\quad\mbox{as}\quad r \rightarrow +0 \\
\\
\mathsf{c}_{R}^2(R-r)(1+[R-r, (R-r)^{\nu+1}]_1)\quad\mbox{as}\quad r \rightarrow R^0
\end{cases}
$$
with
$\displaystyle \mathsf{c}_R^2:=\gamma C_{\rho}^{\gamma-1}\exp\frac{S(R)}{\mathsf{C}_V}$; 
$$g=-\frac{1}{\rho}\frac{dP}{dr}=\begin{cases}
\displaystyle \frac{P_{\mathsf{O}1}}{\rho_{\mathsf{O}}}r(1+[r^2]_1) \quad\mbox{as}\quad r\rightarrow +0 \\
\\
g_{R}+[R-r, (R-r)^{\nu+1}]_1 \quad\mbox{as}\quad r\rightarrow R-0
\end{cases}
$$
with 
$g_R:=\nu \mathsf{c}_R^2$ ; 
$$\mathscr{N}^2=\frac{g}{\gamma \mathsf{C}_V}\frac{dS}{dr}
=\begin{cases}
\frac{1}{2}\mathscr{N}_{\mathsf{O}1}^2r^2+[r^2]_2\quad\mbox{as}\quad r\rightarrow +0 \\
\\
\mathscr{N}_R^2+[R-r, (R-r)^{\nu+1}]_1
\quad\mbox{as}\quad  r\rightarrow R-0,
\end{cases}
$$
with constants $\mathscr{N}_{\mathsf{O}1}^2, \mathscr{N}_R^2$ which are positive thanks to the Assumption \ref{Ass.G1},
and consequently,
\begin{align*}
&E_1=\begin{cases}
[r^2]_0 \quad\mbox{as}\quad r \rightarrow +0 \\
\\
\displaystyle -\frac{4R}{g_R}+[R-r, (R-r)^{\nu+1}]_1 + \\
+\displaystyle (R-r)^{\nu}\Big(4\pi\mathsf{G}C_{\rho}  \Big(\frac{R}{g_R}\Big)^2 + [R-r, (R-r)^{\nu+1}]_1\Big) 
\quad\mbox{as}\quad r\rightarrow R-0
\end{cases} \\
& E_2=\begin{cases}
[r^2]_0 \quad\mbox{as}\quad r \rightarrow +0 \\
\\
\displaystyle -\Big(\frac{R}{g_R}\Big)^2+[R-r, (R-r)^{\nu+1}]_1 \quad\mbox{as}\quad
r\rightarrow R-0.
\end{cases}
\end{align*}

Here we note that we can suppose
\begin{equation}
\mathfrak{N}^2=\mathscr{N}^2\Big[
1+\frac{1}{\mathscr{N}^2}\frac{g}{\mathsf{H}[E]}\Big] \geq 
(1-\varepsilon)\mathscr{N}^2 \geq \frac{1}{C}r^2 ,
\end{equation}
for $\lambda \in [0,\lambda_0]$, thanks to Assumption \ref{Ass.G1},
provided that $\lambda_0$ is sufficiently small, and $0<\varepsilon \ll 1$,
since
$E-1, \displaystyle \frac{dE}{dr} =O(\lambda)r$ so that $\displaystyle \frac{g}{\mathsf{H}[E]} =O(\lambda)r^2$.
So we see
$$
\mathfrak{N}^2=
\begin{cases}
\displaystyle \frac{\mathfrak{N}_{\mathsf{O}1}^2}{2}r^2(1+[r^2]_1)\quad\mbox{as}\quad r\rightarrow +0 \\
\\
\mathfrak{N}_R^2+[R-r, (R-r)^{\nu+1}]_1\quad\mbox{as}\quad r \rightarrow R-0,
\end{cases}
$$
with positive numbers $\mathfrak{N}_{\mathsf{O}1}^2$ and $\mathfrak{N}_R^2$.\\

Anyway we have
\begin{align*}
&\frac{\rho}{\mathsf{H}[\rho]gE}=-\frac{d\rho}{dr}\frac{1}{gE}=
\begin{cases}
\displaystyle \frac{ \rho_{\mathsf{O}1}\rho_{\mathsf{O}} }{ P_{\mathsf{O}1} }
(1+[r^2]_1)\quad\mbox{as}\quad r\rightarrow +0 \\
\\
\displaystyle \frac{\nu}{g_R}(R-r)^{\nu-1}
(1+[R-r, (R-r)^{\nu+1}]_1 + \\
+(R-r)^{\nu}[R-r,(r-r)^{\nu+1}]_0)\quad\mbox{as}\quad r \rightarrow R-0, 
\end{cases}\\
&\frac{g}{r}\frac{\mathsf{c}^2\lambda}{\mathscr{S}_l^2g^2}=\frac{\lambda}{l(l+1)}\frac{r}{g} =
\begin{cases}
\displaystyle \frac{\lambda}{l(l+1)}\frac{ \rho_{\mathsf{O}} }{ P_{\mathsf{O}1} }
(1+[r^2]_1)\quad\mbox{as}\quad r \rightarrow +0 \\
\\
\displaystyle \frac{\lambda}{l(l+1)}\frac{R}{g_R}
(1+[R-r, (R-r)^{\nu+1}]_1)\quad\mbox{as}\quad r \rightarrow R-0.
\end{cases}
\end{align*}

Thus we can claim

\begin{Proposition}
The operators 
$\displaystyle \mathcal{H}_l
\frac{\rho}{\mathsf{H}[\rho]gE},
\mathcal{H}_l
\frac{\rho }{\mathsf{H}[\rho]Er}\frac{\mathsf{c}^2\lambda}{\mathscr{S}_l^2g^2},
\dot{\mathcal{H}}_l
\frac{\rho}{\mathsf{H}[\rho]gE},
\dot{\mathcal{H}}_l
\frac{\rho }{\mathsf{H}[\rho]Er}\frac{\mathsf{c}^2\lambda}{\mathscr{S}_l^2g^2}
$
are bounded linear operator from $L^1([0,R]; rdr)$ to $L^{\infty}([0,R])$ and Lipschitz continuous in $\lambda \in [0,\lambda_0]$, $\lambda_0$ being a fixed small positive number.
\end{Proposition}
-- \\

In order to solve the system of equations \eqref{G30a}, \eqref{G30b}, we suppose the following:\\

{\bf (G) :}\  {\it There is a sufficiently small positive number $\delta_{\mathbf{G}}$ such that }
\begin{equation}
\Big|4\pi\mathsf{G} \frac{\rho}{\mathsf{H}[\rho]gE} \Big|
\leq \delta_{\mathbf{G}}
, \quad
\Big|
-4\pi\mathsf{G}
\frac{\rho }{\mathsf{H}[\rho]Er}\frac{\mathsf{c}^2\lambda}{\mathscr{S}_l^2g^2}\Big|
\leq \delta_{\mathbf{G}}.
\end{equation}

Then 
the system of equations
\eqref{G30a}, \eqref{G30b} is uniquely solved as
\begin{equation}
\begin{bmatrix}
X \\
\\
\dot{X}
\end{bmatrix}
=
4\pi\mathsf{G}
\begin{bmatrix}
\displaystyle 1-4\pi\mathsf{G} \mathcal{H}_l\frac{\rho}{\mathsf{H}[\rho]gE} &
\displaystyle-4\pi\mathsf{G}\mathcal{H}_l
\frac{\rho }{\mathsf{H}[\rho]Er}\frac{\mathsf{c}^2\lambda}{\mathscr{S}_l^2g^2} \\
\\
\displaystyle -4\pi\mathsf{G}\dot{\mathcal{H}}_l\frac{\rho}{\mathsf{H}[\rho]gE} &
\displaystyle 1-4\pi\mathsf{G}\dot{\mathcal{H}}_l
\frac{\rho }{\mathsf{H}[\rho]Er}\frac{\mathsf{c}^2\lambda}{\mathscr{S}_l^2g^2}
\end{bmatrix}^{-1}
\begin{bmatrix}
\mathcal{H}_lY \\
\\
\dot{\mathcal{H}}_lY
\end{bmatrix}. \label{G33}
\end{equation}

\begin{Definition}
We denote the solution  \eqref{G33} of the system of equations \eqref{G30a}\eqref{G30b} for \eqref{G29} by
\begin{equation}
\delta\check{\Phi}=
X=4\pi\mathsf{G}\mbox{\boldmath$\mathcal{H}$}_lY,\qquad 
r\frac{d}{dr}\delta\check{\Phi}=\dot{X}=4\pi\mathsf{G}\dot{\mbox{\boldmath$\mathcal{H}$}_l}Y.\label{G34}
\end{equation} 
\end{Definition}

\begin{Proposition}
The operators $\mbox{\boldmath$\mathcal{H}$}_l, \dot{\mbox{\boldmath$\mathcal{H}$}_l}$ are bounded linear operators from
$L^1([0,R]; rdr)$ to $L^{\infty}([0,R])$ and Lipschitz continuous in $\lambda\in [0,\lambda_0]$ with respect to the operator norm. 
\end{Proposition}

Now let us look at the equation \eqref{G24}. 
The variables $\displaystyle \frac{d}{dr}\delta\check{\Phi}, \delta\check{\Phi}$ in the coefficient $K$ given by \eqref{G24c} are replaced by those 
which are determined from $Y$ through \eqref{G34}.\\

Anyway we introduce the variable $\Psi$ by
\begin{equation}
\eta=W\Psi,
\end{equation}
where
\begin{align}
W&=\exp\Big[-\frac{1}{2}\int^rA\Big]=\mathfrak{Q}^{\frac{1}{2}}
\exp\Big[
2\pi\mathsf{G}\int_0^r\frac{\rho }{\mathsf{H}[\rho]E}\frac{\mathsf{c}^2\lambda}{\mathscr{S}_l^2g^2}\Big] \nonumber \\
&=\frac{\rho^{\frac{1}{2}}gE^{\frac{1}{2}}}{r^2}
\exp\Big[
2\pi\mathsf{G}\int_0^r\frac{\rho }{\mathsf{H}[\rho]E}\frac{\mathsf{c}^2\lambda}{\mathscr{S}_l^2g^2}\Big].
\end{align}
Then the equation \eqref{G24} is reduced to
\begin{equation}
-\frac{d^2\Psi}{dr^2}+
\Big(\frac{1}{4}A^2+\frac{1}{2}\frac{dA}{dr}-B\Big)\Psi-
\frac{\rho K}{W}=0.
\end{equation}
We write this equation as
\begin{equation}
-\frac{d^2\Psi}{dr^2}+b\Psi+f=\frac{1}{\lambda}\kappa\Psi,
\end{equation}
where
\begin{align}
b&=\frac{1}{4}A^2+\frac{1}{2}\frac{dA}{dr}-B-\frac{l(l+1)\mathfrak{N}^2}{\lambda r^2}= \nonumber \\
&=\frac{1}{4}A^2+\frac{1}{2}\frac{dA}{dr}
-\frac{g}{\mathsf{c}^2}\Big(\frac{1}{\mathsf{H}[g]}+\frac{2}{r}\Big)
-\frac{\lambda}{\mathsf{c}^2}
-\frac{l(l+1)}{r^2} + \nonumber \\
&+4\pi\mathsf{G}\Big(\frac{1}{\mathsf{H}[\rho]gE}-\frac{1}{\mathsf{c}^2}\Big), \label{G38}
\end{align}
while
$$
A=\frac{1}{\mathsf{H}[\rho g^2/r^4]}+
\frac{1}{\mathsf{H}[E]}
-4\pi\mathsf{G}\frac{\rho }{\mathsf{H}[\rho]E}\frac{\mathsf{c}^2\lambda}{\mathscr{S}_l^2g^2}, $$
and
\begin{align}
f&=-\frac{\rho K}{W}, \\
\kappa&=\frac{l(l+1)\mathfrak{N}^2}
{r^2}.
\end{align}

The coefficients $b, f, \kappa$ are Lipschitz continuous function of $\lambda \in [0,\lambda_0]$. If we consider this $\lambda$ as a parameter and consider the eigenvalue problem for the eigenvalue $\Lambda$:
\begin{equation}
-\frac{d^2\Psi}{dr^2}+b\Psi+f=\Lambda \kappa \Psi, \label{GEP}
\end{equation}
and if the eigenvalue $\Lambda$ coincides with $1/\lambda$, then this gives the eigenvalue of the original problem under consideration. This is the key idea.\\

Let us perform the Liouville transformation on the eigenvalue problem \eqref{GEP}.
Putting
\begin{align}
x&= \int_0^r\sqrt{\kappa}dr, \qquad y=\kappa^{\frac{1}{4}}\Psi, \nonumber \\
\hat{f}&=\kappa^{-\frac{3}{4}}f, \nonumber \\
q&=\frac{b}{\kappa}+\frac{1}{4}\frac{1}{\kappa}\Big[
\frac{d^2}{dr^2}\log \kappa -\frac{1}{4}\Big(\frac{d}{dr}\log \kappa\Big)^2\Big],
\end{align}
we transform \eqref{GEP} to
\begin{equation}
-\frac{d^2y}{dx^2}+qy+\hat{f}=\Lambda y. \label{GEP2}
\end{equation}
See \cite[p. 275, Theorem 6]{BirkhoffR} or \cite[p.110]{Yosida}.\\

We see
\begin{equation}
x_+:=\int_0^R\sqrt{\kappa}dr=\int_0^R
\sqrt{\frac{l(l+1)\mathfrak{N}^2}{r^2}}dr < \infty
\end{equation}
and
\begin{subequations}
\begin{align}
&x=\sqrt{\frac{l(l+1)\mathfrak{N}_{\mathsf{O}1}^2}{2}}r(1+[r^2]_1)\quad\mbox{as}\quad r\rightarrow +0 \\
&x_+-x=\frac{\sqrt{l(l+1)\mathfrak{N}_R^2}}{R}(R-r)(1+[R-r, (R-r)^{\nu+1}]_1)\quad
\mbox{as}\quad r\rightarrow R-0.
\end{align}
\end{subequations}
The $r$-interval $[0,R]$ is mapped onto the $x$-interval $[0,x_+]$. 

By a tedious calculation, we see
\begin{equation}
q=
\begin{cases}
\displaystyle \frac{l(l+1)}{x^2}(1+[x^2]_1) \quad\mbox{as}\quad x\rightarrow +0 \\
\\
\displaystyle \frac{\nu^2}{(x_+-x)^2}(1+[x_+-x, (x_+-x)^{\nu+1}]_1)\quad\mbox{as}\quad x \rightarrow x_+-0
\end{cases}
\end{equation}
Let us note that $\displaystyle l(l+1) \geq 2 >\frac{3}{4}, \nu^2 >1 > \frac{3}{4}$. \\

We use the following

\begin{Notation}
For any quantities $Q_1, Q_2$, we denote
\begin{equation}
Q_1\vee Q_2:=\max\{ Q_1, Q_2\},\qquad
Q_1\wedge Q_2:=\min\{ Q_1, Q_2\}.
\end{equation}
\end{Notation}

We claim
\begin{Proposition}
There exist constants $\displaystyle K_0, K_1, \frac{3}{4} < K_1 <  l(l+1) \wedge \nu^2$, independent of 
$\lambda \in [0,\lambda_0]$ and $ \mathsf{G}\in [0, \mathsf{G}_0]$ such that
\begin{equation}
\tilde{q}(x;\lambda):=q+K_0 \geq 1+K_1\Big(\frac{1}{x^2}\vee \frac{1}{(x_+-x)^2}\Big)
\end{equation}
for $0<x<x_+$.
\end{Proposition}

In fact, the leading term of $q$ comes from
$$
\frac{1}{\mathsf{H}[\rho g^2/r^4]}=
\begin{cases}
\displaystyle -\frac{2}{r}(1+[r^2]_1)\quad\mbox{as}\quad r \rightarrow +0 \\
\\
\displaystyle \frac{\nu}{R-r}(1+[R-r, (R-r)^{\nu+1}]_1)\quad\mbox{as}\quad r \rightarrow R-0,
\end{cases}.
$$
Actually we note
\begin{align*}
\frac{1}{\mathsf{H}[E]}&=-\frac{1}{E}\Big(\frac{dE_1}{dr}+\lambda\frac{dE_2}{dr}\Big)
\frac{\lambda}{l(l+1)}, \\
\frac{d}{dr}\frac{1}{\mathsf{H}[E]}&=\Big[\Big(\frac{1}{E}\Big(\frac{dE_1}{dr}+\lambda\frac{dE_2}{dr}\Big)\Big)^2-
\frac{1}{E}\Big(\frac{d^2E_1}{dr^2}+\lambda\frac{d^2E_2}{dr^2}\Big)\Big]\frac{\lambda}{l(+1)},
\end{align*}
while
\begin{subequations}
\begin{align}
&r\frac{dE_1}{dr}=\frac{r}{g}\Big(3-4\pi\mathsf{G}\rho\frac{r}{g}\Big)
\Big(-4+8\pi\mathsf{G}\rho\frac{r}{g}\Big)
+4\pi\mathsf{G}r\frac{d\rho}{dr}\Big(\frac{r}{g}\Big)^2, \\
&r\frac{dE_2}{dr}=-2\Big(\frac{r}{g}\Big)^2\Big(3-4\pi\mathsf{G}\rho\frac{r}{g}\Big).
\end{align}
\end{subequations}\\

Moreover we have

\begin{Proposition}\label{Prop.G4}
There exists a constant $M_q$ independent of $\lambda \in [0,\lambda_0], \mathsf{G}\in [0,\mathsf{G}_0]$ such that 
\begin{equation}
\tilde{q}(x;\lambda)=\tilde{q}(x;0)(1+\Omega_q(x,\lambda)\lambda)
\end{equation}
with
\begin{equation}
|\Omega_q(x,\lambda)|\leq M_q, \qquad
\Big|\frac{\partial}{\partial\lambda}\Omega_q(x,\lambda)\Big|\leq M_q
\end{equation}
for $0<x<x_+, 0\leq\lambda \leq \lambda_0,
0\leq \mathsf{G}\leq\mathsf{G}_0$. Here $\Omega_q(x,\lambda)$ is a meromorphic function of $\lambda$. In particular it holds that
\begin{equation}
(1-\varepsilon)\tilde{q}(x;0)\leq
\tilde{q}(x;\lambda)\leq (1+\varepsilon)
\tilde{q}(x;0)
\end{equation}
with $\varepsilon=M_q\lambda_0$, which can be supposed to satisfy
$$0<1-\varepsilon,\qquad
\frac{3}{4}<(1-\varepsilon)K_1.
$$
\end{Proposition}

Actually 
putting
\begin{align*}
[*]:=& \Big[ \frac{1}{4}A^2+\frac{1}{2}\frac{dA}{dr}\Big]
- \Big[\frac{1}{4}A^2+\frac{1}{2}\frac{dA}{dr}\Big]_{\lambda=0} \\
&=\frac{1}{2}\frac{1}{\mathsf{H}[\rho g^2/r^4]}\cdot\Big(\frac{1}{\mathsf{H}[E]}
-4\pi\mathsf{G}\frac{\rho }{\mathsf{H}[\rho]E}\frac{\mathsf{c}^2\lambda}{\mathscr{S}_l^2g^2}\Big)+
\frac{1}{4}\Big(\frac{1}{\mathsf{H}[E]}
-4\pi\mathsf{G}\frac{\rho }{\mathsf{H}[\rho]E}\frac{\mathsf{c}^2\lambda}{\mathscr{S}_l^2g^2}\Big)^2 \\
&+\frac{1}{2}\frac{d}{dr}\Big(\frac{1}{\mathsf{H}[E]}
-4\pi\mathsf{G}\frac{\rho }{\mathsf{H}[\rho]E}\frac{\mathsf{c}^2\lambda}{\mathscr{S}_l^2g^2}\Big),
\end{align*}
we can estimate it, as $ r\rightarrow R-0$,
by
\begin{align*}
&\frac{dE_1}{dr}, \frac{dE_2}{dr}, \frac{d^2E_2}{dr^2}=O(1),\quad \frac{d^2E_1}{dr^2}=
O(1)+O((R-r)^{\nu-2}) , \\
& \frac{\rho }{\mathsf{H}[\rho]E}\frac{\mathsf{c}^2\lambda}{\mathscr{S}_l^2g^2} =O(1)\frac{\lambda}{l(l+1)}, 
\quad \frac{d}{dr}\frac{\rho }{\mathsf{H}[\rho]E}\frac{\mathsf{c}^2\lambda}{\mathscr{S}_l^2g^2} =[O(1)+O((R-r)^{\nu-2})]\frac{\lambda}{l(l+1)}
\end{align*}
to get 
$$[*] =[O(1)+O((R-r)^{\nu-2}]\frac{\lambda}{l(l+1)} =\frac{1}{R-r} O(\lambda)$$
as $r \rightarrow R-0$, and so on.  \\

We put
\begin{equation}
Q_0[y;\lambda]:=\int_0^{x_+}
\Big(\Big|\frac{dy}{dx}\Big|^2+\tilde{q}(x;\lambda)|y|^2\Big)dx
\end{equation}
and
\begin{align}
&\mathfrak{E}:=L^2([0,x_+]; dx) \\
&\mathfrak{E}_1:=\{ y \in \mathfrak{E}\  |\  Q_0[y; 0] <\infty \} .
\end{align}
Note that we know
$$(1-\varepsilon)Q_0[y;0]
\leq Q_0[y;\lambda]
\leq(1+\varepsilon)Q_0[y; 0]. $$\\

Now we have to consider the quadratic form $Q[y; \lambda]$ defined by
\begin{align}
Q[y;\lambda]&:=\int_0^{x_+}\Big(
\Big|\frac{dy}{dx}\Big|^2+\tilde{q}(x;\lambda)|y|^2+\hat{f}y^*\Big)dx \nonumber  \\
&=Q_0[y;\lambda]+\int_0^{x_+}\hat{f}y^*dx.
\end{align}

In order to estimate the perturbation
$Q[y; \lambda]-Q_0[y; \lambda]$ due to the perturbation of self-gravitation,
we  should  analyze 
\begin{align}
\hat{f}&=-\Big(\frac{l(l+1)\mathfrak{N}^2}{r^2}\Big)^{-\frac{3}{4}}
\frac{\rho}{W}\Big[
\Big(\frac{1}{\mathsf{H}[g^2E/r^2]}-4\pi\mathsf{G}
\frac{\rho }{\mathsf{H}[\rho]E}\frac{\mathsf{c}^2\lambda}{\mathscr{S}_l^2g^2}\Big)\frac{\dot{X}}{r}+ \nonumber \\
&+\Big(
\frac{1}{\mathsf{H}[E]}\frac{\mathscr{S}_l^2g}{\mathsf{c}^2\lambda}
-4\pi\mathsf{G}\frac{\rho}{\mathsf{H}[\rho]gE}\Big)X\Big].
\end{align}
Recall that $\dot{X}, X$ are given by
\begin{equation}
X=4\pi\mathsf{G}\mbox{\boldmath$\mathcal{H}$}_lY,\quad \dot{X}=4\pi\mathsf{G}\dot{\mbox{\boldmath$\mathcal{H}$}_l}Y 
\end{equation}
for
\begin{align}
Y&=\frac{1}{\mathsf{H}[\rho]Er}\frac{\mathsf{c}^2\lambda}{\mathscr{S}_l^2g^2}\cdot r\frac{d\eta}{dr}+
\Big(\frac{1}{\mathsf{H}[\rho]gE}-\frac{1}{\mathsf{c}^2}\Big)\eta \nonumber \\
&=
B_1\frac{dy}{dx}+(B_{01}+B_{02}+B_{03}+B_{04})y,
\end{align}
where
\begin{subequations}
\begin{align}
B_{1}&=\frac{W\kappa^{\frac{1}{4}}}{\mathsf{H}[\rho]E}\cdot 
\frac{\mathsf{c}^2\lambda}{\mathscr{S}_l^2g^2}, \label{05.66a}\\
B_{01}&=\frac{1}{4}\frac{W\kappa^{\frac{1}{4}}}{\mathsf{H}[\rho]E}
\frac{1}{\mathsf{H}[\kappa]}\cdot 
\frac{\mathsf{c}^2\lambda}{\mathscr{S}_l^2g^2}, \\
B_{02}&=-\frac{W\kappa^{\frac{1}{4}}}{\mathsf{H}[\rho]}\frac{1}{\mathsf{H}[W]}\cdot
\frac{\mathsf{c}^2\lambda}{\mathscr{S}_l^2g^2}, \\
B_{03}&=\frac{W\kappa^{-\frac{1}{4}}}{\mathsf{H}[\rho]gE}, \\
B_{04}&=-\frac{W\kappa^{-\frac{1}{4}}}{\mathsf{c}^2}. \label{05.66e}
\end{align}
\end{subequations}\\

In this sense we consider that $\hat{f}=\hat{\mathfrak{f}}y$, $\hat{\mathfrak{f}}$ being an operator acting on functions $y$, and we are going to consider
\begin{equation}
Q[y;\lambda]=Q_0[y;\lambda]+(\hat{\mathfrak{f}}y|y)_{\mathfrak{E}}.
\end{equation}
Here the symbol $(\hat{\mathfrak{f}}y|y)_{\mathfrak{E}}$ is a kind of diversion, since we do not claim that $\hat{f} \in \mathfrak{E}$ for any $y \in \mathfrak{E}_1$, but it has a firm meaning as explained later.\\

Let us introduce

\begin{Definition}
We put
\begin{equation}
\langle \rho \rangle (r):=
\begin{cases}
\rho_{\mathsf{O}} \quad&\mbox{for}\quad r=0 \\
\\
\displaystyle \frac{4\pi}{r^3}\int_0^r\rho(\acute{r})\acute{r}^2d\acute{r} \quad&\mbox{for}\quad 0<r \leq R.
\end{cases}
\end{equation}
\end{Definition}

Then $\langle \rho \rangle (r)$ and $1/\langle \rho \rangle (r)$ are positive continuous functions of $r \in [0,R]$, and $\langle \rho \rangle (r)\leq \rho_{\mathsf{O}}$
for $0\leq r\leq R$. Recall that
\begin{equation}
\frac{g}{r}=\frac{1}{r}\frac{d\Phi}{dr}=\frac{4\pi\mathsf{G}}{r^3}\int_0^r\rho(\acute{r})\acute{r}^2d\acute{r} =\mathsf{G}\langle \rho \rangle .
\end{equation}

Suppose\\

{\bf (B0): } \  {\it There is a  small positive number $\delta_{\mathbf{B}}$
such that
\begin{equation}
\lambda_0(1+\rho\langle\rho\rangle^{-1})\langle \rho \rangle^{-1} \leq\delta_{\mathbf{B}}.
\end{equation}
}\\

Recall that $\lambda$ is confined in the interval $[0, \lambda_0]$. Here $\lambda_0$ can be fixed  arbitrarily small so that $\delta_{\mathbf{B}}$ be necessarily small, while $\rho=\bar{\rho}$ is fixed. 
Then,
looking at \eqref{DefEa}, \eqref{DefEb}, \eqref{DefEc}, we see that, for $0\leq \lambda \leq \lambda_0$,
$E-1 \ll 1$ holds, provided $\delta_{\mathbf{B}} \ll 1$. Moreover we put the following condition :\\

{\bf (B1):}\  {\it It holds 
\begin{equation}
-\frac{1}{r}\frac{d\rho}{dr}\cdot 
\langle \rho \rangle^{-1} \leq \delta_{\mathbf{B}}.
\end{equation}
Here $\delta_{\mathbf{B}}$ is a sufficiently small positive number.}\\

Then, by taking $\delta_{\mathbf{B}}$ sufficiently small, we can claim the following\\

\begin{Proposition}
If both {\bf (B0)} and {\bf (B1)} hold
with a sufficiently small $\delta_{\mathbf{B}}$, then the condition {\bf (G)} holds.
\end{Proposition}

In fact we see
\begin{align*}
&\frac{\rho}{\mathsf{H}[\rho]gE} \quad\lesssim\quad -\frac{1}{r}\frac{d\rho}{dr}\frac{r}{g}, \\
&\frac{\rho}{\mathsf{H}[\rho]Er}\frac{\mathsf{c}^2\lambda}{\mathscr{S}_l^2g^2}\quad\lesssim\quad -\frac{1}{r}\frac{d\rho}{dr}\Big(\frac{r}{g}\Big)^2\lambda
\leq -\frac{1}{r}\frac{d\rho}{dr}\frac{r}{g}\cdot \lambda_0\frac{r}{g}.
\end{align*}\\

Let us consider
\begin{equation}
\hat{f}=(C_{11}+C_{12})\dot{X}+(C_{01}+C_{02})X,
\end{equation}
where
\begin{subequations}
\begin{align}
C_{11}&=-\kappa^{-\frac{3}{4}}\frac{\rho}{W}
\cdot \frac{1}{\mathsf{H}[g^2E/r^2]}\frac{1}{r}, \label{05.73a}\\
C_{12}&=
\kappa^{-\frac{3}{4}}\frac{\rho}{W}
\cdot 4\pi\frac{\rho}{\mathsf{H}[\rho]E}\frac{\mathsf{c}^2\lambda}{\mathscr{S}_l^2g^2}
\frac{1}{r}, \\
C_{01}&=-
\kappa^{-\frac{3}{4}}\frac{\rho}{W}
\cdot \frac{1}{\mathsf{H}[E]}
\frac{\mathscr{S}_l^2g}{\mathsf{c}^2\lambda}, \\
C_{02}&=
\kappa^{-\frac{3}{4}}\frac{\rho}{W}
\cdot 4\pi\frac{\rho}{\mathsf{H}[\rho]gE}. \label{05.73d}
\end{align}
\end{subequations}

Put
\begin{equation}
C_{\mathbf{B1}}:=(l+1)(\|xC_{11}\|_{L^{\infty}}+
\|xC_{12}\|_{L^{\infty}})+
\|xC_{01}\|_{L^{\infty}}+
\|xC_{02}\|_{L^{\infty}}.
\end{equation}

Then we have
\begin{equation}
\|x\hat{f}\|_{L^{\infty}}\leq C(2l+1)^{-1}C_{\mathbf{B1}}\|Y\|_{L^1(rdr)},
\end{equation}
since $
|\|\mbox{\boldmath$\mathcal{H}$}_l\||_{L1(rdr)\rightarrow L^{\infty}} \lesssim (2l+1)^{-1}, 
|\|\dot{\mbox{\boldmath$\mathcal{H}$}_l}\||_{L^1(rdr)\rightarrow L^{\infty}} \lesssim (l+1)(2l+1)^{-1}$.

Since 
$$rdr=\Big(\frac{1}{\sqrt{\kappa(r)}}\frac{1}{x}\int_0^x\frac{d\acute{x}}{\sqrt{\kappa(\acute{r})}}\Big)\cdot xdx,$$
we see 
$$ \|Y\|_{L^1(rdr)} \leq \Bigg\|\frac{1}{\sqrt{\kappa}}\Bigg\|_{L^{\infty}}^2\|Y\|_{L^1(xdx)}.$$
Therefore we have
\begin{equation}
\|x\hat{f}\|_{L^{\infty}}
\leq C(2l+1)^{-1}\Bigg\|\frac{1}{\kappa}\Bigg\|_{L^{\infty}}C_{\mathbf{B1}} \label{05.76}
\|Y\|_{L^1(xdx)}.
\end{equation}

Looking at
$$
Y=B_1\frac{dy}{dx}+(B_{01}+B_{02}+B_{03}+B_{04})y
$$
with $B_1, B_{01}, \cdots, B_{04}$ given by \eqref{05.66a} $\sim$ \eqref{05.66e},
we put
\begin{align}
C_{\mathbf{B2}}&:=
\|x_+^{\frac{\nu+1}{2}}(x_+-x)^{-\frac{\nu}{2}+1}B_1\|_{L^{\infty}} + \nonumber \\
&+\sum_{k=1}^4
\|x_+^{\frac{1}{2}}x^2\wedge x_+^{\frac{\nu+3}{2}}(x_+-x)^{-\frac{\nu}{2}+1}B_{0k}\|_{L^{\infty}}.
\end{align}

Then we have
$$\Big|B_1\frac{dy}{dx}\Big|\leq C_{\mathbf{B2}}
x_+^{-\frac{\nu+1}{2}}(x_+-x)^{\frac{\nu}{2}-1} \Big|\frac{dy}{dx}\Big|$$
so that
\begin{align*}
\Big\|B_1\frac{dy}{dx}\Big\|_{L^1(xdx)}&\leq C_{\mathbf{B2}}
x_+^{-\frac{\nu+1}{2}}
\int_0^{x_+}
(x_+-x)^{\frac{\nu}{2}-1}\Big|\frac{dy}{dx}\Big|xdx \\
&\leq 
CC_{\mathbf{B2}}\Big(\int_0^{x_+}\Big|\frac{dy}{dx}\Big|^2dx\Big)^{\frac{1}{2}}
\leq
C C_{\mathbf{B2}}\sqrt{Q_0[y; 0]}.
\end{align*}
And we have
$$|B_0y|\leq C_{\mathbf{B2}}
(x_+^{-\frac{1}{2}}x^{-2}\vee x_+^{-\frac{\nu+3}{2}}(x_+-x)^{\frac{\nu}{2}-1})|y|,$$
$B_0$ being
$B_{01}+B_{02}+B_{03}+B_{04}$, so that
\begin{align*}
\|B_0y\|_{L^1(xdx)}&\leq C_{\mathbf{B2}}
x_+^{-\frac{1}{2}}
\int_0^{x_m}x^{-2}|y|xdx + \\
&+C_{\mathbf{B2}}
 x_+^{-\frac{\nu+3}{2}}
\int_{x_m}^{x_+} (x_+-x)^{\frac{\nu}{2}-1}  |y| xdx \\
&\leq C C_{\mathbf{B2}}\Big(\int_0^{x_m}\frac{1}{x^2}|y|^2dx\Big)^{\frac{1}{2}} \\
&+ C C_{\mathbf{B2}}\Big(\int_{x_m}^{x_+}
\frac{1}{(x_+-x)^2}|y|^2dx\Big)^{\frac{1}{2}} \\
&\leq C' C_{\mathbf{B2}}\sqrt{Q_0[y;0]},
\end{align*}
where $x_m \in ]0,x_+[$ denotes the solution of the equation
$$x_+^{-\frac{1}{2}}x^{-2}=
x_+^{-\frac{\nu+3}{2}}(x_+-x)^{\frac{\nu}{2}-1}.$$
Here we recall 
$$Q_0[y;0] \geq \frac{3}{4}\int_0^{x_+}\frac{1}{x^2}|y|dx, \quad\mbox{and}\quad
\geq \frac{3}{4}\int_0^{x_+}\frac{1}{(x_+-x)^2}|y|dx. $$\\

Thus we can claim the following

\begin{Proposition}
Suppose {\bf (B0), (B1)}. 
If $y \in \mathfrak{E}_1$, then $Y \in L^1([0,x_+]; xdx)$ and there
 exists a constant $M_Y$ independent of $\lambda \in [0,\lambda_0]$ such that
\begin{align}
&\|Y\|_{L^1(xdx)} \leq M_Y\sqrt{Q_0[y;0]} \\
&\|Y|_{\lambda+\delta\lambda}-Y|_{\lambda}\|_{L^1(xdx)} \leq
M_Y|\delta\lambda|\sqrt{Q_0(y;0]}.
\end{align}
Here we can put
\begin{equation}
M_Y=C(2l+1)^{-1}\Bigg\|\frac{1}{\kappa}\Bigg\|_{L^{\infty}}C_{\mathbf{B1}}C_{\mathbf{B2}}
\end{equation}
with a constant $C$ which depends only on $\nu$.
\end{Proposition}

Consequently we have
\begin{Proposition}\label{Prop.G6}
Suppose {\bf (B0), (B1)}. 
If $y \in \mathfrak{E}_1$, then $x\hat{f} \in L^{\infty}([0, x_+])$ and
it holds that 
\begin{align}
&|(\hat{f}|y)_{\mathfrak{E}}| \leq \varepsilon_{\mathbf{B}} Q[y;0] \label{G56}\\
&|(\hat{f}|_{\lambda+\delta\lambda}-\hat{f}|_{\lambda}|y)_{\mathfrak{E}}|\leq
\varepsilon_{\mathbf{B}}|\delta\lambda|Q[y;0].
\end{align}
Here $\varepsilon_{\mathbf{B}}$ is the small positive number given by
\begin{equation}
\varepsilon_{\mathbf{B}}=C(2l+1)^{-1}\Bigg\|\frac{1}{\kappa}\Bigg\|_{L^{\infty}}
C_{\mathbf{B1}}C_{\mathbf{B2}}\|\kappa\|_{L^{\infty}}^{\frac{1}{4}}\sqrt{R}. \label{05.83}
\end{equation}
\end{Proposition}

Here we note $x_+ \leq \|\sqrt{\kappa}\|_{L^{\infty}} R$. \\

Note that we do not claim that $\hat{f}=\hat{\mathfrak{f}}y \in \mathfrak{E}$ for any $y \in \mathfrak{E}_1$ so that
the symbol $(\hat{f}|y)_{\mathfrak{E}}$ is a diversion, which is done by observing that, for any $y \in \mathfrak{E}_1$, 
\begin{align*}
\Big|\int_0^{x_+}\hat{f}y^*dx\Big| &=
\Big|\int_0^{x_+}x\hat{f}\cdot \Big(\frac{y}{x}\Big)^*dx\Big|
\leq \|x\hat{f}\|_{L^{\infty}}\int_0^{x_+}\frac{|y|}{x}dx \\
&\leq \|x\hat{f}\|_{L^{\infty}}\sqrt{x_+}\Big(\int_0^{x_+}\frac{|y|^2}{x^2}dx\Big)^{1/2}
< +\infty.
\end{align*}
which \eqref{G56} means.\\

Here we suppose \\

{\bf (B2):}\  {\it The positive number $\varepsilon_{\mathbf{B}}$ is sufficiently small.}\\

Thanks to Proposition \ref{Prop.G6}, we can claim 

\begin{Proposition}\label{Prop.G7}
Suppose {\bf (B0),(B1),(B2)} with sufficiently small $\delta_{\mathbf{B}},\varepsilon_{\mathbf{B}}, \lambda_0$.
Then it holds that

1) $$Q[y;\lambda] <\infty$$ for $y \in \mathfrak{E}_1$;

2) $$(1-\varepsilon)Q_0[y;0]\leq Q[y;\lambda]
\leq (1+\varepsilon)Q_0[y;0] $$
for $0\leq \lambda\leq \lambda _0$ with $0<\varepsilon \ll 1$;

3)$$|Q[y;\lambda+\delta\lambda]-Q[y;\lambda]|\leq M_Q|\delta\lambda|Q_0[y;0]$$
with a constant $M_Q$ independent of $\lambda \in [0,\lambda_0], 
y \in \mathfrak{E}_1$. 
\end{Proposition}

Since the imbedding $\mathfrak{E}_ 1 \hookrightarrow \mathfrak{E}$ is compact, we can claim that the eigenvalue problem \eqref{GEP} is of the Sturm-Liouville type,
that is, the operator
$\displaystyle  -\frac{d^2}{dx^2}+q + \hat{\mathfrak{f}} $ restricted on $C_0^{\infty}(]0,x_+[)$ admits the Friedrichs extension, a self-adjoint operator in $\mathfrak{E}=L^2([0,x_+])$, and its spectrum consists of simple eigenvalues $(\Lambda_n(\lambda))_{n\in\mathbb{N}}$ such that
$$\Lambda_1(\lambda)<\Lambda_2(\lambda)<\cdots<\Lambda_n(\lambda)<\cdots \rightarrow +\infty.$$
See \cite[Kapiter VII]{CourantH} and  
\cite[p.159, Theorem X.10]{ReedS}.

Of course we are considering $\lambda \in [0,\lambda_0]$, 
$\delta_{\mathbf{B}}$
being sufficiently small. \\

Now the eigenvalue $\Lambda_n(\lambda)$ is given by the Max-Min principle as follows:

For any $V =\{v_1, \cdots, v_n\} \subset \mathfrak{E}$ we put
\begin{equation}
\mathfrak{d}(V;\lambda)=\inf\{\  Q([y;\lambda]\  |\  y \in \mathfrak{E}_1, \|y\|_{\mathfrak{E}}=1,
(y|v_j)_{\mathfrak{E}}=0 \forall v_j \in V\}.
\end{equation}
Then
 it holds that
\begin{equation}
\Lambda_n(\lambda)+K_0=\sup\{\  \mathfrak{d}(V;\lambda)\  |\  V \subset \mathfrak{E}, \sharp V=n\}.
\end{equation}

As for the theory of the Max-Min principle, see e.g., \cite[Chapter 11]{Helffer}. The above characterization of $\Lambda_n$ is given as 
\cite[p. 144, (11.3.1)]{Helffer}.\\

Thanks to Proposition \ref{Prop.G7} we can claim
\begin{Proposition}\label{Prop.G8}
For each $n$ it holds that
\begin{equation}
\frac{1-\varepsilon}{1+\varepsilon}(\Lambda_n(0)+K_0)
\leq\lambda_n(\lambda)+K_0
\leq \frac{1+\varepsilon}{1-\varepsilon}(\Lambda_n(0)+K_0),
\end{equation}
and the function $\lambda \mapsto \Lambda_n(\lambda)$ is continuous on $[0,\lambda_0]$.
\end{Proposition}

Since $\Lambda_n(\lambda_0) \rightarrow +\infty$ as $n \rightarrow \infty$, we can find $n_0\in \mathbb{N}$ such that $\displaystyle \Lambda_n(\lambda_0) >\frac{1}{\lambda_0}$
for $n \geq n_0$. Note that we can suppose
$$\Lambda_n(\lambda)\geq
\frac{1-\varepsilon}{1+\varepsilon}\Lambda_n(0)-\frac{2\varepsilon}{1+\varepsilon}K_0 >0,
$$
for $0\leq \lambda \leq \lambda_0, n \geq n_0$, by replacing $n_0$ by a greater one if necessary. Then by Proposition \ref{Prop.G8} the function $$
\varphi\  : \  \lambda \mapsto \lambda-\frac{1}{\Lambda_n(\lambda)} $$
is continuous on $[0,\lambda_0]$. Since
$$\varphi(0) =-\frac{1}{\Lambda_n(0)} <0,\qquad \varphi(\lambda_0)=\lambda_0
-\frac{1}{\Lambda_n(\lambda_0)} >0,
$$
there exists at least one $\lambda\in ]0,\lambda_0[$ such that $\varphi(\lambda)=0$, that is, 
$$\Lambda_n(\lambda)=\frac{1}{\lambda}.$$
Although we cannot claim that the solution is unique, we can denote a solution by $\lambda_{-n}$ choosing one of them. Then Proposition \ref{Prop.G8} implies
$$\lambda_{-n}=\frac{1}{\Lambda_n(\lambda_{-n})}\leq
\frac{1+\varepsilon}{(1-\varepsilon)\Lambda_n(0)-2\varepsilon K_0}
\rightarrow 0$$
as $n \rightarrow \infty$.

Summing up, we get

\begin{Theorem}\label{Th.G}
Let $l \geq 1$. Suppose Assumptions \ref{Ass.G1}, \ref{Ass.G2}, \ref{Ass.G3}, and suppose 
{\bf (B0), (B1), (B2)} with sufficiently small $\delta_{\mathbf{B}}, \varepsilon_{\mathbf{B}}$.
Then
there exists a sequence $(\lambda_{-n})_{n\in\mathbb{N}}$ of positive eigenvalues of $\vec{L}_l$ such that $\lambda_{-n} \rightarrow 0$ as $n\rightarrow \infty$.
\end{Theorem}

Moreover we claim:

\begin{Theorem}\label{Th.ExG}
There are admissible equilibria for which the conditions of Theorem \ref{Th.G} are satisfied at least for large $l$.
\end{Theorem}

Proof. Let us fix an admissible equilibrium $(\rho_1, S_1)$ which enjoys Assumptions \ref{Ass.G1}, \ref{Ass.G2}, \ref{Ass.G3}. Actually it exists due to the discussion of Section 2, Theorem \ref{Th.1}. Using positive parameters $\tau$, we look at the equilibrium $(\rho, S)$ given by
$$\rho(r)=\tau\rho_1(\tau r),\qquad
S(r)=S_1(\tau r)-\gamma \mathsf{C}_V\log\tau.$$
Of course the corresponding pressure and gravitational potential are $P=P_1, \Phi=\tau^{-1}\Phi_1$ with denoting $r_1=\tau r$. Then we see
$$(1+\rho\langle\rho\rangle^{-1})\langle\rho\rangle^{-1}=\tau^{-1}
(1+\rho_1\langle\rho_1\rangle^{-1})\langle\rho_1\rangle^{-1}
$$ and
$$
-\frac{1}{r}\frac{d\rho}{dr}\langle\rho\rangle^{-1}=
-\tau^2
\frac{1}{r_1}\frac{d\rho_1}{dr_1}\langle\rho_1\rangle^{-1}.
$$
Therefore, taking $\tau$ sufficiently small, and taking $\lambda_0$ so small that $\tau^{-1}\lambda_0$ be small compared with the small $\tau$, the conditions {\bf (B0), (B1)} are satisfied with an arbitrarily small $\delta_{\mathbf{B}}$. Let us fix such an equilibrium. 

Then we note
$$ \frac{1}{M_E}l(l+1) \leq \kappa=\frac{l(l+1)\mathfrak{N}^2}{r^2}\leq M_El(l+1).$$
Here and hereafter $M_E$ stands for various constants depending upon the fixed equilibrium with the fixed $\tau$. Looking at \eqref{05.73a} $\sim$ \eqref{05.73d}, we see
$$C_{\mathbf{B1}}\leq M_E(l(l+1))^{-\frac{1}{4}}.$$
Here we note that
$$x=\int_0^r\sqrt{\kappa}(\acute{r})d\acute{r} \leq M_E(l(l+1))^{\frac{1}{2}}$$ and that
$C_{01}$ reads
$$C_{01}=\kappa^{-\frac{3}{4}}\frac{\rho}{W}
\frac{1}{E}\Big(\frac{dE_1}{dr}+\lambda\frac{dE_2}{dr}\Big)\frac{g}{r^2}.$$
Therefore \eqref{05.76} reads
$$\|x\hat{f}\|_{L^{\infty}}\leq
M_E(2l+1)^{-1}(l(l+1))^{-\frac{5}{4}}\|Y\|_{L^1(xdx)}.
$$
Looking at $B_1, B_{01}, \sim B_{04}$, we see
$$C_{\mathbf{B2}}\leq M_El(l+1).$$
Therefore \eqref{05.83} reads
$$
\varepsilon_{\mathbf{B}}\leq M_E
(2l+1)^{-1}.$$
Clearly $(2l+1)^{-1}$ is small as $l $ is large. This completes the proof. $\square$\\

{\bf Historical Remarks} \\

The limiting case when the terms $\delta\check{\Phi}, d\delta\check{\Phi}/dr$ are neglected in the equations is called `Cowling approximation'. In this case we can forget the perturbation of the self-gravitation, and we can neglect the term $\hat{f}=\hat{\mathfrak{f}}y $
in the eigenvalue problem \eqref{GEP2}. Hence the discussion can be reduced to a quite simpler one. Astrophysicists believe that the eigenvalue problem for this approximation  give, at least approximately, eigenvalues of the problem for $\vec{L}_l$. See, e.g., \cite{LedouxW}, \cite{Cox}, and so on. The priority of this observation may go back to the work by T. G. Cowling, \cite{Cowling}, on November 3, 1941. Actually \cite{Gough} derives the equation \eqref{G18}, but after that, analysis is done by neglecting the term $F$, saying
\begin{quote}
Cowling (1941) showed that, except for modes of low degree $l$ with a numerically small order $n$, the perturbation $\Phi'$ to the gravitational potential has a relatively minir effect on the modes. Although $\Phi'$ must be included in accurate numerical computations of all but the high-degree modes, it has little influence on the basic dynamics, and consequently I will ignore it. Thus I set $F=0$, and eq. (5.4.1) reduces to a single second-order differential equation for $\delta p$.
\end{quote}
See \cite{Gough}, p.439-440. Here $\Phi', \delta p$ of \cite{Gough} correspond to $\delta\check{\Phi},
\Delta \check{P}$ of this article. 

Moreover astrophysicists usually use the so called `planer approximation', which simplifies the analysis of the eigenvalue problem very much. In fact, in the preceding review \cite{DeubnerG} by D. O. Gough himself and in the later book \cite{Aerts} by C. Aerts, J. Christensen-Dalsgaard and D. W. Kurtz, the approximation
$$
{\mathsf{H}[\mathfrak{Q}]}  \approx {\mathsf{H}[\rho]}, \quad
\mathfrak{N}^2  \approx \mathscr{N}^2, $$
and
$$
\frac{\lambda_c}{\mathsf{c}^2}
=\frac{1}{4}A^2+\frac{1}{2}\frac{dA}{dr}-\frac{g}{\mathsf{c}^2}
\Big(\frac{1}{\mathsf{H}[g]}+\frac{2}{r}\Big) 
 \approx \frac{1}{4\mathsf{H}[\rho]^2}\Big(1-2\frac{d}{dr}\mathsf{H}[\rho]\Big)
$$
is adopted. 
If we use this approximation, we can forget the dependence on $\lambda$ of $\mathfrak{N}^2$ and $\lambda_c$.
As D. O. Gough clarifies in \cite[p. 440]{Gough}, this approximation is done by `not taking the spherical geometry fully into account'. In other words, \eqref{G22}, \eqref{G38} `reduce to the above approximation if 
$\displaystyle \frac{1}{\mathsf{H}[g]}+\frac{2}{r} \rightarrow 0$ and
$\displaystyle \frac{1}{\mathsf{H}[\mathfrak{Q}]}\rightarrow\frac{1}{\mathsf{H}[\rho]}$'.
It is said that ` $\mathfrak{N}^2$ and $\lambda_c$ are approximated well except very close to the center of the star'.  So D. O.Gough calls this limits of 
$\lambda_c$ and $\mathfrak{N}^2$ the `planer values'.

On the other hand, the existence proof seriously depends on Assumption \ref{Ass.G1}. We do not know how to deal with the problem when $\mathscr{N}^2$ can take negative values somewhere on the background equilibrium. Namely what happens if we do not require that $d\bar{S}/dr > 0$ throughout $0<r<R$ but admit that 
$d\bar{S}/dr < 0$ so that $ \mathscr{N}^2 < 0$ for $0< R-r \ll 1$ in accordance with the more realistic view due to the asteroseismology? Actually it is asserted that the so called `zone en \'{e}quilibre convectif ' ( $\mathscr{N}^2 < 0$ but $|\mathscr{A}|$ is very small ) appears near the surface in many stellar models. See \cite{LedouxS}. \\

\subsection{ p-modes}

By the same way we can prove the existence of p-modes, say, the existence of a sequence $(\lambda_n)_{n\in\mathbb{N}}$ of positive eigenvalues such that
$\lambda_n \rightarrow +\infty$ as $ n \rightarrow \infty$. But we forget Assumption \ref{Ass.G1}. \\

Let us consider 
\begin{equation}
\mu=\frac{1}{\lambda}\leq \mu_0 =\frac{1}{\lambda^{\mathsf{p}}_0},
\end{equation}
$\lambda^{\mathsf{p}}_0$ being a large positive number, and $\mathsf{G}\leq \mathsf{G}^{\mathsf{p}}_0$,
$\mathsf{G}^{\mathsf{p}}_0$ being a arbitrarily fixed large number. \\

Instead of $E$ we use $E^{\mathsf{p}}$ defined by
\begin{equation}
E=-\Big(\frac{r}{g}\Big)^2\frac{\lambda^2}{l(l+1)}E^{\mathsf{p}},
\end{equation}
that is, 
\begin{equation}
A_{21}=-\rho\lambda E^{\mathsf{p}}.
\end{equation}
in \eqref{Aij} and
\begin{subequations}
\begin{align}
E^{\mathsf{p}}&=1+(E^{\mathsf{p}}_1+\mu E^{\mathsf{p}}_2)\mu, \\
E^{\mathsf{p}}_1&=\Big(\frac{r}{\mathsf{H}[g]}+2\Big)\frac{g}{r}=
4\frac{g}{r}-4\pi\mathsf{G}\rho, \\
E^{\mathsf{p}}_2&=-l(l+1)\Big(\frac{g}{r}\Big)^2.
\end{align}
\end{subequations}

If $\mu_0$ is sufficiently small, then $E^{\mathsf{p}}\geq 1-\varepsilon >0$ for
$0<r<R, 0\leq \mu \leq \mu_0$ with $0 <\varepsilon \ll 1$. 

Instead of $\mathfrak{Q}$ we define $\mathfrak{Q}^{\mathsf{p}}$ by
\begin{equation}
\mathfrak{Q}=-\frac{\lambda^2}{l(l+1)}\mathfrak{Q}^{\mathsf{p}},
\end{equation}
that is, 
\begin{equation}
\mathfrak{Q}^{\mathsf{p}}:=\frac{\rho}{r^2}E^{\mathsf{p}}.
\end{equation}
Then \eqref{G18} reads by using 
\begin{subequations}
\begin{align}
&\frac{1}{\mathsf{H}[\mathfrak{Q}]}=
\frac{1}{\mathsf{H}[\mathfrak{Q}^{\mathsf{p}}]} \\
&\mathfrak{N}^2=\mathscr{N}^2+g\Big(\frac{2}{\mathsf{H}[r/g]}+
\frac{1}{\mathsf{H}[E^{\mathsf{p}}]}\Big).
\end{align}
\end{subequations}
Since we do not suppose Assumption \ref{Ass.G1}, we do not claim that $\mathfrak{N}^2 >0$ but
we have
$$
\mathfrak{N}^2=
\begin{cases}
[r^2]_0 \quad\mbox{as}\quad r \rightarrow +0 \\
[R-r, (R-r)^{\nu+1}]_0\quad\mbox{as}\quad r \rightarrow R-0
\end{cases}
$$
and
there exists  a constant $M$ independent of $\mu \in [0,\mu_0]$ such that
$|\mathfrak{N}^2|\leq M$ for $0<r<R, 0\leq\mu \leq \mu_0$, $\mu_0$ being sufficiently small. 
 Now the equation \eqref{G24} holds valid, while \eqref{G24a}, \eqref{G24b}, \eqref{G24c} should read
\begin{subequations}
\begin{align}
A&=\frac{1}{\mathsf{H}[\rho/r^2]}+\frac{1}{\mathsf{H}[E^{\mathsf{p}}]}+
4\pi\mathsf{G}\frac{1}{\mathsf{H}[\rho]E^{\mathsf{p}}}\mu, \\
B&=\frac{1}{\mathsf{c}^2}\Big(\lambda+g\Big(\frac{1}{\mathsf{H}[g]}+\frac{2}{r}\Big)\Big)
-\frac{l(l+1)}{r^2}(1-\mathfrak{N}^2\mu)+ \nonumber \\
&+4\pi\mathsf{G}\rho\Big(
\frac{l(l+1)g}{\mathsf{H}[\rho]r^2E^{\mathsf{p}}}\mu^2+\frac{1}{\mathsf{c}^2}\Big), \\
K&=\Big(\frac{1}{\mathsf{H}[E^{\mathsf{p}}]}+
4\pi\mathsf{G}\frac{\rho}{\mathsf{H}[\rho]E^{\mathsf{p}}}\mu\Big)\frac{d}{dr}\delta\check{\Phi} + \nonumber \\
&+\Big(-\frac{E^{\mathsf{p}}_1}{g}+\frac{1}{\mathsf{H}[E^{\mathsf{p}}/g]}+
4\pi\mathsf{G}\frac{\rho}{\mathsf{H}[\rho]E^{\mathsf{p}}}\mu\Big)
\frac{l(l+1)g}{r^2}\mu\delta\check{\Phi}.
\end{align}
\end{subequations}

In \eqref{G30a}, \eqref{G30b}, \eqref{G31}, we replace
\begin{subequations}
\begin{align}
\frac{\rho}{\mathsf{H}[\rho]gE}&=-\frac{l(l+1)\rho g}{\mathsf{H}[\rho]E^{\mathsf{p}}r^2}\mu^2, \\
\frac{\rho }{\mathsf{H}Er}\frac{\mathsf{c}^2\lambda}{\mathscr{S}_l^2g^2}&=
-\frac{\rho}{\mathsf{H}[\rho]E^{\mathsf{p}}r}\mu.
\end{align}
\end{subequations}

Then it can be easily seen that, if  $\mu_0$ is sufficiently small, we can solve the equations for $X=\delta\check{\Phi}, \dot{X}=rd\delta\check{\Phi}/dr$ 
to define the operators $\mbox{\boldmath$\mathcal{H}$}_l, \dot{\mbox{\boldmath$\mathcal{H}$}}_l$
for $0 \leq \mu \leq \mu_0, 0\leq \mathsf{G}\leq\mathsf{G}^{\mathsf{p}}_0$. \\

We introduce the variable $\Psi^{\mathsf{p}}$ by
\begin{equation}
\eta=W^{\mathsf{p}}\Psi^{\mathsf{p}},
\end{equation}
where
\begin{align}
W^{\mathsf{p}}&=(\mathfrak{Q}^{\mathsf{p}})^{\frac{1}{2}}
\exp\Big[-2\pi\mathsf{G}\mu\int_0^r\frac{\rho}{\mathsf{H}[\rho]E^{\mathsf{p}}}\Big] \nonumber \\
&=\frac{\rho^{\frac{1}{2}}(E^{\mathsf{p}})^{\frac{1}{2}}}{r}
\exp\Big[-2\pi\mathsf{G}\mu\int_0^r\frac{\rho}{\mathsf{H}[\rho]E^{\mathsf{p}}}\Big].
\end{align}

We write the equation for $\Psi^{\mathsf{p}}$ as
\begin{equation}
-\frac{d^2\Psi^{\mathsf{p}}}{dr^2}+b^{\mathsf{p}}\Psi^{\mathsf{p}}+f^{\mathsf{p}} = \lambda \kappa^{\mathsf{p}}\Psi^{\mathsf{p}},
\end{equation}
where
\begin{subequations}
\begin{align}
b^{\mathsf{p}}&=\frac{1}{4}A^2+\frac{1}{2}\frac{dA}{dr}-B-\frac{\lambda}{\mathsf{c}^2} \nonumber \\
&=\frac{1}{4}A^2+\frac{1}{2}\frac{dA}{dr}
-\frac{g}{\mathsf{c}^2}\Big(\frac{1}{\mathsf{H}[g]}+\frac{2}{r}\Big)
-\frac{l(l+1)}{r^2}(1-\mathfrak{N}^2\mu) \nonumber \\
&-4\pi\mathsf{G}\Big(\frac{l(l+1)g}{\mathsf{H}[\rho]E^{\mathsf{p}}r^2}\mu^2+\frac{1}{\mathsf{c}^2}\Big), \\
f^{\mathsf{p}}&=-\frac{\rho K}{W^{\mathsf{p}}}, \\
\kappa^{\mathsf{p}}&=\frac{1}{\mathsf{c}^2}.
\end{align}
\end{subequations}

Considering $\mu=1/\lambda$ as a parameter we look at the eigenvalue problem for the eigenvalue $\Lambda^{\mathsf{p}}$:
\begin{equation}
-\frac{d^2\Psi^{\mathsf{p}}}{dr^2}+b^{\mathsf{p}}\Psi^{\mathsf{p}}+f^{\mathsf{p}}=\Lambda^{\mathsf{p}}\kappa^{\mathsf{p}}\Psi^{\mathsf{p}}. \label{P13}
\end{equation}
If $\Lambda^{\mathsf{p}}=1/\mu$, then this gives the eigenvalue $\lambda=1/\mu$ of the original problem. \\

By the Liouville transformation
\begin{align}
& x^{\mathsf{p}}=\int_0^r\sqrt{\kappa^{\mathsf{p}}}dr,\qquad y^{\mathsf{p}}=(\kappa^{\mathsf{p}})^{\frac{1}{4}}\Psi^{\mathsf{p}}, \nonumber \\
&\hat{f}^{\mathsf{p}}=(\kappa^{\mathsf{p}})^{-\frac{3}{4}}f^{\mathsf{p}}, \nonumber \\
&q^{\mathsf{p}}=\frac{b^{\mathsf{p}}}{\kappa^{\mathsf{p}}}+
\frac{1}{4}\frac{1}{\kappa^{\mathsf{p}}}\Big[
\frac{d^2}{dr^2}\log \kappa^{\mathsf{p}} -\frac{1}{4}\Big(\frac{d}{dr}
\log \kappa^{\mathsf{p}}\Big)^2\Big],
\end{align}
we transform \eqref{P13} to
\begin{equation}
-\frac{d^2y^{\mathsf{p}}}{dx^{\mathsf{p}2}}+q^{\mathsf{p}}y^{\mathsf{p}}+\hat{f}^{\mathsf{p}}=\Lambda^{\mathsf{p}}y^{\mathsf{p}}. \label{P14}
\end{equation}

We see
\begin{align*}
&x^{\mathsf{p}}_+:=\int_0^R\sqrt{\frac{1}{\mathsf{c}^2}}dr < \infty, \\
&x^{\mathsf{p}}=\sqrt{\frac{\rho_{\mathsf{O}}}{\gamma P_{\mathsf{O}}}}r(1+[r^2]_1)\quad\mbox{as}\quad r\rightarrow +0, \\
&x^{\mathsf{p}}_+-x^{\mathsf{p}}=\frac{2}{\mathsf{c}_R}\sqrt{R-r}(1+
[R-r,(R-r)^{\nu+1}]_1)\quad\mbox{as}\quad r \rightarrow R-0.
\end{align*}

The singularity of $q^{\mathsf{p}}$ turns out to be
$$
q^{\mathsf{p}}=
\begin{cases}
\displaystyle \frac{2}{(x^{\mathsf{p}})^2}(1+[(x^{\mathsf{p}})^2]_1)\quad\mbox{as}\quad x^{\mathsf{p}} \rightarrow +0 \\
\\
\displaystyle \frac{(2\nu+1)(2\nu+3)}{4}\frac{1}{(x^{\mathsf{p}}_+-x^{\mathsf{p}})^2}(1+
[(x^{\mathsf{p}}_+-x^{\mathsf{p}})^2, (x^{\mathsf{p}}_+-x^{\mathsf{p}})^{2(\nu+1)}]_1) \\
\quad\mbox{as}\quad x^{\mathsf{p}} \rightarrow x^{\mathsf{p}}_+.
\end{cases}
$$

Note that $$2 > \frac{3}{4}, \qquad\frac{(2\nu+1)(2\nu+3)}{4}>\frac{15}{4}>\frac{3}{4}.$$
We can claim that there exist constants
 $\displaystyle K^{\mathsf{p}}_0, K^{\mathsf{p}}_1, \frac{3}{4}<K^{\mathsf{p}}_1 < 2$ independent of $\mu, \mathsf{G}$ such that
$$\tilde{q}^{\mathsf{p}}(x^{\mathsf{p}},; \mu)=q^{\mathsf{p}}+K^{\mathsf{p}}_0 \geq
1+K^{\mathsf{p}}_1\Big(
\frac{1}{(x^{\mathsf{p}})^2} \vee \frac{1}{(x^{\mathsf{p}}_+-x^{\mathsf{p}})^2}\Big).
$$
Moreover
$$\tilde{q}^{\mathsf{p}}(x^{\mathsf{p}};\mu)=
\tilde{q}^{\mathsf{p}}(x^p;0)
(1+\Omega_q^{\mathsf{p}}(x^{\mathsf{p}},\mu)\mu) $$
with
$$|\Omega_q^{\mathsf{p}}(x^{\mathsf{p}},\mu)|\leq M_q^{\mathsf{p}},\quad
\Big|\frac{\partial}{\partial \mu}\Omega^{\mathsf{p}}_q(x^{\mathsf{p}},\mu)\Big|\leq M_q^{\mathsf{p}}
$$
for $0 < x< x^{\mathsf{p}}_+, 0\leq\mu\leq \mu_0, 0\leq \mathsf{G}\leq \mathsf{G}^{\mathsf{p}}_0$ with
a constant $M_q^{\mathsf{p}}$ such that
$\varepsilon:=M^{\mathsf{p}}_q\mu_0$ satisfies
$$0<1-\varepsilon,\qquad
\frac{3}{4} <(1-\varepsilon )K^{\mathsf{p}}_1.$$\\

We put
\begin{align}
&Q_0^{\mathsf{p}}[y^{\mathsf{p}};\mu]:=\int_0^{x^{\mathsf{p}}_+}
\Big(\Big|\frac{dy^{\mathsf{p}}}{dx^{\mathsf{p}}}\Big|^2+\tilde{q}^{\mathsf{p}}(x^{\mathsf{p}};\mu)|y^{\mathsf{p}}|^2\Big)dx^{\mathsf{p}}, \\
&\mathfrak{E}^{\mathsf{p}}:=L^2([0,x^{\mathsf{p}}_+], dx^{\mathsf{p}}) \\
&\mathfrak{E}^{\mathsf{p}}_1:= \{ y^{\mathsf{p}}\in \mathfrak{E}^{\mathsf{p}}\  |\  Q_0^{\mathsf{p}}[y^{\mathsf{p}};0] <\infty \}
\end{align}
and
\begin{equation}
Q^{\mathsf{p}}[y^{\mathsf{p}};\mu]:=
\int_0^{x^{\mathsf{p}}_+}\Big(
\Big|\frac{dy^{\mathsf{p}}}{dx^{\mathsf{p}}}\Big|^2+
({q}^{\mathsf{p}}+K^{\mathsf{p}}_0)|y^{\mathsf{p}}|^2
+\hat{f}^{\mathsf{p}}(y^{\mathsf{p}})^*\Big)dx^{\mathsf{p}}.
\end{equation}\\

If $y^{\mathsf{p}} \in \mathfrak{E}^{\mathsf{p}}_1$, then 
$Y \in L^1(x^{\mathsf{p}}dx^{\mathsf{p}})$ and
\begin{subequations}
\begin{align}
&\|Y\|_{L^1(x^{\mathsf{p}}dx^{\mathsf{p}})}\leq M^{\mathsf{p}}_Y\sqrt{Q^{\mathsf{p}}_0(y^{\mathsf{p}};0]}, \\
&\|Y|_{\mu+\delta\mu}-Y|_{\mu}\|_{L^1(x^{\mathsf{p}}dx^{\mathsf{p}})}\leq
M^{\mathsf{p}}_Y|\delta\mu|\sqrt{Q^{\mathsf{p}}_0[y^{\mathsf{p}};0]}.
\end{align}
\end{subequations}
Here
\begin{align}
Y&=-\frac{W^{\mathsf{p}}(\kappa^{\mathsf{p}})^{\frac{1}{4}}}{\mathsf{H}[\rho]E^{\mathsf{p}}}\Big[\mu\frac{dy^{\mathsf{p}}}{dx^{\mathsf{p}}}+ \nonumber \\
&+\Big(\frac{1}{4}\frac{1}{\mathsf{H}[\kappa^{\mathsf{p}}]}\mu+
\frac{l(l+1)g}{r^2}\mu^2+
\frac{\mathsf{H}[\rho]E^{\mathsf{p}}}{\mathsf{H}[W^{\mathsf{p}}]}+
\frac{\mathsf{H}[\rho]E^{\mathsf{p}}}{\mathsf{c}^2}\Big)
(\kappa^{\mathsf{p}})^{-\frac{1}{2}}y^{\mathsf{p}}\Big]. \label{Yp}
\end{align}\\

Consequently, if $y^{\mathsf{p}} \in \mathfrak{E}^{\mathsf{p}}_1$, then 
$x^{\mathsf{p}}\hat{f}^{\mathsf{p}}\in L^{\infty}$ and we have estimates by
$\sqrt{ Q^{\mathsf{p}}_0[y^{\mathsf{p}};0] }$ and so on. And
we have, for $y^{\mathsf{p}}\in \mathfrak{E}^{\mathsf{p}}_1$,
\begin{subequations}
\begin{align}
&(1-\varepsilon)Q^{\mathsf{p}}_0[y^{\mathsf{p}};0]\leq
Q^{\mathsf{p}}[y^{\mathsf{p}};\mu]\leq
(1+\varepsilon)Q^{\mathsf{p}}_0[y^{\mathsf{p}};0], \\
&|Q^{\mathsf{p}}[y^{\mathsf{p}};\mu+\delta\mu]-
Q^{\mathsf{p}}[y^{\mathsf{p}};\mu]|\leq
M^{\mathsf{p}}_Q|\delta\mu|Q^{\mathsf{p}}_0[y^{\mathsf{p}}; 0].
\end{align}
\end{subequations}

Here, by estimating the coefficients of \eqref{Yp}, we see that
 it is sufficient to take $\mu_0$ sufficiently small, and no other restrictions are necessary. 
\\

Therefore the eigenvalue problem \eqref{P14} is of the Sturm-Liouville type, and we have a sequence of simple eigenvalues $(\Lambda^{\mathsf{p}}_n(\mu))_{n\in\mathbb{N}}$ such that
$$\Lambda^{\mathsf{p}}_1(\mu)<\Lambda^{\mathsf{p}}_2(\mu)<\cdots<
\Lambda^{\mathsf{p}}_n(\mu)<\cdots \rightarrow +\infty.$$

Let $n^{\mathsf{p}}_0$ be such that $\displaystyle \Lambda^{\mathsf{p}}_n(\mu_0) >\frac{1}{\mu_0}$ for
$n \geq n^{\mathsf{p}}_0$. Since 
$$\varphi^{\mathsf{p}}\  :\  \mu \mapsto \mu-\frac{1}{\Lambda^{\mathsf{p}}_n(\mu)}$$
is continuous and
$$\varphi^{\mathsf{p}}(0)=-\frac{1}{\Lambda^{\mathsf{p}}_n(0)} <0,\qquad
\varphi^{\mathsf{p}}(\mu_0) =\mu_0-\frac{1}{\Lambda^{\mathsf{p}}_n(\mu_0)} >0,
$$
there exists $\mu_n \in ]0,\mu_0[$ such that 
$\displaystyle \mu_n=\frac{1}{\Lambda^{\mathsf{p}}_n(\mu_n)}$. Then 
$\displaystyle \lambda_n=\frac{1}{\mu_n}$ forms a sequence of eigenvalues such that 
$\lambda_n \rightarrow +\infty$ as $ n\rightarrow \infty$. Summing up, we have

\begin{Theorem} \label{Th.P}
Let $l \geq 1$. Suppose Assumptions \ref{Ass.G2}, \ref{Ass.G3}.
Then there exists a sequence 
$(\lambda_n)_{n\in\mathbb{N}}$ of positive eigenvalues of $\vec{L}_l$ such that
$\lambda_n \rightarrow +\infty$ as $n \rightarrow \infty$.
\end{Theorem}

\section{Formulation as the first order 4-dimensional system of ordinary differential equations.}

We are considering the eigenvalue problem
\begin{equation}
\vec{L}_l\vec{V}=\lambda\vec{V}. \label{S01}
\end{equation}
In this section we keep supposing $l \geq 1$, and we consider $\lambda\not=0$. 
Also we consider the inhomogeneous problem
\begin{equation}
(\vec{L}_l-\lambda)\vec{V}=\vec{f}, \label{S02}
\end{equation}
where 
\begin{equation}
\vec{f}=\begin{bmatrix}
f^r \\
\\
f^h
\end{bmatrix}
\in \mathfrak{X}_l. \label{S03}
\end{equation}
In this section we use a formulation of he problem as a first order system of
ordinary differential equations. We want to examine the argument by J. Eisenfeld, 1969, \cite{Eisenfeld}.\\

We introduce the variables
\begin{subequations}
\begin{align}
& y_1=V^r \\
& y_2=\frac{1}{r}(\delta\check{P}-\rho gV^r) \\
& y_3=\frac{1}{r}\delta\check{\Phi} \\
& y_4=\frac{d}{dr}\check{\Phi}.
\end{align}
\end{subequations}

Keeping in mind that 
$$\delta\check{\Phi}=-4\pi\mathsf{G}\mathcal{H}_l(\delta\check{\rho})
$$
implies
$$\Big[\frac{d^2}{dr^2}+\frac{2}{r}\frac{d}{dr}-
\frac{l(l+1)}{r^2}\Big]\delta\check{\Phi}=4\pi\mathsf{G}\delta\check{\rho},
$$ we can derive the equation
\begin{equation}
r\frac{d\vec{y}}{dr}=A(r,\lambda)\vec{y}+\vec{h}, \label{S05}
\end{equation}
where 
\begin{align}
&\vec{y}=
\begin{bmatrix}
y_1 \\
y_2 \\
y_3 \\
y_4
\end{bmatrix},
\qquad A(r,\lambda)=
\begin{bmatrix}
a_{11} & a_{12} & a_{13} & a_{14} \\
a_{21} & a_{22} & a_{23} & a_{24} \\
a_{31} & a_{32} & a_{33} & a_{34} \\
a_{41} & a_{42} & a_{43} & a_{44} \\
\end{bmatrix}, \\
&\vec{h}=
\begin{bmatrix}
h_1 \\
h_2 \\
h_3 \\
h_4
\end{bmatrix}
=
\begin{bmatrix}
-Lf^h \\
\displaystyle \rho\Big(f^r+L\frac{g}{r}f^h\Big) \\
0 \\
0
\end{bmatrix},
\end{align}
while here and hereafter we denote
\begin{equation}
L=\frac{l(l+1)}{\lambda}
\end{equation}
and
\begin{align}
& a_{11}=-2+L\frac{g}{r},
\quad a_{12}=(-r^2+L\mathsf{c}^2)\frac{1}{\mathsf{c}^2\rho},
\quad a_{13}=L,
\quad a_{14}=0, \nonumber \\
&a_{21}=\rho\Big[-L\Big(\frac{g}{r}\Big)^2+g\Big(\frac{1}{\mathsf{H}[g]}+\frac{2}{r}\Big)+\lambda\Big], \quad 
a_{22}=-L\frac{g}{r}-1,
\nonumber \\
&a_{23}=-L\rho\frac{g}{r},
\quad a_{24}=-\rho
\nonumber \\
&a_{31}=a_{32}=0,\quad a_{33}=-1,\quad a_{34}=1, \nonumber \\
&a_{41}=4\pi\mathsf{G}r\frac{\rho}{\mathsf{H}[\rho]},
\quad
a_{42}=4\pi\mathsf{G}\frac{r^2}{\mathsf{c}^2}, \quad a_{43}=l(l+1),
\quad a_{44}=-2.
\end{align}\\

The boundary conditions at $r=0$, which corresponds to $\vec{V}\in 
\overset{\circ}{\mathfrak{W}}_l $ and $\delta\check{\Phi}=-4\pi\mathsf{G}\mathcal{H}_l(\delta\check{\rho})$, are satisfied when
\begin{subequations}
\begin{align}
&\int_0^{R/2}|y_1|^2 r^2dr <\infty,\quad
\int_0^{R/2}|y_2|^2r^4dr <\infty, \\
&y_3=O(r^{1}),\quad y_4=O(r^{-\frac{1}{2}}),
\end{align}
\end{subequations}
and
the boundary conditions at $r=R$ are satisfied when
\begin{subequations}
\begin{align}
&\int_{R/2}^R|y_1|^2\rho dr <\infty,\quad
\int_{R/2}^R|y_2|^2\frac{1}{\mathsf{c}^2\rho}dr <\infty, \\
&y_1=O((R-r)^{\kappa}) \quad\mbox{with}\quad \kappa >-\frac{\nu}{2} \\
&y_3(R)=\lim_{r\rightarrow R-0}y_3(r)
\quad \mbox{and}\quad
y_4(R)=\lim_{r\rightarrow R-0}y_4(r)
\quad\mbox{exist} \nonumber \\
\mbox{and}& \nonumber \\
&(l+1)y_3(R)+y_4(R)=0. \label{S11c}
\end{align}
\end{subequations}

Here we note that $\vec{V}\in \mathfrak{W}_l$ if and only if
$\vec{V} \in \mathfrak{X}_l$ and
$\displaystyle \delta\check{\rho} \in L^2\Big(\frac{\mathsf{c}^2}{\rho}r^2dr\Big)$, while
$$\delta\check{\rho}=-\frac{d\rho}{dr}y_1+\frac{r}{\mathsf{c}^2}y_2 $$
and
$\displaystyle -\frac{d\rho}{dr}y_1 \in  L^2\Big(\frac{\mathsf{c}^2}{\rho}r^2dr\Big) $
if $y_1 \in L^2(\rho r^2dr)$. Recall Proposition \ref{EqVxieta}, while 
$$\xi=y_1, \quad \eta=r y_2, \quad  \delta\check{\Phi}=r y_3 $$
satisfy \eqref{AX2} if $y_1 \in L^2(\rho r^2dr), y_2\in L^2(r^2dr / \mathsf{c}^2\rho) \subset
L^2(r^2dr/\rho), y_3 \in L^2(\rho r^2dr)$.

Also we note that \eqref{S11c}, 
which means 
$$ r\frac{d}{dr}\delta\check{\Phi}=-(l+1)\delta\check{\Phi}\quad\mbox{at}\quad r=R, $$ 
comes from that $\delta\check{\Phi}=-4\pi\mathsf{G}\mathcal{H}_l(\delta\check{\rho})$ should satisfy
 $$
\delta\check{\Phi}(r)=\frac{C}{r^{l+1}} \quad\mbox{for}\quad r \geq R.
$$
Otherwise, the corresponding $\delta\check{\Phi}$
for which $\displaystyle r y_3$ is the candidate
might be equal to  $-4\pi\mathsf{G}\mathcal{H}_l(\delta\check{\rho})+C'r^{l}$ with $C'\not=0$.\\

\subsection{Eigenvalue problem}

We are going to analyze the homogeneous problem, say, the eigenvalue problem
\begin{equation}
r\frac{d\vec{y}}{dr}=A(r,\lambda)\vec{y}. \label{6.3}
\end{equation}\\

Now let us suppose  Assumption \ref{Ass.G2} and consider the system \eqref{6.3} at $r=+0$. Under the Assumption \ref{Ass.G2}, we see
\eqref{6.3} reads
\begin{equation}
r\frac{d\vec{y}}{dr}=\Big(K_0+\sum_{m\geq 1}r^{2m}K_m\Big)\vec{y}, \label{6.10}
\end{equation}
where
\begin{equation}
K_0=
\begin{bmatrix}
-2+Lg_O & \frac{L}{\rho_{\mathsf{O}}} & L & 0 \\
\rho_{\mathsf{O}}(-Lg_O^2+g_O+\lambda) & -Lg_O-1 & -L\rho_{\mathsf{O}}g_O &
 -\rho_{\mathsf{O}} \\
0 & 0 & -1 & 1 \\
0 & 0 & l(l+1) & -2
\end{bmatrix}
\end{equation}
and
$\sum r^{2m}K_m$ is a convergent matrix-valued power series in $r^2$ with positive radius of convergence. Here $$g_O:=\frac{P_{\mathsf{O}1}}{\rho_{\mathsf{O}}}=\lim_{r\rightarrow +0}\frac{g}{r},
$$ while
$$
P_{\mathsf{O}1}=-\lim_{r\rightarrow +0}\frac{1}{r}\frac{dP}{dr}.
$$
The eigenvalues of $K_0$ are $l-1, -(l+2)$, which are double. 

Thus, putting
\begin{equation}
z=r^2,\quad \vec{y}=
\begin{bmatrix}
L & 0 & L & 0 \\
-\rho_{\mathsf{O}}(Lg_O-l-1) & -\rho_{\mathsf{O}} & -\rho_{\mathsf{O}}(Lg_O+l) & 
-\rho_{\mathsf{O}} \\
0 & 1 & 0 & 1 \\
0 & l & 0 & -(l+1)
\end{bmatrix}
\vec{w},
\end{equation}
we have a system
\begin{equation}
z\frac{d\vec{w}}{dz}=\Big(\mathbf{R}+
\sum_{m\geq 0}z^{m+1}A_m\Big)\vec{w} \label{6.14}
\end{equation}
with
\begin{equation}
\mathbf{R}=\mathrm{diag}(\rho_1, \rho_1, \rho_2, \rho_2),
\end{equation}
where 
\begin{equation} \rho_1=\frac{l-1}{2}, \qquad \rho_2=-\frac{l+2}{2}.
\end{equation}
Since $\rho_1-\rho_2=l+\frac{1}{2}$ is not an integer, 
\cite[Chapter 4, Theorem 4.1]{CoddingtonL} gives a fundamental matrix $\Phi_w$ of the system \eqref{6.14} of the form
\begin{equation}
\Phi_w(z)=\Big(I_4+\sum_{m\geq 1}z^mP_m\Big)z^{\mathbf{R}},
\end{equation}
where $\sum z^mP_m $ is a convergent matrix-valued power series. 
Note that
$$ z^{\mathbf{R}}=\mathrm{diag}(z^{\rho_1}, z^{\rho_1}, z^{\rho_2}, z^{\rho_2}).
$$

As result, we have a fundamental system of solutions 
\begin{align}
\mbox{\boldmath$\Phi$}_O&=
\begin{bmatrix}
L & 0 & L & 0 \\
-\rho_{\mathsf{O}}(Lg_O-l-1) & -\rho_{\mathsf{O}} & -\rho_{\mathsf{O}}(Lg_O+l) & 
-\rho_{\mathsf{O}} \\
0 & 1 & 0 & 1 \\
0 & l & 0 & -(l+1)
\end{bmatrix}
(r^2)^{\mathbf{R}} \nonumber \\
&(I_4+\sum_{m\geq 1}z^mP_m) \nonumber \\
&=[ \mbox{\boldmath$\varphi$}_{O1} \quad \mbox{\boldmath$\varphi$}_{O2} \quad  \mbox{\boldmath$\varphi$}_{O3} \quad \mbox{\boldmath$\varphi$}_{O4}]. \label{60019}
\end{align}

We see that
\begin{subequations}
\begin{align}
&\mbox{\boldmath$\varphi$}_{O1}=
\begin{bmatrix}
L+r[r^2]_0 \\
\alpha+ r[r^2]_0 \\
[r^2]_1 \\
[r^2]_1
\end{bmatrix}
r^{l-1}, \label{60020a} \\
&\mbox{\boldmath$\varphi$}_{O2} =
\begin{bmatrix}
r[r^2]_0 \\
-\rho_{\mathsf{O}}+r[r^2]_0 \\
1+[r^2]_1 \\
l+[r^2]_1
\end{bmatrix}
r^{l-1}, \label{60020b} \\
&\mbox{\boldmath$\varphi$}_{O3}=
\begin{bmatrix}
L+[r^2]_1+r^{2l+1}[r^2]_1 \\
\beta+[r^2]_1+r^{2l+1}[r^2]_1 \\
[r^2]_1+r^{2l+1}[r^2]_1 \\
[r^2]_1+r^{2l+1}[r^2]_1
\end{bmatrix}
r^{-(l+2)}, \label{60020c} \\
&\mbox{\boldmath$\varphi$}_{O4}=
\begin{bmatrix}
[r^2]_1+r^{2l+1}[r^2]_1 \\
-\rho_{\mathsf{O}}+[r^2]_1+r^{2l+1}[r^2]_1 \\
1+[r^2]_1+r^{2l+1}[r^2]_1 \\
-(l+1)+[r^2]_1+r^{2l+1}[r^2]_1
\end{bmatrix}
r^{-(l+2)} \label{60020d}
\end{align}
\end{subequations}
Here we denote
\begin{equation}
\alpha=-\rho_{\mathsf{O}}(Lg_O-l-1),\qquad
\beta=-\rho_{\mathsf{O}}(Lg_O+l).
\end{equation}

Here \eqref{60020a}\eqref{60020b} are not obvious, since \eqref{60019}
tels us 
\begin{subequations}
\begin{align}
&\mbox{\boldmath$\varphi$}_{O1}=\begin{bmatrix}
Lr^{l-1}+Lp_{11}r^{l-1}+Lp_{31}r^{-(l+2)}, \\
\alpha r^{l-1}+\alpha p_{11}r^{l-1}+\beta p_{31}r^{-(l+2)}, \\
p_{21} r^{l-1}+p_{41}r^{-(l+2)}, \\
lp_{21}r^{l-1}-(l+1)p_{41}r^{-(l+2)}
\end{bmatrix}
, \\
&\mbox{\boldmath$\varphi$}_{O2}=\begin{bmatrix}
Lp_{12}r^{l-1}+Lp_{32}r^{-(l+2)}, \\
-\rho_{\mathsf{O}}+\alpha p_{12}r^{l-1}+\beta p_{32}r^{-(l+2)}, \\
r^{l-1}+p_{22} r^{l-1}+p_{42}r^{-(l+2)}, \\
lr^{l-1}+lp_{22}r^{l-1}-(l+1)p_{42}r^{-(l+2)}
\end{bmatrix}, \\
&\mbox{\boldmath$\varphi$}_{O3}=\begin{bmatrix}
Lr^{-(l+2)}+Lp_{13}r^{l-1}+Lp_{33}r^{-(l+2)}, \\
\beta r^{-(l+2)}+\alpha r^{l-1}+\alpha p_{13}r^{l-1}+\beta p_{33}r^{-(l+2)}, \\
p_{23} r^{l-1}+p_{43}r^{-(l+2)}, \\
lp_{23}r^{l-1}-(l+1)p_{43}r^{-(l+2)}
\end{bmatrix}
,
\\
&\mbox{\boldmath$\varphi$}_{O4}=\begin{bmatrix}
Lp_{14}r^{l-1}+Lp_{34}r^{-(l+2)}, \\
-\rho_{\mathsf{O}}+\alpha p_{14}r^{l-1}+\beta p_{34}r^{-(l+2)}, \\
r^{-(l+2)}+p_{21} r^{l-1}+p_{41}r^{-(l+2)}, \\
-(l+1)r^{-(l+2)}+lp_{21}r^{l-1}-(l+1)p_{41}r^{-(l+2)}
\end{bmatrix}
.
\end{align}
\end{subequations}
Here we denote
$$\sum_{m\geq 1}r^{2m}P_m=(p_{ij}),\qquad
p_{ij}=p_{ij}(r)=[r^2]_1.
$$

In fact, we can show \eqref{60020a} as follows:

For $\vec{y}=\mbox{\boldmath$\varphi$}_{O1}$, \eqref{60019}  implies
\begin{align*}
&y_1=Lr^{l-1}+Lp_{11}r^{l-1}+Lp_{31}r^{-(l+2)}, \\
&y_2=\alpha r^{l-1}+\alpha p_{11}r^{l-1}+\beta p_{31}r^{-(l+2)}, \\
&y_3=p_{21} r^{l-1}+p_{41}r^{-(l+2)}, \\
&y_4=lp_{21}r^{l-1}-(l+1)p_{41}r^{-(l+2)}
\end{align*}

We claim $p_{41}=0$. Otherwise, suppose $p_{41}=cr^{2k}(1+[r^2]_1), k \geq 1, c\not=0$. 
Looking at the equation
$$ r\frac{dy_3}{dr}=-y_3+y_4, $$
we see
$$r\frac{dy_3}{dr}\sim c(2k-(l+2))r^{2k-(l+2)},
$$ and
$$-y_3+y_4 \sim (-c-(l+1))r^{2k-(l+2)},
$$ 
provided that $2k-(l+2)\not=0$. Then $2k-(l+2)=-l-(l+1)$, that is, $k=0$,
provided that $2k-(l+2)\not=0$. A contradiction. Hence $2k-(l+2)=0$
so that $p_{41}=cr^{l+2}(1+[r^2]_1)$. Then
\begin{align*}
&y_3=p_{21}r^{l-1}+c(1+r^2]_1), \\
&y_4=lp_{21}r^{l-1}-(l+1)c(1+[r^2]_1).
\end{align*}
Then the equation
$$ r\frac{dy_3}{dr}=-y_3+y_4 $$
deduces $ (-1-(l+1))c=0$, say, $c=0$, a contradiction. Summing up, we can claim $p_{41}=0$
and
\begin{align*}
&y_1=Lr^{l-1}+Lp_{11}r^{l-1}+Lp_{31}r^{-(l+2)}, \\
&y_2=\alpha r^{l-1}+\alpha p_{11}r^{l-1}+\beta p_{31}r^{-(l+2)}, \\
&y_3=p_{21} r^{l-1}, \\
&y_4=lp_{21}r^{l-1}.
\end{align*}
Suppose $p_{31}\not=0$. Then $p_{31}=cr^{2k}(1+[r^2]_1),
k \geq 1, c\not=0$. Look at the equation
$$r\frac{dy_1}{dr}=a_{11}y_1+a_{12}y_2+L y_3. $$
We can see that, if $l-1>2k-(l+2)$, then $2k-(l+2)=-(l+2)$, say, $k=0$, a contradiction. 
Therefore $l-1 \leq 2k-(l+2)$, say, $k \geq l+1$, and $p_{31}=r^{2l+2}[r^2]_0$,
so that $p_{31}r^{-(l+2)}=r[r^2]_0r^{l-1}$. Of course this is the case if $p_{31}=0$.
Thus we can claim \eqref{60020a}.

In the same way we can show $p_{42}=0, p_{32}=r^{2l+2}[r^2]_0$ and \eqref{60020b}. In order to claim \eqref{60020c}\eqref{60020d}, it is sufficient to note that $ 
r^{l-1}=r^{2l+1}\cdot r^{(l+2)}$.\\

It is easy to see that only $\mbox{\boldmath$\varphi$}_{O1}, \mbox{\boldmath$\varphi$}_{O2}$ satisfy the boundary conditions.
Therefore we have
\begin{equation}
\vec{y}=C_{01}\mbox{\boldmath$\varphi$}_{O1}+C_{02}\mbox{\boldmath$\varphi$}_{O2},
\end{equation}
with constants $C_{01}, C_{02}$
in order that $\vec{y}$ gives $\vec{V} \in \mathsf{D}(\vec{L}_l)$. \\

Let us consider the system \eqref{6.3} at the boundary point $r=R-0$. Here we suppose
 Assumption \ref{Ass.G3} and 
\begin{Assumption}\label{Ass.nuQ}
 The index $\displaystyle \nu=\frac{1}{\gamma-1}$ is a rational number. Let $\nu=N/D, N, D $ being mutually prime natural numbers.
\end{Assumption}

We use the variable $s$ defined by
\begin{equation}
s=(R-r)^{1/D} \qquad \Leftrightarrow \quad R-r=s^D.
\end{equation}

Using the variable $\vec{Y}=(Y_1,Y_2, Y_3, Y_4)^{\top}$ defined by
\begin{equation}
y_1=Y_1,\quad y_2=s^NY_2,\quad y_3=Y_3,\quad y_4=Y_4.
\end{equation}
the system \eqref{6.3}
reads
\begin{equation}
s\frac{d\vec{Y}}{ds}=\tilde{A}(s,\lambda)\vec{Y}, \label{6024}
\end{equation}
where
\begin{equation}
\tilde{A}(s,\lambda)=
\begin{bmatrix}
\tilde{a}_{11} & \tilde{a}_{12} &\tilde{a}_{13} &\tilde{a}_{14}  \\
\tilde{a}_{21} &\tilde{a}_{22} &\tilde{a}_{23} &\tilde{a}_{24} \\
\tilde{a}_{31} &\tilde{a}_{32} &\tilde{a}_{33} &\tilde{a}_{34} \\
\tilde{a}_{41} &\tilde{a}_{42} &\tilde{a}_{43} &\tilde{a}_{44} 
\end{bmatrix}
\end{equation}
with
\begin{align}
&\tilde{a}_{11}=-\frac{D}{R-s^D}s^Da_{11}, \quad
\tilde{a}_{12}=-\frac{D}{R-s^D}s^{D+N}a_{12}, \nonumber \\
&\tilde{a}_{13}=-\frac{D}{R-s^D}s^Da_{13}, \quad
\tilde{a}_{14}=-\frac{D}{R-s^D}s^Da_{11},  \nonumber \\
&\tilde{a}_{21}=-\frac{D}{R-s^D}s^{D-N}a_{21}, \quad
\tilde{a}_{22}=-N-\frac{D}{R-s^D}s^Da_{22}, \nonumber \\
&\tilde{a}_{23}=-\frac{D}{R-s^D}s^{D-N}a_{23}, \quad
\tilde{a}_{24}=-\frac{D}{R-s^D}s^{D-N}a_{11}, \nonumber \\
&\tilde{a}_{31}=-\frac{D}{R-s^D}s^Da_{31}, \quad
\tilde{a}_{32}=-\frac{D}{R-s^D}s^{D+N}a_{32}, \nonumber \\
&\tilde{a}_{33}=-\frac{D}{R-s^D}s^Da_{33}, \quad
\tilde{a}_{34}=-\frac{D}{R-s^D}s^Da_{34}, \nonumber \\
&\tilde{a}_{41}=-\frac{D}{R-s^D}s^{D}a_{41}, \quad
\tilde{a}_{42}=-\frac{D}{R-s^D}s^{D+N}a_{42}, \nonumber \\
&\tilde{a}_{43}=-\frac{D}{R-s^D}s^Da_{43}, \quad
\tilde{a}_{44}=-\frac{D}{R-s^D}s^Da_{44}. 
\end{align}

We see that
\begin{equation}
\tilde{A}(s, \lambda)=K_0^R+\sum_{m\geq 1}s^mK_m^R,
\end{equation}
where
\begin{equation}
K_0^R=
\begin{bmatrix}
0 & \sigma & 0 & 0 \\
0 & -N & 0 & 0 \\
0 & 0 & 0 & 0 \\
0 & 0 & 0 & 0
\end{bmatrix}.
\end{equation}
Here
\begin{equation}
\sigma:=\frac{DR}{\mathsf{c}_R^2C_{\rho}}
\end{equation}

Introducing $\vec{w}$ by
\begin{equation}
\vec{Y}=
\begin{bmatrix}
1 & 0 & 0 & \sigma \\
0 & 0 & 0 & -N \\
0 & 1 & 0 & 0 \\
0 & 0 & 1 & 0
\end{bmatrix}
\vec{w},
\end{equation}
the equation \eqref{6024} reads
\begin{equation}
s\frac{d\vec{w}}{ds}=(\mathbf{R}^R+\sum_{m\geq 1}s^mA_m^R)\vec{w} \label{6031}
\end{equation}
with 
\begin{equation}
\mathbf{R}^R=\mathrm{diag}(0, 0, 0, -N).
\end{equation}
Hence, applying the recipe prescribed in the proof of
\cite[p.120, Chapter 4, Lemma]{CoddingtonL}, we have a fundamental matrix $\Phi_w^R$
of the equation \eqref{6031} of the form
\begin{equation}
\Phi_w^R=(I+\sum_{m\geq 1}s^mP_m^R)\mathrm{diag}(1,1,1,s^{-N}),
\end{equation}
and the fundamental matrix $\Phi_Y^R$ of the equation
\eqref{6024} of the form
\begin{align}
\Phi_Y^R&=
\begin{bmatrix}
1 & 0 & 0 & \sigma \\
0 & 0 & 0  & -N \\
0 & 1 & 0 & 0 \\
0 & 0 & 1 & 0
\end{bmatrix}
\Phi_w^R \nonumber \\
&=
\begin{bmatrix}
1 & 0 & 0 & \sigma s^{-N} \\
0 & 0 & 0 & Ns^{-N} \\
0 & 1 & 0 & 0 \\
0 & 0 & 1 & 0
\end{bmatrix}
(I+\sum_{s\geq 1}s^mP_m^R). \label{6034}
\end{align}
This gives the fundamental matrix $\mbox{\boldmath$\Phi$}_R$ of the equation
\eqref{6.3}
\begin{equation}
\mbox{\boldmath$\Phi$}_R=
\mathrm{diag}(1,s^N, 1,1)\Phi_Y^R=
[\mbox{\boldmath$\varphi$}_{R1} \quad \mbox{\boldmath$\varphi$}_{R2} \quad \mbox{\boldmath$\varphi$}_{R3} \quad
\mbox{\boldmath$\varphi$}_{R4}]. \label{6035}
\end{equation}\\

\begin{Remark}
The above reduction of \eqref{6024} to \eqref{6031} seriously depends on the assumption that $\displaystyle \nu=\frac{1}{\gamma-1}$ be rational. If  $\nu$ is irrational, it may be impossible to find a variable $s$ by which the Briot-Bouquet type singularity at $s=0 \Leftrightarrow R-r=0$ turns out to be of the canonical form like \eqref{6031},
in which the principal part is diagonal and the remainder terms are analytic in $s$.  In fact the coefficients
$$a_{11}=-2+L\frac{g}{r},\quad a_{12}=(-r^2+L\mathsf{c}^2)\frac{1}{\mathsf{c}^2\rho},\quad\mbox{etc}
$$
have singularity of the form
$$\Psi(R-r,(R-r)^{\nu+1})=\sum_{j_1,j_2}c_{j_1j_2}(R-r)^{j_1}(R-r)^{(\nu+1)j_2} $$
with $\frac{\partial}{\partial X_1}\Psi(X_1,X_2)\not=0, \frac{\partial}{\partial X_2}\Psi(X_1,X_2)\not=0$, generally. Therefore, when $\nu$ is irrational,  the singularity is essentially transcendental (not algebraic) and we may be unable to find a variable by which coefficients turn out to be analytic simultaneously. Even when $\nu$ is rational, Assumption \ref{Ass.G3} is inevitable.  In fact neither  transformation of the Briot-Bouquet type singularity of \eqref{6024} to the canonical form \eqref{6031} nor specification of asymptotic behavior of the elements of the fundamental matrix of solutions is known when $\tilde{A}(s,\lambda)=K_0^R+\mathbf{K}(s)$
with $\mathbf{K}(\cdot) \in C^{\infty}, \mathbf{K}(0)=O$, which is smooth but can be not analytic, since the eigenvalues of $K_0^R$ are $0, 0, 0, -N(<0, \in \mathbb{Z})$. The theory of  diagonalizing and asymptotic behavior of solutions usually assumes that the remainder terms are analytic, except for the case with all positive, or all  negative eigenvalues, that is, the case of a strict source or a sink. 
In this sense the discussion of  \cite{Eisenfeld} seems to be too naive. 
\end{Remark}

We see that
\begin{subequations}
\begin{align}
&\mbox{\boldmath$\varphi$}_{R1}=
\begin{bmatrix}
1+[s]_1 \\
s^N[s]_1 \\
0\\
0
\end{bmatrix},
\label{6036a} \\
&\mbox{\boldmath$\varphi$}_{R2}=
\begin{bmatrix}
[s]_1 \\
s^N[s]_1 \\
1+[s]_1\\
[s]_1
\end{bmatrix},
\label{6036b} \\
&\mbox{\boldmath$\varphi$}_{R3}=
\begin{bmatrix}
[s]_1 \\
s^N[s]_1 \\
[s]_1\\
1+[s]_1
\end{bmatrix},
\label{6036c} \\
&\mbox{\boldmath$\varphi$}_{R4}=
\begin{bmatrix}
s^{-N}(\sigma+[s]_1) \\
N(1+[s]_1) \\
[s]_1\\
[s]_1
\end{bmatrix}.
\label{6036d} 
\end{align}
\end{subequations}
Here \eqref{6036a}\eqref{6036b}\eqref{6036c} are not obvious, since \eqref{6034}\eqref{6035}
tel us 
\begin{subequations}
\begin{align}
&\mbox{\boldmath$\varphi$}_{R1}=
\begin{bmatrix}
1+p_{11}+\sigma s^{-N}p_{41}\\
Np_{41}\\
p_{21}\\
p_{31}
\end{bmatrix}, \\
&\mbox{\boldmath$\varphi$}_{R2}=
\begin{bmatrix}
\sigma s^{-N}p_{42}\\
N p_{42}\\
1+p_{22}\\
p_{32}
\end{bmatrix}, \\
&\mbox{\boldmath$\varphi$}_{R3}=
\begin{bmatrix}
\sigma s^{-N}p_{43}\\
N p_{43}\\
p_{23}\\
1+p_{33}
\end{bmatrix}, \\
&\mbox{\boldmath$\varphi$}_{R4}=
\begin{bmatrix}
s^{-N}(\sigma+p_{44})\\
N(1+p_{44})\\
p_{24}\\
p_{34}
\end{bmatrix},
\end{align}
\end{subequations}
where
$$\sum_{m\geq 1}s^mP_m^R=(p_{ij}),\qquad
p_{ij}=p_{ij}(s)=[s]_1.$$
Hence \eqref{6036a}\eqref{6036b}\eqref{6036c} claim that
$p_{41}, p_{42}, p_{43} = s^N[s]_1$. This can be proved by using the equation.,
say,
$$ s\frac{dY_1}{ds}=\tilde{a}_{11}Y_1 + 
\tilde{a}_{12}Y_2 + \tilde{a}_{13}Y_3 + \tilde{a}_{14}Y_4 $$
which holds for
$$ \vec{Y}=
\begin{bmatrix}
1+p_{11}+\sigma s^{-N}p_{41} \\
Ns^{-N}p_{41} \\
p_{21} \\
p_{31}
\end{bmatrix}
$$
which corresponds to $\mbox{\boldmath$\varphi$}_{R1}$ and so on. In fact, suppose that $p_{41}=s^{\mu}(C+[s]_1)$ with $C\not=0, 1\leq \mu <N$. Then we would have
$$s\frac{dY_1}{ds}\sim \sigma C(\mu-N)s^{\mu-N} $$
and
$$
\tilde{a}_{11}Y_1 + 
\tilde{a}_{12}Y_2 + \tilde{a}_{13}Y_3 + \tilde{a}_{14}Y_4
\sim \sigma C N s^{\mu-N}, $$
a contradiction.
Suppose that $p_{41}=s^N(C+[s]_1)$ with $C\not=0$. Then we would have
$$s\frac{dY_1}{ds}=[s]_1 $$
and
$$
\tilde{a}_{11}Y_1 + 
\tilde{a}_{12}Y_2 + \tilde{a}_{13}Y_3 + \tilde{a}_{14}Y_4 =\sigma C N + [s]_1,$$
a contradiction. Hence we see $p_{41}=s^{N+1}[s]_0$. The same proof can be done for
$p_{42}, p_{43}$. 
Using the equation 
$$s\frac{dY_3}{ds}=-\frac{R}{R-s^D}s^D(-Y_3+Y_4), $$
we can claim that $p_{21}=p_{31}=0$. In fact, otherwise, $Y_3=p_{21}=C_1s^{\mu}(1+[s]_1),
Y_4=p_{31}=C_2s^{\mu}(1+[s]_1)$, with $\mu \geq 1, (C_1)^2+(C_2)^2\not=0 $; then
$$\mu C_1s^{\mu}(1+[s]_1)=
-\frac{D}{R-s^D}s^{D+\mu}(-C_1+C_2+[s]_1);$$
since $D \geq 1$, this requires $-C_1+C_2=0$; if $C_1\not=0$, then $\mu=D+\mu+k$ with $\exists k \geq 1$, a contradiction; therefore $C_1=C_2=0$, a contradiction. 
\\

If $\vec{y}$ satisfies the boundary condition, then we have
\begin{equation}
\vec{y}=C_{R1}\mbox{\boldmath$\varphi$}_{R1}+C_{R2}(\mbox{\boldmath$\varphi$}_{R2}-(l+1)\mbox{\boldmath$\varphi$}_{R3}).
\end{equation}\\

Summing up, $\vec{y}$ should be
\begin{equation}\vec{y}=C_{01}\mbox{\boldmath$\varphi$}_{O1}+C_{02}\mbox{\boldmath$\varphi$}_{O2}=
C_{R1}\mbox{\boldmath$\varphi$}_{R1}+C_{R2}
(\mbox{\boldmath$\varphi$}_{R2}-(l+1)\mbox{\boldmath$\varphi$}_{R3}) \label{6.30}
\end{equation}
in order that $\vec{y}$ gives $\vec{V}\in \mathsf{D}(\vec{L}_l)$. Conversely, if so, then the corresponding $\vec{V} $ belongs to $\mathfrak{W}_l$
and we see that $V^r= y_1$ is bounded, therefore, thanks to
Proposition \ref{Prop.2}, $\vec{V}$ belongs to $\overset{\circ}{\mathfrak{W}}_l$ 
so to $\mathsf{D}(\vec{L}_l)$ as an eigenfunction.\\

The condition \eqref{6.30} reads
$$D(r,\lambda):=\mathrm{det}(\mbox{\boldmath$\varphi$}_{O1}(r,\lambda),
\mbox{\boldmath$\varphi$}_{O2}(r,\lambda),
\mbox{\boldmath$\varphi$}_{R1}(r,\lambda),
\mbox{\boldmath$\varphi$}_{R2}
(r,\lambda)-(l+1)\mbox{\boldmath$\varphi$}_{R3}(r,\lambda)) =0. $$
But the condition $D(r,\lambda)=0$ is independent of $r$. So, fixing $r_0 \in ]0,R[$, we can consider $f(\lambda)=D(r_0,\lambda)$, which is a holomorphic function of
$\lambda \in \mathbb{C}\setminus \{0\}$. Now we have that 
$f(\lambda)=0$ if and only if $\lambda$ is an eigenvalue of \eqref{6.1}. Since $f(\lambda)\not=0$ for 
$ \lambda \notin [ -C, +\infty[ $, which
$\lambda$ belongs to the resolvent set of $\vec{L}_l$, $f$ is not an identical $0$. Therefore the zeros of $f$  cannot accumulate to a value in $\mathbb{C}\setminus \{0\}$.  Thus we can claim

\begin{Theorem}
Suppose the Assumptions \ref{Ass.G2}, \ref{Ass.G3}, \ref{Ass.nuQ}. If there exist eigenvalues to the  eigenvalue problem
\eqref{6.3}, then they are at most countably many eigenvalues located on the real axis, and cannot accumulate to a value $\not=0$.
\end{Theorem}

\subsection{Inhomogeneous problem}

Now let us consider the inhomogeneous equation \eqref{S02}, supposing that 
$\lambda\not= 0$ and that $\lambda$ is not an eigenvalue.
In other words, considering
\begin{equation}
\mbox{\boldmath$\Phi$}:=[\mbox{\boldmath$\varphi$}_{O1} \quad \mbox{\boldmath$\varphi$}_{02} \quad \mbox{\boldmath$\varphi$}_{R1} \quad 
\mbox{\boldmath$\varphi$}_{R2}-(l+1)\mbox{\boldmath$\varphi$}_{R3} ],
\end{equation}
we suppose that 
\begin{equation}
\mathrm{det}\mbox{\boldmath$\Phi$} \not=0.
\end{equation}
Let $\vec{f}$ be given in $\mathfrak{W}_l$ and suppose that $\|\vec{f}\|_{\mathfrak{W}_l} \leq 1$. Thus the given $\vec{h}=(h_1, h_2, 0, 0)^{\top}$ is supposed to satisfy
$\|h_1\|_{L^2(\rho r^2dr)} \leq C,
\|h_2\|_{L^2(r^2dr /\rho)} \leq C$. Here and hereafter $C$ stands for various constant independent of $\vec{f}$ such that $\|\vec{f}\|_{\mathfrak{W}_l} \leq 1$
and denote by $O(1)$ a quantity $Q$ such that $|Q|\leq  C$, $C$ being  independent of $\vec{f}$ . Later we shall  use the following

\begin{Proposition}\label{Prop.PartInt}
For $(f^r, f^h)^{\top}\in \mathfrak{W}_l$, it holds that
\begin{equation}
\rho(f^r-lf^h)=-\frac{1}{l+1}rf^{[\delta\rho]}-
\frac{1}{l+1}r^l\frac{d}{dr}(r^{-l+1}\rho f^r),
\end{equation}
where
\begin{equation}
f^{[\delta\rho]}:=-\frac{1}{r^2}\frac{d}{dr}(r^2\rho f^r)+
\frac{l(l+1)}{r}\rho f^h
\quad \in\quad L^2(\frac{c^2}{\rho}r^2dr),
\end{equation}
while
$$\|\vec{f}\|_{\mathfrak{W}_l}^2=
\|\vec{f}\|_{\mathfrak{X}_l}^2+
\|f^{[\delta\rho]}\|_{L^2(\frac{c^2}{\rho}r^2dr)}^2.
$$
\end{Proposition}

The solution of 
\eqref{S02} should be of the form
\begin{equation}
\vec{y}(r)=\mbox{\boldmath$\Phi$}(r)\int_{R/2}^r\mbox{\boldmath$\Phi$}(\acute{r})^{-1}
\vec{h}(\acute{r})\frac{d\acute{r}}{\acute{r}}+\mbox{\boldmath$\Phi$}(r)\mathbf{c},
\end{equation}
where $\mathbf{c}$ is a suitable constant vector, which should be determined 
from $\vec{h}$ of $\vec{f}$
so that the corresponding $\vec{V}$ to $\vec{y}$ belong to 
$\mathsf{D}(\vec{L}_l)$.

In \cite[p.365]{Eisenfeld} J. Eisenfeld claims that the constant vector
$\mathbf{c}=(c_1,c_2, c_3, c_4)^{\top}$ should be chosen as
\begin{align}
& c_1=-\int_{R/2}^Rk_1(r)\frac{dr}{r},\qquad 
c_2=-\int_{R/2}^Rk_2(r)\frac{dr}{r}, \nonumber \\
&c_3=\int_0^{R/2}k_3(r)\frac{dr}{r},\qquad
c_4=\int_0^{R/2}k_4(r)\frac{dr}{r},
\end{align}
where
$k_j=k_i(r), i=1,2,3,4,$ is the $i$-th component of the vector
$\mbox{\boldmath$\Phi$}(r)^{-1}\vec{h}(r)$, that is,
\begin{equation}
\vec{k}=
\begin{bmatrix}
k_1 \\
k_2 \\
k_3 \\
k_4
\end{bmatrix}=
\mbox{\boldmath$\Phi$}(r)^{-1}\vec{h}(r)
\end{equation}
In other words, it is claimed that the solution $\vec{y}$ should be given as
\begin{align}
&\vec{y}(r)=\Big(-\int_r^Rk_1(\acute{r})\frac{d\acute{r}}{\acute{r}}\Big)\mbox{\boldmath$\varphi$}_{O1}(r) +
\Big(-\int_r^Rk_2(\acute{r})\frac{d\acute{r}}{\acute{r}}\Big)\mbox{\boldmath$\varphi$}_{O2}(r) + \nonumber \\
&+\Big(\int_0^rk_3(\acute{r})\frac{d\acute{r}}{\acute{r}}\Big)\mbox{\boldmath$\varphi$}_{R1}(r) +
\Big(\int_0^rk_4(\acute{r})\frac{d\acute{r}}{\acute{r}}\Big)(\mbox{\boldmath$\varphi$}_{R2}(r)-(l+1)\mbox{\boldmath$\varphi$}_{R3}(r)). \label{6.48}
\end{align}
However there we can find no persuasive argument on why the constants $c_j$
should be chosen so, or, why
the definite integrals $$\int_{R/2}^Rk_1(r)\frac{dr}{r},\quad\int_{R/2}^Rk_2(r)\frac{dr}{r},\quad\int_0^{R/2}k_3(r)\frac{dr}{r}, \quad
\int_0^{R/2}k_4(r)\frac{dr}{r} $$ are well-determined as finite numbers.
 In order to determine $\mathbf{c}$ so, we need take into account the asymptotic behaviors of $\mbox{\boldmath$\Phi$}(r)$
and $\mbox{\boldmath$\Phi$}(r)^{-1}$ as $r \rightarrow +0$ and as $r\rightarrow R-0$. \\

We have two fundamental matrices $ \mbox{\boldmath$\Phi$}_O(r;\lambda)$ and $\mbox{\boldmath$\Phi$}_R(r;\lambda)$, say,
\begin{align*}
&\mbox{\boldmath$\Phi$}_O=[\mbox{\boldmath$\varphi$}_{O1}\quad \mbox{\boldmath$\varphi$}_{O2}\quad \mbox{\boldmath$\varphi$}_{O3}\quad \mbox{\boldmath$\varphi$}_{O4}], \\
&\mbox{\boldmath$\Phi$}_R=[\mbox{\boldmath$\varphi$}_{R1}\quad \mbox{\boldmath$\varphi$}_{R2}\quad \mbox{\boldmath$\varphi$}_{R3}\quad \mbox{\boldmath$\varphi$}_{R4}].
\end{align*}
Then there exists a non-singular constant matrix $\mathcal{C}(\lambda)$ such that
$$\mbox{\boldmath$\Phi$}_R(r;\lambda)=\mbox{\boldmath$\Phi$}_O(r;\lambda)\mathcal{C}(\lambda). $$
See \cite[Chapter 1, Theorem 7.3]{CoddingtonL}. So we use

\begin{Definition}
We denote
\begin{align}
&\mbox{\boldmath$\Phi$}_O=(\mbox{\boldmath$\varphi$}_{O1}, \mbox{\boldmath$\varphi$}_{O2}, \mbox{\boldmath$\varphi$}_{O3}, \mbox{\boldmath$\varphi$}_{O4}), \\
&\mbox{\boldmath$\Phi$}_R=(\mbox{\boldmath$\varphi$}_{R1}, \mbox{\boldmath$\varphi$}_{R2}, \mbox{\boldmath$\varphi$}_{R3}, \mbox{\boldmath$\varphi$}_{R4}). \\
&\mbox{\boldmath$\Phi$}_R(r;\lambda)=\mbox{\boldmath$\Phi$}_O(r;\lambda)\mathcal{C}(\lambda), \\
&\mathcal{C}(\lambda)=(c_{ij}(\lambda))_{i,j=1,2,3,4}, \\
\mbox{that is,} & \nonumber \\
&\mbox{\boldmath$\varphi$}_{Rj}=\sum_{i=1}^4c_{ij}(\lambda)\mbox{\boldmath$\varphi$}_{Oi}, \qquad j=1,2,3,4, \\
&\mathcal{C}(\lambda)^{-1}=(\gamma_{ji}(\lambda))_{j,i=1,2,3,4}, \\
\mbox{that is,}  &\nonumber \\
&\mbox{\boldmath$\varphi$}_{Oi}=\sum_{j=1}^4\gamma_{ji}(\lambda)\mbox{\boldmath$\varphi$}_{Rj},\qquad i=1,2,3,4.
\end{align}
\end{Definition}

Then we have
\begin{align}
\mbox{\boldmath$\Phi$}&=\mbox{\boldmath$\Phi$}_O
\begin{bmatrix}
1 & 0 & c_{11} & c_{12}-(l+1)c_{13} \\
0 & 1 & c_{21} & c_{22}-(l+1)c_{23} \\
0 & 0 & c_{31} & c_{32}-(l+1)c_{33} \\
0 & 0 & c_{41} & c_{42}-(l+1)c_{43}
\end{bmatrix}  \label{6.56}\\
&=
\mbox{\boldmath$\Phi$}_R
\begin{bmatrix}
\gamma_{11} & \gamma_{12} & 1 & 0 \\
\gamma_{21} & \gamma_{22} & 0 & 1 \\
\gamma_{31} & \gamma_{32} & 0 & -(l+1) \\
\gamma_{41} & \gamma_{42} & 0 & 0 
\end{bmatrix}. \label{6.57}
\end{align}

Here we suppose that  $\lambda$ is not an eigenvalue. Then $\mathrm{det}\mbox{\boldmath$\Phi$} \not=0$ and it
means
\begin{align}
&\mathrm{det}
\begin{bmatrix}
c_{31} & c_{32}-(l+1)c_{33} \\
c_{41} & c_{42}-(l+1)c_{43}
\end{bmatrix} 
\not=0, \\
\mbox{and} & \nonumber \\
&\mathrm{det}
\begin{bmatrix}
\gamma_{31}+(l+1)\gamma_{21} & \gamma_{32}+(l+1)\gamma_{22} \\
\gamma_{41} & \gamma_{42}
\end{bmatrix}
\not= 0.
\end{align}\\

Now, by dint of \eqref{60020a} - \eqref{60020d}, we have
\begin{equation}
\mbox{\boldmath$\Phi$}_O=
\begin{bmatrix}
(L+\mathsf{o})r^{l-1} & \mathsf{o}r^{l-1} & (L+\mathsf{o})r^{-(l+2)} & \mathsf{o}r^{-l+2} \\
(\alpha +\mathsf{o})r^{l-1} & (-\rho_{\mathsf{O}}+\mathsf{o})r^{l-1} & (\beta+\mathsf{o})r^{-(l+2)} & (-\rho_{\mathsf{O}}+\mathsf{o})r^{-(l+2)} \\
\mathsf{o}r^{l-1} & (1+\mathsf{o})r^{l-1} & \mathsf{o}r^{-(l+2)} & (1+\mathsf{o})r^{-(l+2)} \\
\mathsf{o} r^{l-1} & (l+\mathsf{o})r^{l-1} & \mathsf{o} r^{-(l+2)} & (-(l+1)+\mathsf{o})r^{-(l+2)} 
\end{bmatrix}. \label{6.60}
\end{equation}
This implies
\begin{align}
&\mbox{\boldmath$\Phi$}_O^{-1}= \nonumber \\
&\begin{bmatrix}
\Big(\frac{\beta}{\beta-\alpha}\frac{1}{L}+\mathsf{o}\Big)r^{-(l-1)} & \Big(-\frac{1}{\beta-\alpha}+\mathsf{o}\Big)r^{-(l-1)} & \Big(-\frac{\rho_{\mathsf{O}}}{\beta-\alpha}+\mathsf{o}\Big)r^{-(l-1)} & \mathsf{o} r^{-(l-1)} \\
\mathsf{o} r^{-(l-1)} & \mathsf{o}r^{-(l-1)} & \Big(\frac{l+1}{2l+1}+\mathsf{o}\Big)r^{-(l-1)} & \Big(\frac{1}{2l+1}+\mathsf{o}\Big)r^{-(l-1)} \\
\Big(-\frac{\alpha}{\beta-\alpha}\frac{1}{L}+\mathsf{o}\Big)r^{l+2} & \Big(\frac{1}{\beta-\alpha}+\mathsf{o}\Big)r^{l+2} & \Big(\frac{\rho_{\mathsf{O}}}{\beta-\alpha}+ \mathsf{o}\Big) r^{l+2} & \mathsf{o}r^{l+2} \\
\mathsf{o}r^{l+2} & \mathsf{o}r^{l+2} & \Big(\frac{l}{2l+1}+\mathsf{o}\Big)r^{l+2} & \Big(-\frac{1}{2l+1}+\mathsf{o}\Big)r^{l+2}
\end{bmatrix}
\label{6.61}.
\end{align}
In \eqref{6.60}, \eqref{6.61}, the symbol $\mathsf{o}$ stands for $O(r)$. 

By dint of \eqref{6036a} - \eqref{6036d}, we have
\begin{equation}
\mbox{\boldmath$\Phi$}_R=
\begin{bmatrix}
1+\mathsf{o} & \mathsf{o} &
 \mathsf{o} & s^{-N}(\sigma +\mathsf{o}) \\
O(s^{N+1}) & O(s^{N+1}) & O(s^{N+1}) & N(1+\mathsf{o}) \\
0 & 1+\mathsf{o} & \mathsf{o} & \mathsf{o} \\
0 & \mathsf{o} & 1+\mathsf{o} & \mathsf{o}
\end{bmatrix} \label{6.62}
\end{equation}
This implies
\begin{equation}
\mbox{\boldmath$\Phi$}_R^{-1}=
\begin{bmatrix}
1+\mathsf{o} & -s^{-N}\Big(\frac{\sigma}{N}+\mathsf{o}\Big) & \mathsf{o} & \mathsf{o} \\
\mathsf{o}O(s^{N+1}) & \mathsf{o} & 1+\mathsf{o} & \mathsf{o} \\
\mathsf{o}O(s^{N+1}) & \mathsf{o} & \mathsf{o} & 1+\mathsf{o} \\
O(s^{N+1}) & \frac{1}{N}+\mathsf{o} & \mathsf{o} & \mathsf{o}
\end{bmatrix} \label{6.63}
\end{equation}
In \eqref{6.62},\eqref{6.63} the symbol $\mathsf{o}$ stands for various $O(s)$. \\

Let us write \eqref{6.56} as
\begin{equation}
\mbox{\boldmath$\Phi$}=\mbox{\boldmath$\Phi$}_O
\begin{bmatrix}
I & \mathcal{A}_O \\
O & \mathcal{B}_O
\end{bmatrix},
\end{equation}
where 
\begin{equation}
\mathcal{A}_O=
\begin{bmatrix}
c_{11} & c_{12}-(l+1)c_{13} \\
c_{21} & c_{22}-(l+1)c_{23}
\end{bmatrix}, \quad
\mathcal{B}_O=
\begin{bmatrix}
c_{31} & c_{32}-(l+1)c_{33} \\
c_{41} & c_{42}-(l+1)c_{43}
\end{bmatrix}.
\end{equation}
Since we are supposing that $\lambda$ is ont an eigenvalue, $\mathcal{B}_O$ is invertible.

Let us write \eqref{6.57} as 
\begin{align}
&\mbox{\boldmath$\Phi$}= \nonumber \\
&=\mbox{\boldmath$\Phi$}_R
\begin{bmatrix}
I & O \\
-(l+1)J & I
\end{bmatrix}
\begin{bmatrix}
I & \mathcal{A}_R \\
O & \mathcal{B}_R
\end{bmatrix}
\begin{bmatrix}
O & I \\
I & O
\end{bmatrix}, \label{6.66}
\end{align}
where 
\begin{align}
&J = 
\begin{bmatrix}
0 & 1 \\
0 & 0
\end{bmatrix}, \\
&\mathcal{A}_R=
\begin{bmatrix}
\gamma_{11} & \gamma _{12} \\
\gamma _{21} & \gamma_{22}
\end{bmatrix}, \quad
\mathcal{B}_R=
\begin{bmatrix}
\gamma_{31}+(l+1)\gamma_{21} & \gamma_{32}+ (l+1)\gamma_{22} \\
\gamma_{41} & \gamma_{42}
\end{bmatrix}.
\end{align}
Since $ \lambda$ is not an eigenvalue, $\mathcal{B}_R$ is invertible. 

Then we have
\begin{align}
\mbox{\boldmath$\Phi$}^{-1}&=
\begin{bmatrix}
I & -\mathcal{A}_O\mathcal{B}_O^{-1} \label{6.69}\\
O & \mathcal{B}_O^{-1}
\end{bmatrix}
\mbox{\boldmath$\Phi$}_O^{-1} \\
&=
\begin{bmatrix}
(l+1)\mathcal{B}_R^{-1}J & \mathcal{B}_R^{-1} \\
I & -\mathcal{A}_R\mathcal{B}_R^{-1}
\end{bmatrix} \label{6.70}
\mbox{\boldmath$\Phi$}_R^{-1}.
\end{align}\\

Let us denote
\begin{equation}
\vec{\mathfrak{h}}_O=
\begin{bmatrix}
\mathfrak{h}_O^1 \\
\\
\mathfrak{h}_O^2 \\
\\
\mathfrak{h}_O^3 \\
\\
\mathfrak{h}_O^4
\end{bmatrix}
:=\mbox{\boldmath$\Phi$}_O^{-1}\vec{h}.
\end{equation}
Then \eqref{6.61} reads
\begin{align}
&\mathfrak{h}_O^1 =\Big[\Big(\frac{\beta}{\beta-\alpha}\frac{1}{L}+O(r)\Big)h_1+
\Big(-\frac{1}{\beta-\alpha}+O(r)\Big)h_2\Big]r^{-(l-1)} \nonumber \\
&\mathfrak{h}_O^2=[O(r)h_1+O(r)h_2] r^{-(l-1)} \nonumber \\
&\mathfrak{h}_O^3=
\Big[\Big(-\frac{\alpha}{\beta-\alpha}\frac{1}{L}+O(r)\Big)h_1+
\Big(\frac{1}{\beta-\alpha}+O(r)\Big)h_2\Big]r^{l+2} \nonumber  \\
&\mathfrak{h}_O^4=[O(r)h_1+O(r)h_2]r^{l+2}. \label{6.72}
\end{align}

Let us denote
\begin{equation}
\vec{\mathfrak{h}}_R=
\begin{bmatrix}
\mathfrak{h}_R^1 \\
\\
\mathfrak{h}_R^2 \\
\\
\mathfrak{h}_R^3 \\
\\
\mathfrak{h}_R^4
\end{bmatrix}
:=\mbox{\boldmath$\Phi$}_R^{-1}\vec{h}.
\end{equation}
Then \eqref{6.63} reads
\begin{align}
&\mathfrak{h}_R^1=(1+O(s))h_1-s^{-N}\Big(\frac{\sigma}{N}+O(s)\Big)h_2, \nonumber \\
&\mathfrak{h}_R^2=O(s^{N+2})h_1+O(s)h_2 \nonumber \\
&\mathfrak{h}_R^3=O(s^{N+2})h_1+O(s)h_2 \nonumber \\
&\mathfrak{h}_R^4=O(s^{N+1})h_1+\Big(\frac{1}{N}+O(s)\Big)h_2. \label{6.74}
\end{align}

Let us look at
\begin{equation}
\vec{k}=\mbox{\boldmath$\Phi$}^{-1}\vec{h}=
\begin{bmatrix}
(l+1)\mathcal{B}_R^{-1}J & \mathcal{B}_R^{-1} \\
I & -\mathcal{A}_R\mathcal{B}_R^{-1}
\end{bmatrix}
\vec{\mathfrak{h}}_R.
\end{equation}
Keeping in mind that 
$$J
\begin{bmatrix}
\mathfrak{h}_R^1 \\
\\
\mathfrak{h}_R^2
\end{bmatrix}
=
\begin{bmatrix}
\mathfrak{h}_R^2\\
\\
0
\end{bmatrix},
$$
which does not contain $\mathfrak{h}_R^1$, we see that, for $i=1,2$, \eqref{6.74}
implies
$$|k_i|\leq C (s^{N+1}|h_1|+|h_2|) $$
so that
\begin{equation}
K_i(r):=-\int_r^Rk_i(\acute{r})\frac{d\acute{r}}{\acute{r}},\qquad i=1,2,
\end{equation}
are well defined and enjoy the estimate
\begin{align}
|K_i(r)|&\leq C\int_0^s|h_1(\acute{r})|\acute{s}^{N+D}d\acute{s}+
C\int_0^s|h_2(\acute{r})|\acute{s}^{D-1}d\acute{s} \nonumber \\
&\leq C\sqrt{\int_0^s \acute{s}^{N+D+1}d\acute{s}}
\sqrt{\int |h_1|^2\acute{s}^{N+D-1}d\acute{s}}  \nonumber \\
&+C\sqrt{\int_0^s \acute{s}^{D+N-1}d\acute{s}}
\sqrt{\int |h_2|^2\acute{s}^{-N+D-1}d\acute{s}} \nonumber \\
& \leq C' s^{\frac{N+D}{2}}.
\end{align}
Here $C'$ is independent of $\vec{h}$ which comes from $\vec{f}$ such that 
$\|\vec{f}\|_{\mathfrak{X}_l} \leq 1$.
In this sense 
\begin{equation}
c_i=-\int_{R/2}^Rk_i(r)\frac{dr}{r}=K_i\Big(\frac{R}{2}\Big),\qquad i=1,2
\end{equation}
are well defined as finite numbers.

On the other hand, look at
$$
\vec{k}=
\begin{bmatrix}
I & -\mathcal{A}_O\mathcal{B}_O^{-1} \\
O & \mathcal{B}^{-1}
\end{bmatrix}
\begin{bmatrix}
\mathfrak{h}_O^1 \\
\\
\mathfrak{h}_O^2 \\
\\
\mathfrak{h}_O^3 \\
\\
\mathfrak{h}_O^4 
\end{bmatrix}.
$$
By \eqref{6.72}, we see that, for $i=3,4$,
$$
|k_i|\leq C(|h_1|+|h_2|)r^{l+2}.$$
Therefore
\begin{equation}
K_i(r):=\int_0^rk_i(\acute{r})\frac{d\acute{r}}{\acute{r}},\qquad i=3,4,
\end{equation}
are well defined and enjoy the estimates
\begin{align}
|K_i(r)|&\leq C\sqrt{\int_0^r\acute{r}^{2l}d\acute{r}}
\sqrt{\int |h_1|^2\acute{r}^2d\acute{r}} + \nonumber \\
&+ C\sqrt{\int_0^r\acute{r}^{2l}d\acute{r}}\sqrt{\int|h_2|^2\acute{r}^2d\acute{r}} \nonumber  \\
&\leq C' r^{\frac{2l+1}{2}}.
\end{align}
In this sense 
\begin{equation}
c_i=\int_0^{R/2}k_i(r)\frac{dr}{r}=K_i\Big(\frac{R}{2}\Big),\qquad i=3,4
\end{equation}
are well defined as finite numbers.\\

Let us examine the behaviors of the solution $\vec{y}$ given by \eqref{6.48}:
\begin{align}
&\vec{y}(r)=\Big(-\int_r^Rk_1(\acute{r})\frac{d\acute{r}}{\acute{r}}\Big)\mbox{\boldmath$\varphi$}_{O1}(r) +
\Big(-\int_r^Rk_2(\acute{r})\frac{d\acute{r}}{\acute{r}}\Big)\mbox{\boldmath$\varphi$}_{O2}(r) + \nonumber \\
&+\Big(\int_0^rk_3(\acute{r})\frac{d\acute{r}}{\acute{r}}\Big)\mbox{\boldmath$\varphi$}_{R1}(r) +
\Big(\int_0^rk_4(\acute{r})\frac{d\acute{r}}{\acute{r}}\Big)(\mbox{\boldmath$\varphi$}_{R2}(r)-(l+1)\mbox{\boldmath$\varphi$}_{R3}(r)). 
\end{align}

Let us consider the behavior of $\vec{y}$ as $r \rightarrow +0$. 

As for
$$
\begin{bmatrix}
k_1 \\
\\
k_2
\end{bmatrix}
=\begin{bmatrix}
\mathfrak{h}_O^1 \\
\\
\mathfrak{h}_O^2
\end{bmatrix}
-\mathcal{A}_O\mathcal{B}_O^{-1}
\begin{bmatrix}
\mathfrak{h}_O^3 \\
\\
\mathfrak{h}_O^4
\end{bmatrix},
$$
we need a sharp estimate of
$$
\mathfrak{h}_O^1=\Big[\Big(\frac{\beta}{\beta-\alpha}\frac{1}{L}+O(r)\Big)h_1+
\Big(-\frac{1}{\beta-\alpha}+O(r)\Big)h_2\Big]r^{-(l-1)}.
$$
Keeping in mind Proposition \ref{Prop.PartInt}, we observe
\begin{align*}
&\frac{\beta}{\beta-\alpha}\frac{1}{L}h_1 -\frac{1}{\beta-\alpha}h_2  = \\
&=\frac{1}{\beta-\alpha}\Big[\Big(-\rho_{\mathsf{O}}Lg_O+gL\frac{g}{r}\Big)f^h-
\rho f^r +\rho_{\mathsf{O}}lf^h\Big] \\
&=\frac{1}{\beta-\alpha}\Big[
[r^2]_1f^h
-\rho f^r +\rho lf^h\Big] \\
&=\frac{1}{\beta-\alpha}\Big[
[r^2]_1f^h
+\frac{1}{l+1}rf^{[\delta\rho]}+\frac{1}{l+1}r^l\frac{d}{dr}(r^{-l+1}\rho f^r)\Big].
\end{align*}

Then we have

\begin{align*}
\mathfrak{h}_O^1&=O(1)(|h_1|+|h_2|+ |f^{[\delta\rho]}|)r^{-l+2}
+\frac{1}{(\beta-\alpha)(l+1)}r\frac{d}{dr}(r^{-l+1}\rho f^r), \\
k_1(r)&=O(1)(|h_1|+|h_2|+ |f^{[\delta\rho]}|)r^{-l+2}
+\frac{1}{(\beta-\alpha)(l+1)}r\frac{d}{dr}(r^{-l+1}\rho f^r), \\
K_1(r)& =
O\Big(\int_r^R(|h_1|+|h_2|+|f^{[\delta\rho]}|)\acute{r}^{-l+1}d\acute{r}\Big)+ \\
&-\frac{1}{(\beta-\alpha)(l+1)}\int_r^R\frac{d}{d\acute{r}}(\acute{r}^{-l+1}\rho f^r(\acute{r}))d\acute{r} \\
&=O(r^{\frac{-2l+1}{2}})+
\frac{1}{(\beta-\alpha)(l+1)}(r^{-l+1}\rho f^r +O(1)), \\
K_1\mbox{\boldmath$\varphi$}_{O1,i}&=O(r^{-\frac{1}{2}})+O(1)
f^r, \\
\mbox{and} & \\
&\|K_1\mbox{\boldmath$\varphi$}_{O1,i}\|_{L^2(r^2dr)}  \leq C
\end{align*}
for $i=1,2$.

On the other hand,
$$|k_2(r)|\leq Cr^2(|h_1|+|h_2|)r^{-(l-1)} $$
implies
$$\|K_2\mbox{\boldmath$\varphi$}_{O1,i}\|_{L^2(r^2dr)} \leq C $$
for $i=1,2$.

As for 
$$
\begin{bmatrix}
k_3 \\
\\
k_4
\end{bmatrix}
=\mathcal{B}_O^{-1}
\begin{bmatrix}
\mathfrak{h}_O^3 \\
\\
\mathfrak{h}_O^4
\end{bmatrix}
, $$
 we have
\begin{align*}
|k_i(r)|&\leq C\Big[
(|h_1|+|h_2|)
\Big] r^{l+2}, \\
|K_i(r)|&=\Big|\int_0^R
k_i(\acute{r})\frac{d\acute{r}}{\acute{r}}\Big| \\
&\leq C r^{\frac{2l+1}{2}}
\end{align*}
for $i=3,4$.
Since $|\mbox{\boldmath$\varphi$}_{Ri}|\leq Cr^{-(l+1)}$, we see
$$
\|K_i \mbox{\boldmath$\varphi$}_{Ri}\|_{L^2(r^2dr)}\leq C
$$
for $i=3,4$. 

Summing up, 
$ \|y_1\|_{L^2(r^2dr)}$ and $\|y_2\|_{L^2(r^2dr)}$ are estimated by
$ C$. \\

Next we consider the behavior of $\vec{y}$ defined by \eqref{6.48} as $r \rightarrow R-0$.

We know that, for $j=1,2$,
$$ K_j(r)=-\int_r^Rk_j(\acute{r})\frac{d\acute{r}}{\acute{r}} $$
enjoys
$$|K_j(r)|\leq C s^{\frac{N+D}{2}}. $$
Since $|\mbox{\boldmath$\varphi$}_{Oj, 1}|\leq C s^{-N}$, we have
$$
|K_j \mbox{\boldmath$\varphi$}_{Oj, 1}|\leq C s^{\frac{-N+D}{2}}\quad\in\quad L^2(s^{N+D-1}ds), \quad j=1,2.
$$
Since $|\mbox{\boldmath$\varphi$}_{Oj,2}|\leq C$, we have
$$|K_j\mbox{\boldmath$\varphi$}_{Oj,2}|\leq C s^{\frac{N+D}{2}} \quad \in \quad L^2(s^{-N-1}ds),\quad j=1,2.
$$

As for $j=3,4$, we look at
$$
\begin{bmatrix}
k_3 \\
\\
k_4
\end{bmatrix}
=
\begin{bmatrix}
\mathfrak{h}_R^1 \\
\\
\mathfrak{h}_R^2
\end{bmatrix}
-\mathcal{A}_R\mathcal{B}_R^{-1}
\begin{bmatrix}
\mathfrak{h}_R^3 \\
\\
\mathfrak{h}_R^4
\end{bmatrix}
. $$
\eqref{6.74} implies
\begin{align*}
|k_3|&\leq C(|h_1|+s^{-N}|h_2|), \\
|k_4|&\leq C(s^{N+1}|h_1|+|h_2|).
\end{align*}
Thereforer
$$K_3(r)=\int_0^rk_3(\acute{r})\frac{d\acute{r}}{\acute{r}} $$
enjoys
$$
|K_3|\leq C s^{\frac{-N+D}{2}} $$
and
$$
K_4(r)=\int_0^rk_4(\acute{r})\frac{d\acute{r}}{\acute{r}} $$
enjoys
$$
|K_4|\leq C.
$$
Therefore
\begin{align*}
|K_3\mbox{\boldmath$\varphi$}_{R1,1}|&\leq C s^{\frac{-N+D}{2}}\quad \in\quad L^2(s^{N+D-1}ds), \\
|K_3\mbox{\boldmath$\varphi$}_{R1, 2}|&\leq C s^{\frac{N+D+2}{2}}\quad
\in \quad L^2(s^{-N-1}ds), \\
|K_4(\mbox{\boldmath$\varphi$}_{R2,1}-(l+1)\mbox{\boldmath$\varphi$}_{R3, 1})|&
\leq Cs \quad \in \quad L^2(s^{N+D-1}ds), \\
|K_4(\mbox{\boldmath$\varphi$}_{R2,2}-(l+1)\mbox{\boldmath$\varphi$}_{R3, 2})|&
\leq Cs^{N+1} \quad \in \quad L^2(s^{-N-1}ds).
\end{align*}

Summing up, we can estimate $\|y_1\|_{L^2(\rho r^2dr)}$ and 
$\|y_2\|_{L^2(r^2dr/\mathsf{c}^2\rho)}$ by $C$. 

Moreover the estimate obtained above 
$$ |y_1|\leq C s^{\frac{-N+D}{2}} \leq C' (R-r)^{\frac{-\nu+1}{2}} $$
enjoys the condition sufficient for $\vec{V}$ to belong to 
$\overset{\circ}{\mathfrak{W}}_l$, thanks to Proposition \ref{Prop.2}\\

Summing up, we can claim the following

\begin{Theorem}
Let $l \geq 1$. 
Suppose the Assumptions \ref{Ass.G2}, \ref{Ass.G3}, \ref{Ass.nuQ}. 
Let $\lambda \not=0$ be not an eigenvalue of \eqref{S01}. Suppose that 
$\vec{f} \in \mathfrak{W}_l$. Then the problem \eqref{S02}
admits a solution $\vec{V} \in \mathsf{D}(\vec{L}_l)$ 
such that $\|\vec{V}\|_{\mathfrak{W}_l}\leq C\|\vec{f}\|_{\mathfrak{W}_l}$. Therefore the spectrum $\sigma(\vec{L}_l)$ of the operator $\vec{L}_l$ as a self-adjoint operator in $\mathfrak{W}_l$ is such that $\sigma(\vec{L}_l)\setminus \{0\}$  
consists of countable many eigenvalues which do not accumulate to a value $\not=0, \not=\infty$.
\end{Theorem}

As a corollary we can claim

\begin{Theorem}
The eigenfunctions of $\vec{L}_l$ form a complete orthogonal system of $\mathfrak{W}_l$.
\end{Theorem}

For a proof of the completeness of eigenfunctions, see \cite[p.905, X.3.4 Theorem]{DunfordS}, which can be applied 
to the unbounded self-adjoint operator $\vec{L}_l$ thanks to \cite[p.177, Chapter III, Theorem 6.15]{Kato}. \\

However, the absence of continuous spectra for the operator $\vec{L}_l$ considered as an operator in  $\mathfrak{X}_l$, which is asserted in \cite{Eisenfeld}, is doubtful, although we have not yet found a counter example.

\vspace{10mm}

{\bf \large Acknowledgment}\\ 

The idea to prove the existence of g-modes was obtained 
during the stay of the author at the Department of Mathematics, National University of Singapore in March 4-11, 2020. The author expresses his sincerely deep  thanks to Professor Shih-Hsien Yu for the invitation and stimulating discussions, and to the Department of Mathematics, National University of Singapore for the hospitality and the financial support. The author expresses his sincerely deep thanks to the anonymous referee who, having read the manuscript carefully, gave helpful comments to ameliorate the presentation. This work is  supported by JSPS KAKENHI Grant Number JP18K03371. \\

{\bf\large Appendix}\\

{\bf 1).}\   We consider the solution of
$$
-\frac{1}{r^2}\frac{d}{dr}r^2\frac{du}{dr}=4\pi\mathsf{G}f^{\rho}(u),\quad
u=u_{\mathsf{O}}+O(r^2)\quad\mbox{as}\quad r\rightarrow +0.
$$

As shown in \cite{TM1984}, $u$ is monotone decreasing and either 
1) $u(r)>0$ for $\forall r \geq 0$ and $u(r) \rightarrow 0$ as $r \rightarrow +\infty$, or 2) there is a finite zero $R$, namely $u(R)=0$, and 
$$ u(r)=\mu\Big(\frac{1}{r}-\frac{1}{R}\Big)\quad \mbox{for}\quad r \geq R$$
with $\displaystyle \mu=-r^2\frac{du}{dr}\Big|_{r=R} >0$. 
We are going to prove that 2) is the case if $\frac{4}{3}<\gamma$ or if $\frac{6}{5}<\gamma$ and $u_{\mathsf{O}}$ is small.

Let $\frac{4}{3}<\gamma<2$. Note that
$$f^{\rho}(u) \sim \Big(\frac{\gamma-1}{\mathsf{A}\gamma}\Big)^{\frac{1}{\gamma-1}}u^{\frac{1}{\gamma-1}}\quad\mbox{as}\quad u \rightarrow +0.$$
We see 
$$\int_0^14\pi\mathsf{G}f^{\rho}(u)u^{-4}du =+\infty $$
for $\displaystyle \frac{1}{\gamma-1}<3$. Thus $u$ cannot remain positive on $[0,+\infty[$ but should have a finite radius, by \cite[Theorem 1]{TM1984}.

Let $\frac{6}{5}<\gamma <2$. By the change of variables $(r,u) \mapsto (\xi, U)$ defined by
$$u=u_{\mathsf{O}}U,\quad r=\frac{1}{\sqrt{4\pi\mathsf{G}}}
\Big(\frac{\mathsf{A}\gamma}{\gamma-1}\Big)^{\frac{1}{2(\gamma-1)}}
u_{\mathsf{O}}^{-\frac{2-\gamma}{2(\gamma-1)}}\xi, $$
we have 
$$
-\frac{1}{\xi^2}\frac{d}{d\xi}\xi^2\frac{dU}{d\xi}=F(U),\quad U=
1+O(\xi^2) \quad\mbox{as}\quad \xi \rightarrow +0.$$
Here 
$$F(U):=(U\vee 0)^{\frac{1}{\gamma-1}}(1+\Lambda_{\rho}(u_{\mathsf{O}}U)).$$
Note that $F(U)$ converges to $(U\vee 0)^{\frac{1}{\gamma-1}}$ uniformly on $U \in ] -\infty,1]$ as $u_{\mathsf{O}} \rightarrow +0$. So, $U(\xi)$ converges to $\theta(\xi)$ uniformly as
$u_{\mathsf{O}} \rightarrow +0$, where $\theta(\xi)$ is the Lane-Emden function of index $\frac{1}{\gamma-1}$, namely, the solution of
$$-\frac{1}{\xi^2}\frac{d}{d\xi}\xi^2\frac{d\theta}{d\xi}=(\theta\vee 0)^{\frac{1}{\gamma-1}},\quad \theta =1+O(\xi^2)\quad\mbox{as}\quad \xi \rightarrow +0.$$
Since $\frac{6}{5}<\gamma$, it is known that $\theta(\xi)$ has the finite zero
$\xi_1(\frac{1}{\gamma-1})$. Therefore $\theta(\xi_0) <0$ for $\xi_0:=2\xi_1(\frac{1}{\gamma-1})$. Since $U(\xi)$ tends to $\theta(\xi)$ uniformly as
 $u_{\mathsf{O}}\rightarrow 0$, we have
$U(\xi_0)<0$ provided that $u_{\mathsf{O}}$ is sufficiently small. Then there is $\Xi$ such that
$0<\Xi <\xi_0$ and $U(\Xi)=0$. This means that $u(R)=0$, where
$$R=
\frac{1}{\sqrt{4\pi\mathsf{G}}}
\Big(\frac{\mathsf{A}\gamma}{\gamma-1}\Big)^{\frac{1}{2(\gamma-1)}}
u_{\mathsf{O}}^{-\frac{2-\gamma}{2(\gamma-1)}}\Xi. $$\\

{\bf 2).}\   Let us consider a function $f$ of the form
$$f(u)=Ku^{\nu}(1+\sum_{k\geq 1}f_ku^k) $$
as $u \rightarrow +0$, where $K$ is a positive constant, $1<\nu <+\infty$, $\sum f_ku^k$
is a convergent power series, while $f(u)>0$ for $u >0$. Suppose $u=u(r), 0<r <R,$
satisfies $u >0, du/dr <0,$
$$\frac{d^2u}{dr^2}+\frac{2}{r}\frac{du}{dr}+f(u)=0$$
on $0<r<R$, 
$u \rightarrow +0$ as $r\rightarrow R-0$, and the limit
$\displaystyle  \lim_{r\rightarrow R-0}\frac{du}{dr} $ exists to be finite and strictly negative.

Then we have the expansion
$$u=C\frac{R-r}{R}
\Big[
1+
\sum_{k_1+k_2+k_3\geq 1}b_{k_1k_2k_3}
\Big(\frac{R-r}{R}\Big)^{k_1}
\Big(C'\Big(\frac{R-r}{R}\Big)^{\nu+1}\Big)^{k_2}
\Big(C\frac{R-r}{R}\Big)^{k_3}\Big],$$
where $C$ is a positive constant, $C'=R^2KC^{\nu-1}$, and
$\sum b_{k_1k_2k_3}X_1^{k_1}X_2^{k_2}X_3^{k_3}$ is a convergent triple power series.\\

Let us sketch the proof. First note that there is a convergent power series 
$$G(u)=\sum_{k\geq 1}G_ku^k$$
such that
$$u\frac{d}{du}\log f(u)=\nu+G(u).$$

Putting
$$x_1=-\Big(\frac{r}{u}\frac{du}{dr}\Big)^{-1},
\quad
x_2=\frac{r^2f(u)}{u}\Big(\frac{r}{u}\frac{du}{dr}\Big)^{-2},
\quad x_3=u, $$
we get the autonomous system
\begin{align*}
&u\frac{dx_1}{du}=(1-x_1+x_2)x_1, \\
&u\frac{dx_2}{du}=(\nu+1-4x_1+2x_2+G(x_3))x_2, \\
&u\frac{dx_3}{du}=x_3.
\end{align*}
Here $(x_1,x_2, x_3)$ is considered as functions of $u >0$, and we have
$(x_1,x_2,x_3) \rightarrow (0,0,0)$ as $u \rightarrow +0$. 

Then there is a transformation of variables $(x_1,x_2,x_3) \leftrightarrow (\xi_1,\xi_2,\xi_3)$ of the form
$$x_j=\xi_j(1+[\xi_1,\xi_2,\xi_3]_1),\quad j=1,2,3,$$
which reduce the system to
$$u\frac{d\xi_1}{du}=\xi_1,\quad
u\frac{d\xi_2}{du}=(\nu+1)\xi_2 \quad u\frac{d\xi_3}{du}=\xi_3.$$
Here and hereafter $[X_1, X_2, X_3]_1$ generally stands for various convergent triple power series
of the form
$\sum_{k_1+k_2+k_3\geq 1}a_{k_1k_2k_3}X_1^{k_1}X_2^{k_2}X_3^{k_3}$.
Take a general solution
$$\xi_1=C_1u,\quad \xi_2=C_2u^{\nu+1},\quad \xi_3=u. $$
putting $C=1/C_1$, we get the desired expansion, since the integration of
$$-\frac{u}{r}\frac{dr}{du}=x^1 $$
gives
\begin{align*}
\log \frac{R}{r}&=-\log \Big(1-\frac{R-r}{R}\Big)=\frac{R-r}{R}\Big(1+
\Big[\frac{R-r}{R}\Big]_1\Big) \\
&=C_1u(1+[C_1u,C_2u^{\nu+1},u]_1).
\end{align*}

\end{document}